\documentclass[letterpaper,11pt,reqno]{amsart} 
\usepackage[margin=1in]{geometry}
\usepackage[T1]{fontenc}
\usepackage{pgfplots,lmodern,microtype}
\usepackage{setspace}
\DeclareMathOperator{\cl}{trace}
\pgfplotsset{compat=1.18}
\usepackage{systeme,mathrsfs,xfrac} 
\usepackage[colorlinks=true,linkcolor=blue,citecolor=blue,urlcolor=blue,pdfencoding=auto,psdextra]{hyperref} 
\usepackage{amsmath,amssymb,amsthm,amsfonts,amsbsy,latexsym,dsfont,color} 
\usepackage[noadjust,sort]{cite}
\usepackage[foot]{amsaddr}
\usepackage{comment}
\usepackage{tikz}
\usetikzlibrary{positioning,arrows.meta,calc,bending,decorations.pathreplacing,backgrounds}

\definecolor{base}{HTML}{243449}
\definecolor{particleA}{HTML}{168C91}
\definecolor{particleB}{HTML}{E38B2C}
\definecolor{particleC}{HTML}{8B61B2}
\definecolor{particleD}{HTML}{CF5575}

\usepackage[textsize=tiny]{todonotes}
\usepackage{regexpatch}
\makeatletter
\xpatchcmd{\@todo}{\setkeys{todonotes}{#1}}{\setkeys{todonotes}{inline,#1}}{}{}
\makeatother

\usepackage{enumitem}

\newenvironment{enumeratea}{\begin{enumerate}[label=\textup{\alph*)},leftmargin=2em,itemsep=1ex]}{\end{enumerate}}
\newenvironment{enumeraten}{\begin{enumerate}[label=\textup{\arabic*)},leftmargin=2em,itemsep=1ex]}{\end{enumerate}}
\newenvironment{enumerateA}{\begin{enumerate}[label=\textup{(\Alph*)},leftmargin=2em,itemsep=1ex]}{\end{enumerate}}

\newtheorem{thm}{Theorem}[section]

\newtheorem{lem}[thm]{Lemma}

\newtheorem{cor}[thm]{Corollary}
\newtheorem{prop}[thm]{Proposition}

\theoremstyle{definition}
\newtheorem{defn}[thm]{Definition}
\newtheorem{rem}[thm]{Remark}
\newtheorem{remark}[thm]{Remark}

\newtheorem{claim}[thm]{Claim}

\renewcommand{\le}{\leqslant}  
\renewcommand{\ge}{\geqslant} 
\renewcommand{\leq}{\leqslant}  
\renewcommand{\geq}{\geqslant}

\newcommand{\ind}{\mathds{1}}
\newcommand{\eps}{\varepsilon}
\newcommand{\To}{\longrightarrow}
\newcommand{\norm}[1]{\left\Vert#1\right\Vert}
\newcommand{\abs}[1]{\left\vert#1\right\vert}

\newcommand{\ie}{\emph{i.e.,}}

\newcommand{\equald}{\stackrel{\mathrm{d}}{=}}
\newcommand{\probc}{\stackrel{\mathrm{P}}{\longrightarrow}}

\def\qed{ \hfill $\blacksquare$}  
\let\ga=\alpha \let\gb=\beta  \let\gd=\delta 
 \let\gh=\eta   \let\gk=\kappa \let\gl=\lambda     \let\gn=\nu   \let\gr=\rho \let\gs=\sigma  
\let\gx=\chi  
\let\gC=\Gamma \let\gD=\Delta  \let\gL=\Lambda

\newcommand{\cB}{\mathcal{B}}\newcommand{\cC}{\mathcal{C}}
\newcommand{\cE}{\mathcal{E}}\newcommand{\cF}{\mathcal{F}}
\newcommand{\cG}{\mathcal{G}}\newcommand{\cH}{\mathcal{H}}\newcommand{\cI}{\mathcal{I}}
\newcommand{\cK}{\mathcal{K}}\newcommand{\cL}{\mathcal{L}}
\newcommand{\cN}{\mathcal{N}}\newcommand{\cO}{\mathcal{O}}
\newcommand{\cP}{\mathcal{P}}\newcommand{\cR}{\mathcal{R}}
\newcommand{\cS}{\mathcal{S}}\newcommand{\cT}{\mathcal{T}}
\newcommand{\cX}{\mathcal{X}}
\newcommand{\cZ}{\mathcal{Z}}  
\newcommand{\vone}{\mathbf{1}}

\newcommand{\mvJ}{\boldsymbol{J}}

\newcommand{\mvY}{\boldsymbol{Y}}\newcommand{\mvZ}{\boldsymbol{Z}}

\newcommand{\mvt}{\boldsymbol{t}}

\newcommand{\mvgt}{\boldsymbol{\tau}}

\newcommand{\bB}{\mathbb{B}}

\newcommand{\dN}{\mathds{N}}

\newcommand{\dR}{\mathds{R}}

\newcommand{\dZ}{\mathds{Z}} 

\newcommand{\sB}{\mathscr{B}}

\newcommand{\sL}{\mathscr{L}}

\DeclareMathOperator{\E}{\mathds{E}}
\DeclareMathOperator{\pr}{\mathds{P}}

\DeclareMathOperator{\var}{Var}

\DeclareMathOperator{\ber}{Ber}
\DeclareMathOperator{\bin}{Bin}

\newcommand{\wh}[1]{\widehat{#1}}

\newcommand{\RNum}[1]{\uppercase\expandafter{\romannumeral #1\relax}}

\newcommand{\plim}{\ensuremath{\stackrel{\pr}{\rightarrow}}}		

\newcommand{\dd}{\mathrm d}

\newcommand{\floor}[1]{\lfloor #1\rfloor}

\newcommand{\invisible}[1]{}

\DeclareMathOperator{\unif}{Unif}

\newcommand{\wY}{\widehat{Y}}
\newcommand{\wX}{\widehat{X}}
\newcommand{\wZ}{\widehat{Z}}
\newcommand{\wE}{\widehat{\E}}
\newcommand{\wpr}{\widehat{\pr}}
\newcommand{\ball}{\mathsf{B}}
\newcommand{\ballsp}{\mathfrak{G}}

\DeclareMathOperator{\dist}{d}

\newcommand{\edeg}{\nu^{\mathrm{deg}}}
\newcommand{\nb}{\mathrm{Nbhd}}

\DeclareMathOperator{\tv}{TV}

\newcommand{\Fplus}{\mathcal{F}^{\scriptscriptstyle+}} 
\newcommand{\Ptv}[1]{\norm{#1}_{\operatorname{TV}}}
\DeclareMathOperator{\loc}{\mathrm{loc}}
\DeclareMathOperator{\parent}{par}

\begin{document}
\title[Duplication-divergence random graphs I]{Uniform duplication-divergence random graphs I:\\ Degrees and Local Limits}
\author[Dey]{Partha S.~Dey$^{\ast}$}
\author[Maitra]{Neeladri Maitra$^{\ast}$}
\author[Manthripragada]{Jaya M.~Manthripragada$^{\ast}$}
\author[Terlov]{Grigory Terlov$^{\dagger}$}

\date{\today}

\address{$^{\ast}$Department of Mathematics, University of Illinois Urbana--Champaign}
\address{$^{\dagger}$Department of Statistics and Operational Research, University of North Carolina, Chapel Hill}
\email{\{psdey,nmaitra,jm66\}@illinois.edu, gterlov@unc.edu }

\begin{abstract}
    We study a growing random graph in which each new vertex chooses a uniformly random parent, joins that parent with probability $r\in(0,1]$, and independently joins each of the parent's current neighbors with probability $p\in(0,1)$. Let $p_c$ be the unique solution of $p_c e^{p_c}=1$.  For $0<p<p_c$, we prove that the empirical degree distribution converges in probability to the stationary distribution $\pi$ of an explicit reset-and-growth Markov chain. Its  upper tail satisfies $\pi((y,\infty))\sim C_{p,r}y^{-\kappa}$, where $C_{p,r}>0$ and  $\gk>0$ solves $p^\gk+p\gk=1$. This confirms predictions from physics (Lambiotte, Krapivsky, Bhat and Redner, \textit{Physical Review Letters}, 2016). We further show that for every fixed radius $\ell$ and $p<p_c$, the empirical distribution of induced rooted $\ell$-neighborhoods converges to the invariant law of a graph-valued reset-and-growth Markov process. We give an explicit construction of the limiting graph using a marked ancestral genealogy. Its root component is almost surely finite for $r<1$, and infinite with exactly one end for $r=1$. Moreover, for every $0<p<1$, both the fixed-vertex degrees and the maximum degree grow on the scale $t^p$, and the vertex of maximum degree is eventually unique and persistent.
\end{abstract}

\maketitle

\begin{center}
    \vspace*{-.5cm}\includegraphics[height=9cm]{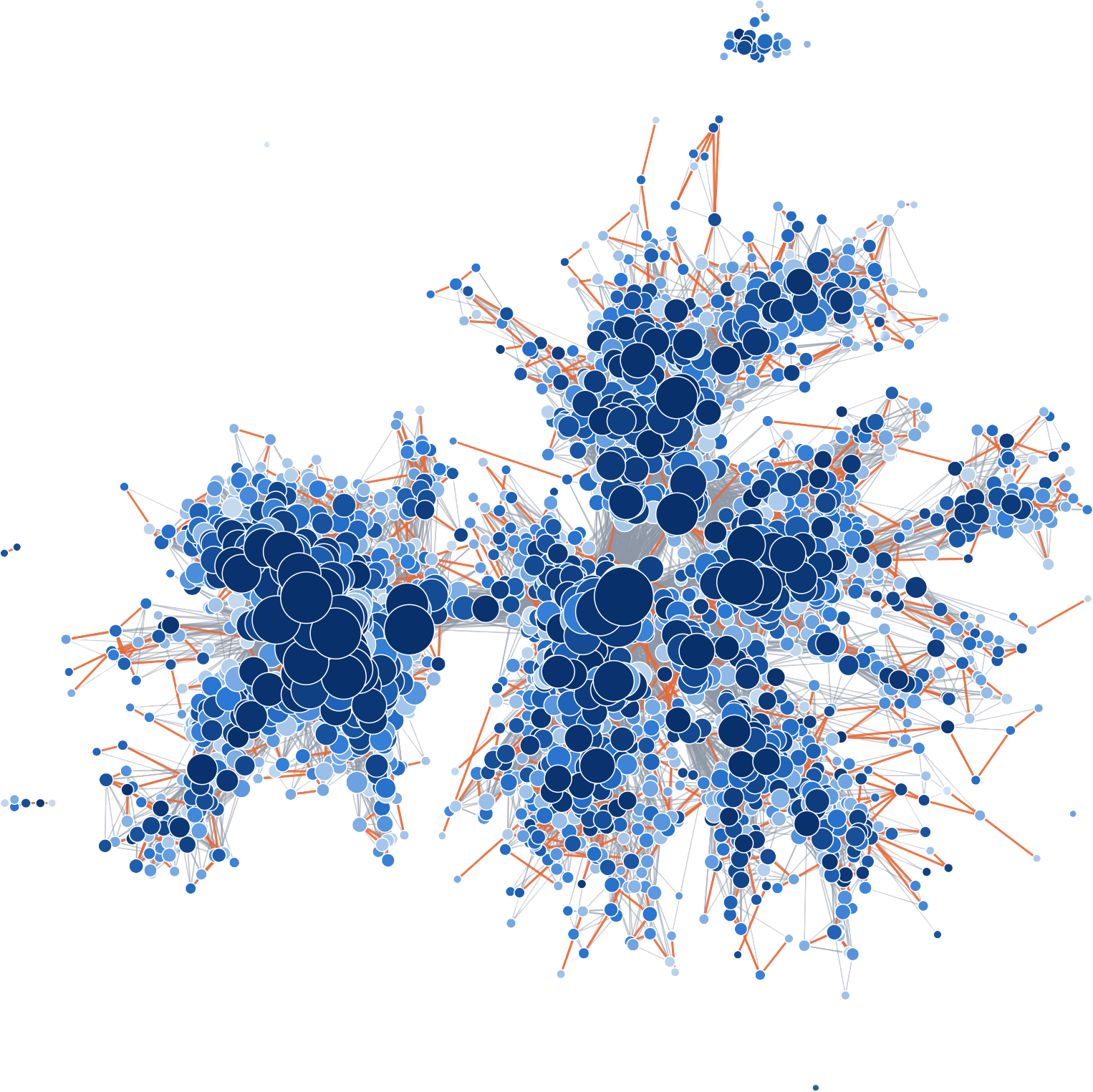}
    
    {{A network grown via the \emph{duplication-divergence} mechanism with 3000 vertices,\\ parent link probability $r=0.8$, and edge duplication probability $p=0.6$. }}
\end{center}

\setcounter{tocdepth}1\tableofcontents


\section{Introduction}\label{sec:intro}

The mathematical study of random graphs with evolution driven by edge \textit{duplication} or \emph{copying} draws motivation from the following three scientific directions. 
\begin{itemize}[leftmargin=1em]
\item \textbf{Gene duplication.} Starting in the early 1920s, the appearance of new genes in an organism's genome was experimentally explained by \emph{gene duplication}. Throughout the 20th century, it has been argued to be a key driver of evolution ~\cite{sturtevant1925effects, hogben1932causes, bridges1936bar, muller1936bar, stephens1951possible, ingram1961gene,ohno2013evolution}. 

This relates to protein-protein interaction (PPI) networks, that are among the most important biological networks to understand. With the emergence of complete PPI maps in the early 2000's~\cite{uetz2000comprehensive, ito2001comprehensive, jeong2001lethality}, it was often observed that these interaction networks are \emph{scale-free}, i.e., the degree distribution is controlled by a few \textit{hubs}. The scale-freeness of PPI networks have been also historically explained via gene duplication, see e.g.,~\cite{vazquez2003modeling,sole2002model,PSS03, emp_evolve_social_net, chung2003}.  For empirical evidence supporting the emergence of power-law degree tails (giving rise to scale-free networks) in PPI networks, we also refer the reader to~\cite{wagner2003global, teichmann2004gene, middendorf2005inferring}. For general interests in the biological motivations, let us also point the reader to the surveys~\cite{zhang2003evolution, taylor2004duplication}.

\item \textbf{Internet and social networks.} A second motivation comes from attempts to understand connectivity patterns in the internet and social networks. A new webpage or user may discover connections by consulting an existing webpage or user, and retaining some of those connections. Motivated by this mechanism, growing random graph models with the duplication mechanism modeling the web-graph have been suggested~\cite{kleinberg1999web, kumar2000stochastic}. In the context of social networks, the tendency of individuals forming connections with `friends of friends' is called \emph{triadic closure}. For empirical evidence of this phenomenon in various social networks, we refer the reader to~\cite{emp_evolve_social_net, leskovec2008microscopic}. Motivated by capturing this phenomenon, ~\cite{jackson2007meeting} proposes random graphs evolving via edge duplication as models for social network evolution.

\item \textbf{Network densification.} A third motivation comes from a purely network science perspective. \emph{Densification} of a network refers to the phenomenon of its average degree growing over time as its size increases, with new nodes joining. For empirical evidence of this phenomenon across different networks, we refer the reader to~\cite{leskovec2005graphs}. 
It is worth mentioning that previously physicists have looked at the densification phenomenon under the term `accelerating growth'~\cite{acc_growth}. The natural question: 
\[
\emph{What mechanisms cause networks to densify over time as they evolve?}
\]
Random graphs growing via the duplication mechanism have been offered as a plausible explanation for densification in networks~\cite{Ispolatov2005, lattanzi2009affiliation, LPBR16, Bhat_16, falkenberg2021heterogeneous}. The key idea is that if the duplication probability exceeds some model-dependent threshold, the network densifies, with the edge-count growing superlinearly in the network size. 
\end{itemize}

\noindent Motivated by the relevance of the duplication mechanism in qualitatively different networks that come from a variety of settings, we find it timely to rigorously study this model and its structural properties. The present paper is the first in a series of works accomplishing exactly that. 
Here we focus on degree distribution and local limits, thereby rigorously establishing power law predictions from \cite{LPBR16} and extending it to a more general setting. In addition we also discover several surprising behavior of the model, that to our knowledge has not been noticed before.

\subsection{Organization of the rest of the paper} 
Section~\ref{sec:main} defines the model and states the main results. We review related works in Section~\ref{sec:literature}. Section~\ref{sec:chain} studies the limiting degree chain, including recurrence and the stationary tail. Section~\ref{sec:gen} presents a general criterion for empirical convergence. Section~\ref{sec:local} connects this chain to empirical degrees and rooted neighborhoods. Section~\ref{sec:limit-law} contains proofs for an exact construction of a graph from the limit law. Section~\ref{sec:max-degree},
treats fixed and maximum degrees, and proves persistence of the hub. 
In Appendix~\ref{appendix: c>0 for p<1/2}, we present another brief proof of the positivity of a constant $C$ related to the upper-tail probability of the degree of a typical vertex, assuming $p\le 1/2$ (the general argument for any $p<p_c$ is in Section \ref{sec: proof of C>0 p<pc}).

\subsection*{Declaration on the use of generative AI}

No AI tools were used to aid with research or writing. The only uses of AI were to generate TikZ and Python code for our figures.
\subsection*{Acknowledgments}

G.T.\ thanks Istv\'an Kov\'acs, Bal\'asz R\'ath, and Sidney Redner for insightful conversations about the model and its applications. He also thanks P.S.D.\ for hospitality and arranging the visit to UIUC in April 2026, where all authors could discuss this project in-person.
G.T.\ was supported in part by the NSF RTG grant DMS-2134107.

\section{Model definition and main results}\label{sec:main}

\subsection{Model definition}

Fix two parameters $p,r\in[0,1]$ and a non-negative function $f$ on $\dZ_+$. For any graph $G$ and a vertex $u$ in it, denote by $D(u;G)$ the degree of $u$ in $G$. The graph evolves over discrete time steps as follows. 
\begin{itemize}
    \item At time $t=1$, we start with the graph $G_1$ consisting of the single vertex $1$.
    \item At time $t\ge 2$, given the graph $G_{t-1}$ at time $(t-1)$, a new vertex $t$ is added to the network.
    \item Conditional on $G_{t-1}$, the new vertex $t$ samples a \emph{parent} vertex $u_t$ from the graph $G_{t-1}$ with probability mass proportional to $f(D(u_t;G_{t-1}))$. 
    \item The vertex $t$ connects to the vertex $u_t$ via an edge with probability $r$ and to each neighbor of $u_t$ in $G_{t-1}$ with probability $p$. All these edge indicators are conditionally independent given $G_{t-1}$ and $u_t$.
\end{itemize}

\begin{remark}
    When $r=0$, we start the duplication-divergence model from a deterministic graph $G_{n_0}$ containing at least one edge, since starting from a singleton would leave the graph edgeless at every step. All subsequent updates follow the same rules as above. In the uniform-parent case $f\equiv1$, this is the classical \textit{partial duplication} model.
\end{remark}

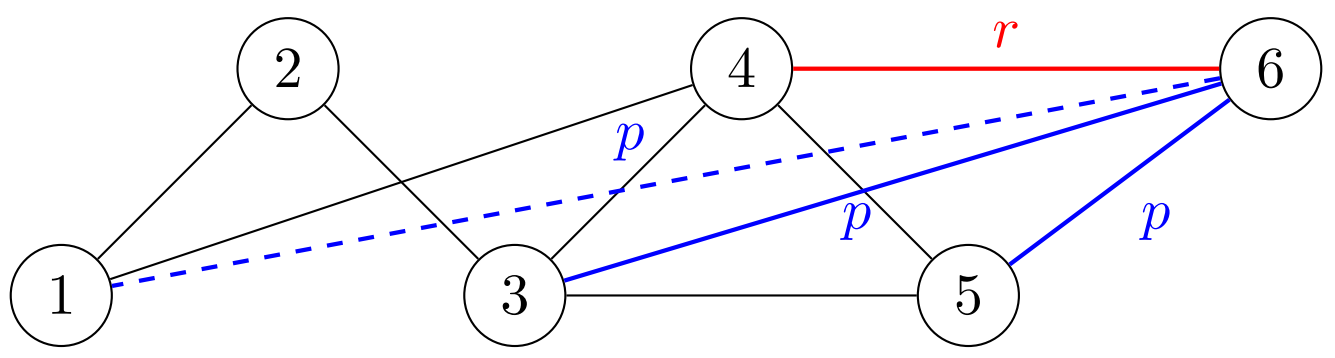
\begin{figure}[htbp]
\centering
\begin{tikzpicture}[scale=1, >=stealth]

  \node[draw, circle] (v1) at (0,0) {$1$};
  \node[draw, circle] (v2) at (1.5,1.5) {$2$};
  \node[draw, circle] (v3) at (3,0) {$3$};
  \node[draw, circle] (v4) at (4.5,1.5) {$4$};
  \node[draw, circle] (v5) at (6,0) {$5$};

  \draw (v1) -- (v2);
  \draw (v2) -- (v3);
  \draw (v3) -- (v4);
  \draw (v4) -- (v5);
  \draw (v4) -- (v1);
\draw (v5) -- (v3);
  \node[draw, circle] (v6) at (8,1.5) {$6$};

  \draw[red, thick] (v6) -- (v4)
      node[midway, above] {$r$};

  \draw[blue, thick] (v6) -- (v5)
      node[midway, below right] {$p$};
    \draw[blue, dashed, thick] (v6) -- (v1)
      node[midway, above left] {$p$};
  \draw[blue, thick] (v6) -- (v3)
      node[midway, below left] {$p$};

\end{tikzpicture}
\caption{A duplication-divergence step. Vertex $4$ is chosen as the parent, and the new vertex $6$ attaches to $4$ with probability $r$ (red edge), and to each neighbor of $4$ with probability $p$ (blue edges). In this example, the edges $\{6,5\}$ and $\{6,3\}$ occur, while the edge $\{6,1\}$ does not (dashed).}
\label{fig:duplication}
\end{figure}

\begin{figure}[htbp]
    \begin{center}
        \fbox{\includegraphics[width=.42\linewidth]{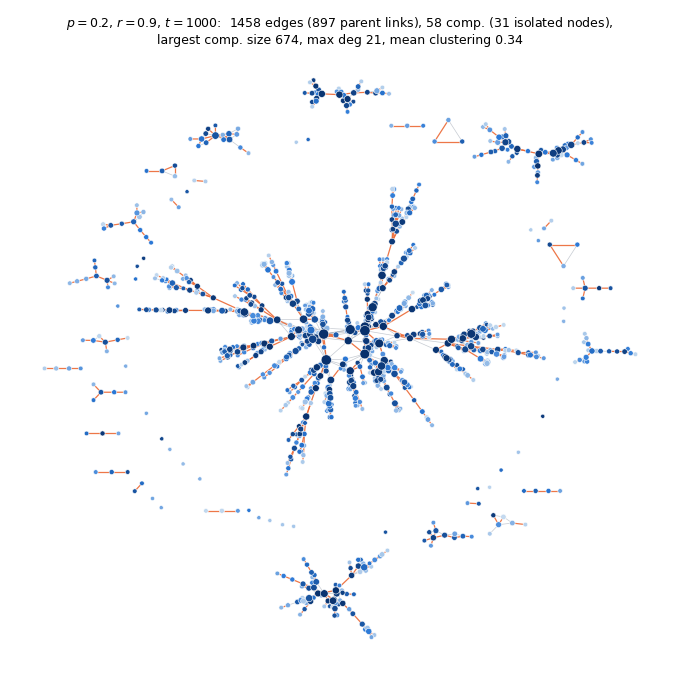}}%
        \fbox{\includegraphics[width=.42\linewidth]{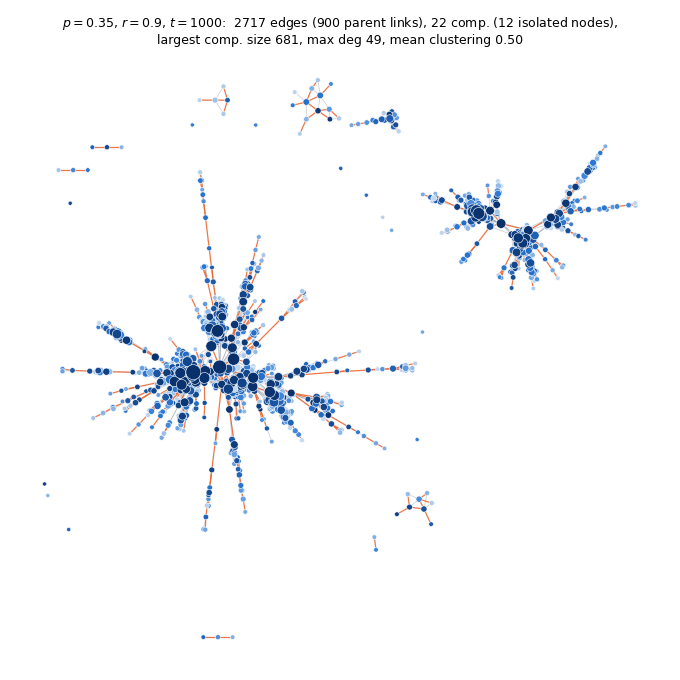}}\\
        \vspace*{-.3ex}\fbox{\includegraphics[width=.42\linewidth]{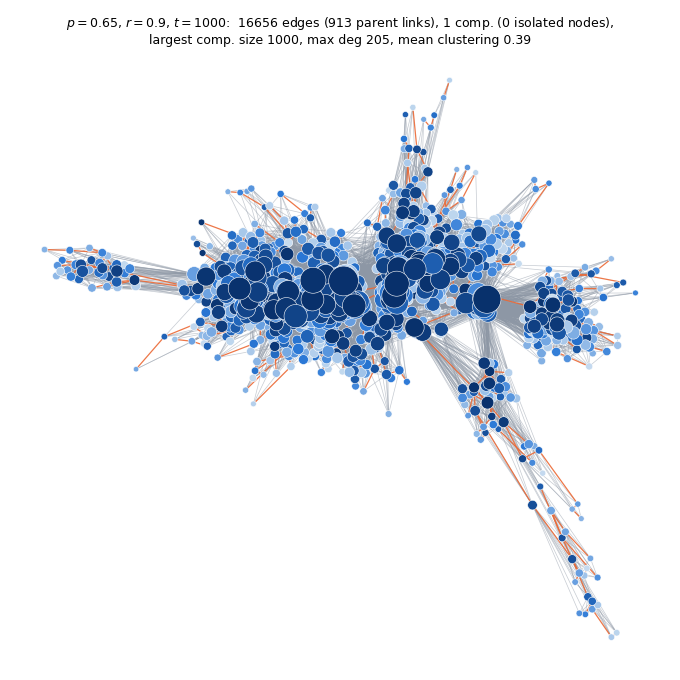}}%
        \fbox{\includegraphics[width=.42\linewidth]{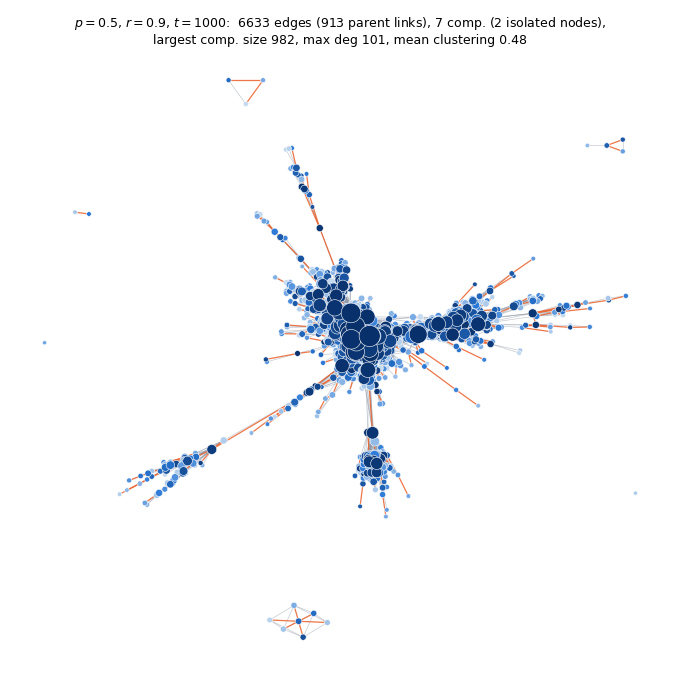}}
    \end{center}
     \caption{Realizations at time $t=1000$ with parent-link probability $r=0.9$ and $p\in\{.2,.35,.5,.65\}$ in clockwise order starting from the top left.  Red edges are parent links, and gray edges arise by duplication. Bigger vertices have higher degrees. 
     }
     \label{fig:copy-model-simulation}
\end{figure}

In this paper, we take $f\equiv 1$, so that the parent vertex is chosen uniformly from the existing vertices. We refer to this uniform-parent model as the \emph{uniform duplication-divergence model}.
We call $p$ the \textit{duplication probability} and $r$ the \textit{parent-link probability}. Unless a boundary case is stated explicitly, we assume that $p\in(0,1)$ and $r\in(0,1]$. See Figure~\ref{fig:copy-model-simulation} for simulated pictures for different values of $p$ and $r=0.9$.

We write $G_t=G_t(p,r)$ for the graph on $[t]$ after $t$ steps. Let $\deg_{t}(i)=D(i;G_{t})$ denote the degree of vertex $i\in[t]$ at time $t\ge 1$. 
For $j\ge 0$, let $N_{j}(t)$ be the degree count given by 
\begin{align}\label{def:degree-count}
 N_j(t):=\abs{\{v\in[t]:\deg_t(v)=j\}}
\end{align}
and let $\edeg_t$ be the empirical degree distribution measure at time $t\ge 1$ given by 
\begin{align}\label{eq:edeg}
 \edeg_t:=\frac1t\sum_{v=1}^t\gd_{\deg_t(v)},
\end{align}
which we can identify as the random probability vector $\bigl({N_j(t)}/t\bigr)_{j\geq0}$. We now assume that $p\in(0,1)$ and $r\in(0,1]$.  Unlike the classical model with $r=0$, the edgeless state is no longer absorbing: parent links continually create edges, while duplication amplifies the edges already present.  
Let $p_c\approx 0.5671$ be the unique solution to the equation
\begin{align}\label{def:pc}
    p_c e^{p_c}=1.
\end{align} 
For $p<p_c$, let  $\gk=\gk_p\in (0, \infty)$ be the unique positive solution to the equation 
\begin{align}\label{def:kappa}
    p^\gk +\gk p=1.
\end{align}

\begin{remark}
    Note that, $\gk\downarrow0$ as $p\uparrow p_c$. To see this, observe that by the definition of $p_c$, we have $\log p_c=-p_c$, so that the derivative of the function $p_c^\kappa+\kappa p_c$ in $\kappa \in [0,\infty)$ equals $p_c(1-p_c^\kappa)>0$. In other words, this function takes the value $1$ at $\kappa=0$ and stays larger than $1$ for all $\kappa>0$.
\end{remark}

Now, we are ready to state the main results of this article.

\subsection{The typical and the empirical degree distribution}
First, we present results for the limiting degree distribution of a typical vertex. This verifies predictions for the special case $r=1$ from \cite{LPBR16}. For $r>0$, the parent-link mechanism changes the limiting degree dynamics and leads to a different power-law exponent. Recall $\kappa=\kappa_p$ for $p<p_c$ from~\eqref{def:kappa}. Below by `$\implies$' we denote convergence in distribution.

\begin{thm}[Limiting degree distribution and power-law tail]\label{thm:degree-limit} 
Assume that  $0<r\le 1$. Let $U_T$ be a uniformly chosen vertex at time $T$ from the graph. If $0<p<p_c$, then
\[
\deg_T(U_T) \implies D_\infty\text{ as }T\to\infty.
\]
where $D_\infty$ is a non-trivial random variable supported on $\dZ_+$ and satisfies
\[
\pr(D_\infty > y) = (C+o(1)) \cdot y^{-\gk}\text{ as } y\to\infty
\]
for some constant $C=C_{p,r}\in (0,\infty)$. In particular, $\E D_\infty<\infty$ if and only if $0<p<1/2$. Moreover, for $p_c<p<1$, we have 
\[
\log (1+\deg_T(U_T)) /\log T\probc \log(pe^p)>0 \text{ as } T\to\infty.
\]
\end{thm}

The proof is presented in Section~\ref{ssec:pf-degree-limit}. It identifies $D_\infty$ concretely as the stationary limit of a reset-and-growth Markov chain.

\begin{defn}[The limiting reset-and-growth chain]\label{def:limit-chain}
For $x\in\dZ_+$, let $H_x$ be the $\{x,x+1,\ldots\}$--valued random variable with
\begin{align}\label{def:H_x}
\pr(H_x\ge k)=\prod_{i=x}^{k-1}\frac{r+pi}{1+r+pi},
\qquad k> x.
\end{align}
The chain $\mvY=(Y_n)_{n\ge0}$ on $\dZ_+$ is
defined by
\begin{align}\label{eq:XY-chain}
\begin{split}
    \textup{Reset step: }& (X_{n+1}\mid Y_n=y)\sim\ber(r)+\bin(y,p),\\
    \text{ and }
\textup{Growth step: }& (Y_{n+1}\mid X_{n+1}=x)\sim H_x.
\end{split}
\end{align}
\end{defn}

 $H_x$ admits the following equivalent description. Consider a pure-birth process started at $x$ with birth rate $(r+pi)$ when in state $i$. Then $H_x$ is the size of this process at time $E$, where $E$ is an independent exponential random variable of rate $1$. Intuitively, the reset step describes the degree at birth of a copy of a degree-$y$ parent. The growth step then runs this birth process for an independent exponential amount of time. 

\begin{thm}\label{thm:limit-chain}
We have the following for $0<p<1, 0<r\le 1$.
\begin{enumeratea}
\item If $0<p<p_c$, then $\mvY$ is positive recurrent with a unique stationary
distribution $D_\infty\sim\pi$. With $\gk$ as defined in~\eqref{def:kappa}, we have
\[
\lim_{y\to\infty}y^{\gk}\pr\bigl(D_\infty>y\bigr)=C\text{ for a finite }C>0.
\]
\item If $p>p_c$, then $\mvY$ is transient and $Y_n\to\infty$ almost surely.
\end{enumeratea}
\end{thm}
The proof is given in Section~\ref{ssec:pf-limit-chain}. We use Goldie's implicit renewal theorem~\cite[Corollary 2.4]{goldie91} to show the existence of the limit $C\ge 0$. The strict positivity of $C$ is established separately in Section~\ref{sec: proof of C>0 p<pc} by a Kesten--type excursion argument, and by a shorter convexity argument for $p\le1/2$ in Appendix~\ref{appendix: c>0 for p<1/2}.

Next, we study the evolution of the empirical degree distribution $\edeg_t$ over time.  Note that $\edeg_t$ is the conditional law of $\deg_t(U_t)$ given $G_t$. Moreover, $\edeg_t$ is supported on the set $\{0,1,\ldots,t-1\}$ for all $t$. For any function $\varphi:\dZ_+\to \dR$, let us write
\[
\edeg_t\varphi := \frac1t\sum_{v=1}^t \varphi(\deg_t(v)).
\]

\begin{thm}[Empirical degree distribution]\label{thm:edegree-lln}
Let $p\in(0,1)$, $r\in(0,1]$, and $p<p_c$.  Then as $t \to \infty$
\[
 \bigl\|\edeg_t-\pi\bigr\|_{\tv}\plim0.
\]
In particular, for every $j\in\dZ_+$, we have
$
 N_j(t)/t\probc \pi_j,
$
and, more generally, $\edeg_t\varphi\probc \pi\varphi$ for every bounded $\varphi\colon\dZ_+\to\dR$. We also have 
\[
\var\bigl(t \edeg_t\varphi\bigr)
\le  C\norm{\varphi}_\infty^2\cdot 
\begin{cases}
    t & \text{ if } 0<p<1/4,\\
    t\log t &\text{ if } p=1/4,\\
    t^{4p} &\text{ if } 1/4<p<1/e,\\
    t^{2p+2p\log\abs{e\log p}/\abs{\log p} } &\text{ if } 1/e\le p<p_c.
\end{cases}
\]
\end{thm}

The proof is presented in Section~\ref{ssec:pf-edegree-lln}. It combines an exact backward martingale representation of degree counts with a coupling of two copies of the discretized degree chain whose difference, raised to a small power, contracts at the rate $1-p^\ga-\ga p>0$ for some $\ga>0$. Such an $\ga$ is  available for every $p<p_c$.  The range $p<p_c$ is therefore the same as in Theorem~\ref{thm:degree-limit}. 
In particular, it includes the values $1/2\leq p<p_c$ at which the average degree 
diverges in probability while the empirical degree distribution still converges.

\subsection{Local neighborhood convergence}
We now extend the empirical degree distribution convergence result to induced rooted neighborhoods of every fixed radius  $\ell\ge 2$. For a finite graph $G$ with a distinguished root $o$ and an integer $\ell\geq1$, let $\ball_\ell(G,o)$ be the induced rooted subgraph on the vertices at graph distance at most $\ell$ from $o$, and let $\ballsp_\ell$ be the countable set of isomorphism classes of finite rooted graphs of radius at most $\ell$.  The empirical rooted $\ell$-neighborhood distribution is
\begin{align}\label{def:empirical-l-ball}
 \gn_t^{(\ell)}
 :=\frac1t\sum_{o\leq t}\gd_{\ball_\ell(G_t,o)},
\end{align}
viewed as a random probability measure on $\ballsp_\ell$. We call $\edeg_t$ the \emph{root-degree projection} of $\gn_t^{(1)}$. An induced radius-one ball also records edges between neighbors of the root, and therefore $\ballsp_\ell$ contains more information than the degree distribution alone. 
\begin{thm}[Empirical $\ell$-neighborhood distribution]\label{thm:l-ball-limit}
Fix $\ell\geq1$. Assume that $0<p<p_c$ and $0<r\leq1$.  There is an explicit graph-valued reset-and-growth Markov process on $\ballsp_\ell$, with generator $\sL_\ell$ explicitly given in~\eqref{def:l-ball-generator}, which has a unique stationary law $\Pi_\ell$, and 
\begin{align}\label{eq:k-ball-tv-limit}
	\left\|\gn_t^{(\ell)}-\Pi_\ell\right\|_{\tv}\plim0.
\end{align}
Moreover, the family is projectively consistent. Truncating a $\Pi_\ell$-distributed ball at radius $j\leq\ell$ recovers $\Pi_j$, and the root degree under $\Pi_\ell$ has distribution $\pi$.
\end{thm} 

The proof is presented in Section~\ref{ssec:pf-l-ball-limit}. As a consequence of Theorem~\ref{thm:l-ball-limit}, we establish the existence of the local weak limit of the duplication--divergence model graph $G_t$. For the definition of the \emph{local topology} in the statement below, we refer the reader to \cite[Chapter 2]{RGCN_2}. Also, when we write $(G_t,v)$ for some vertex $v \in G_t$, what we mean is the rooted graph (identified up to rooted isomorphism classes), that is the connected component of $v$ in $G_t$, and rooted at $v$. 
\begin{cor}[Existence of the local weak limit]\label{cor:local-weak-limit}
There exists a unique probability law $\Pi$ on locally finite rooted connected graphs such that, for $(G_\infty,o_\infty)\sim\Pi$,
\[
    B_\ell(G_\infty,o_\infty)\sim\Pi_\ell\qquad\text{for every }\ell\ge1.
\]
If $U_t$ is a uniformly chosen vertex of $G_t$, then $(G_t,U_t)\Longrightarrow(G_\infty,o_\infty)$ in the local topology.  
Moreover,
\[
    \gn_t^{\loc}
    :=\frac1t\sum_{v=1}^t\delta_{(G_t,v)}
\]
converges weakly in probability to $\Pi$. The degree of $o_\infty$ in $G_\infty$ has distribution $\pi$.
\end{cor}

For $0<p<p_c$, $0<r\le1$, let $\Pi$ be the locally finite rooted local limit, with root-degree law $\pi$. We construct the whole graph jointly. Its genealogy is the ancestral (or \emph{fringe}) limit of the random recursive tree, see~\cite{Aldous1991} and~\cite[Theorem~5.25 and Example~6.1]{HolmgrenJanson2017}. We decorate this genealogy with parent-link and duplication marks to obtain the law of the limiting graph, which we describe next.

\subsection{Construction of the limiting graph.}\label{sec:lim_graph_constn} 
\begin{figure}[!htbp]
    \centering
\includegraphics[width=1\linewidth]{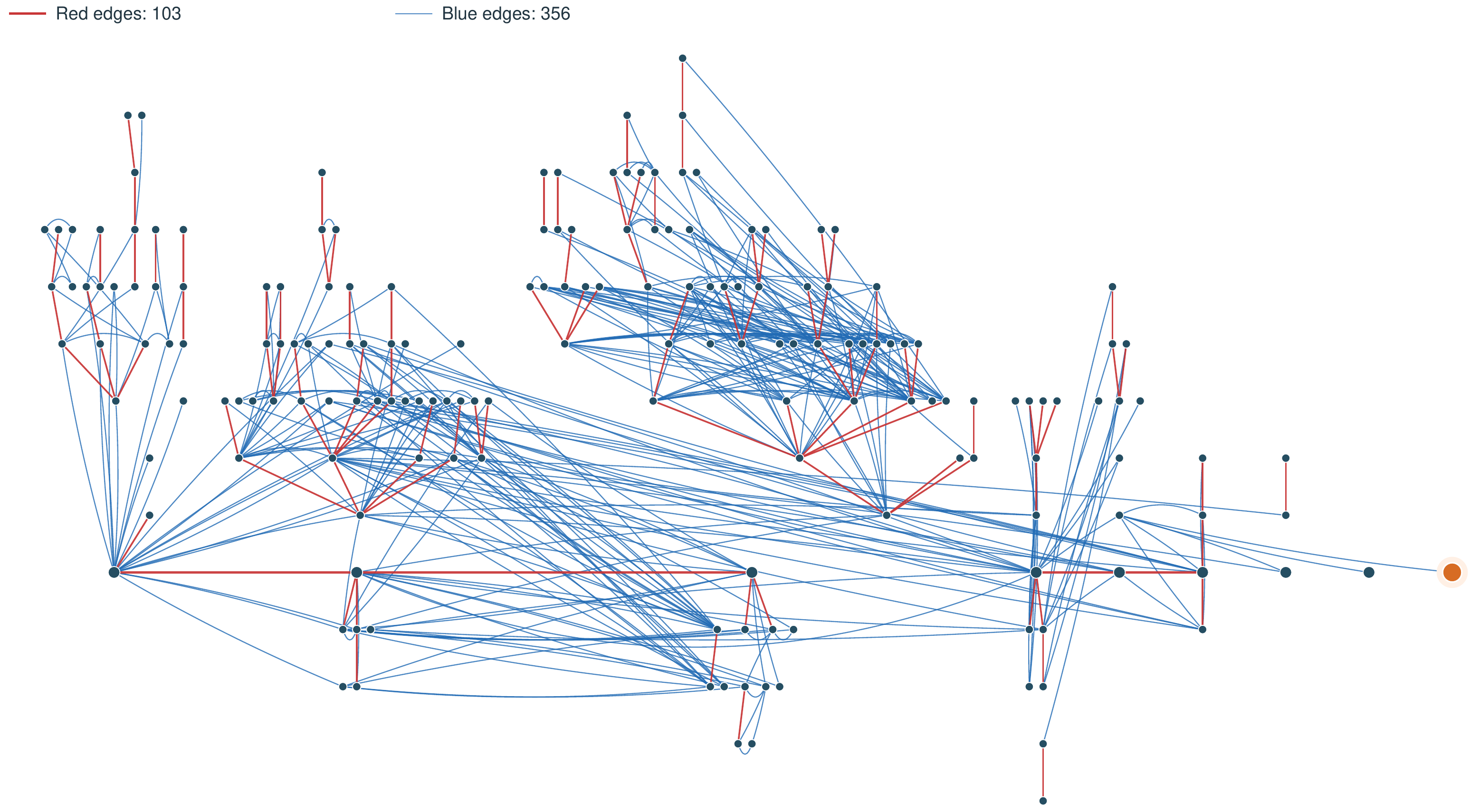}
    \caption{A sample graph from the constructed limit law for $p=.4, r=.5$. The orange vertex is $a_0$.}
    \label{fig:limit-law}
\end{figure}

Here, we give a rigorous joint construction of the limiting random graph in the local limit sense. See Figure~\ref{fig:limit-law} for a simulated picture for $p=.4,r=.5$.
\begin{enumerateA}
    \item \textbf{The Ancestral Tree $\cT$.} \label{def:local-lim A}
Let $(E_k)_{k\ge0}$ be i.i.d.~rate-one exponential random variables. Define $A_k:=\sum_{j=0}^k E_j$ for $k\ge 0$. We will first construct a genealogical tree $\cT$ with vertex labels in $(-\infty, 0]$, as follows. If a vertex has label $l$, then let us think $-l \geq 0$ to be its age.
\begin{enumerate}[label=$\bullet$,leftmargin=4ex]
    \item Construct an infinite \emph{spine} of vertices $a_0,a_1,a_2,\dots$ where vertex $a_k$ has label $-A_k$ (or age $A_k$). Think of these vertices placed on $(-\infty, 0]$ at their respective labels.
    \item For any vertex $a_k$ in the spine, $k \geq 0$, drop a finite rate one pure birth (Yule) process from it evolved up to $A_k$ units of continuous time. Let us call this finite genealogical tree $\mvt_k$, and for $k_1\neq k_2$, $\mvt_{k_1}$ is independent of $\mvt_{k_2}$. If a vertex in $\mvt_k$ has birth time $t$, then assign it the label $-A_k+t$.
    \item The ancestral spine $a_0,a_1,\dots$ together with the finite bushes $\mvt_0, \mvt_1,\mvt_2,\cdots$ respectively hanging off of them, where each vertex is labeled via the procedure described above, constitute the entire genealogical tree $\cT$. Note that in this tree each vertex has label in $(-\infty,0]$ by construction, and each vertex can be identified by its unique label.  
\end{enumerate}
We call $\cT$ the \emph{ancestral tree}.  For two vertices $x,y\in \cT$, we write $y\prec_{\cT} x$ if $y$ is younger than $x$.

    \item \textbf{Marks $\mvJ=(J_{x,y})_{x,y\in\cT}$ on the vertex pairs.} \label{def:local-lim B}
Conditional on the ancestral tree $\cT$, assign mutually independent marks
\begin{align*}
 	R_{\{x,y\}}&\sim\ber(r)\quad \text{for 
    parent--child pairs in $\cT$},\\
	P_{\{x,y\}}&\sim\ber(p)\quad \text{for all other pairs $x,y \in \cT$}.
\end{align*}
Let us also think of parent-child pairs $\{x,y\}$ with $R_{\{x,y\}}=1$ as \textit{red} edges on $\cT$. For any $x\in \cT$, define by $\parent(x)$ the parent of $x$ in $\cT$. In particular, we have $\parent(a_k)=a_{k+1}$ for $k\ge 0$. For parent--child pairs, \ie\ for pairs $x,y\in \cT$ with $d_\cT(x,y)=1$, define $J_{x,y}=J_{y,x}=R_{\{x,y\}}$. Otherwise, do the following.
\begin{enumerate}[label=$\bullet$,leftmargin=4ex]
    \item For any $x,y \in \cT$ with $d_\cT(x,y)=2$ and where $x\prec_\cT y$ (\ie\ the age of $x$ is lesser than that of $y$), define \begin{align}\label{eq:spine-pair-recursion}
	 J_{x,y}=J_{y,x}=P_{\{x,y\}}J_{\parent(x),y}.
\end{align} Observe that this is well defined since $\parent(x)$ is a child of $y$ in $\cT$, so that $d_\cT(\parent(x),y)=1$, and for such pairs $J_{\{\parent(x),y\}}$ is already defined.
\item We proceed inductively. Namely, given that $J_{x,y}$ is already defined for all pairs $x,y \in \cT$ with $d_\cT(x,y)=k$, for any pair $x',y'$ with $d_\cT(x',y')=k+1$ and $x'$ younger than $y'$, we define
\begin{align}
	 J_{x',y'}=J_{y',x'}=P_{\{x',y'\}}J_{\parent(x'),y'}.
\end{align}
As usual, this recursion makes sense since $J_{\parent(x'),y'}$ is already defined (since $d_\cT(\parent(x'),y')=k$).
\item We continue this procedure inductively and indefinitely, at each step considering vertex pairs with distance more by one than from the previous step. Since a younger vertex cannot be an ancestor of an older vertex, replacing the younger endpoint by its parent reduces genealogical distance by one. Thus this recursion is well defined by induction on $d_\cT(x,y)$. The younger endpoint may change during the recursion.
\end{enumerate}

For each unordered pair $e=\{x,y\}$ of vertices from $\cT$, let $\cl_\cT(e)$ be its \emph{trace}, a set of unordered vertex pairs from $\cT$, defined as follows. Repeatedly replace the younger endpoint by its parent until reaching a parent--child pair in $\cT$, and record every pair encountered, including the initial and terminal pairs. This trace contains $(d_\cT(x,y)-1)$ many non-parent--child pairs and exactly one parent--child pair. For a finite set $F$ of vertex pairs from $\cT$, define
\begin{align}\label{def:CRF}
\begin{split}
    \cl_\cT(F) &:=\cup_{e\in F}\cl_\cT(e),\\
    \cR(F) &:=\{e\in\cl_\cT(F):e\text{ is a parent--child pair}\},\\
    \cC(F) &:=\cl_\cT(F)\setminus\cR(F).
 \end{split}
\end{align}
Thus, for $F=\{\{x,y\}\}$, we have
$
\abs{\cR(F)}=1\text{ and }\abs{\cC(F)}=d_\cT(x,y)-1.
$
Moreover, by construction, it is easy to verify that, 
\[
J_{x,y}=\prod_{e\in \cl_\cT(\{x,y\})} M_e
\]
where $M_{\{u,v\}}=P_{\{u,v\}}$ if $d_\cT(u,v)>1$ and $M_{\{u,v\}}=R_{\{u,v\}}$ if $d_\cT(u,v)=1$.

    \item \textbf{The Rooted Connected Graph $(\cG_o,o)$.} Let $\cG$ be the graph with adjacency indicators $J_{x,y}$ (in other words, the vertex set is the same as that of $\cT$ and the edges are the union of the red and the blue edges). We think of this graph being rooted at $o:=a_0$. Let $\cG_o$ be the connected component of $o$ in $\cG$ and $o$ be the root of $\cG_o$. \label{def:local-lim C}
\end{enumerateA}

\begin{figure}[!htbp]
    \centering
    \begin{tikzpicture}[
  x=1.2cm, y=1cm,
  line width=1pt,
  line cap=round,
  vertex/.style={
    circle, draw, fill=white!95!black,
    minimum size=7mm, inner sep=0pt, line width=.2pt, font=\large
  },
  trace/.style={
    draw=blue!50!cyan,
    line width=0.9pt,
    dash pattern=on 3pt off 5pt
  }
]
  \node[vertex] (a) at (0,0)     {$a$};
  \node[vertex] (b) at (2.5,0)  {$b$};
  \node[vertex] (c) at (5,0)  {$c$};
  \node[vertex] (x) at (0,-1.5) {$x$};
  \node[vertex] (y) at (1.5,-1.5) {$y$};
  \node[vertex] (z) at (3.5,-1.5) {$z$};

  \draw (-1,-0.0) -- (a);
  \draw (a) -- (b) -- (c) -- (6,0.0);
  \draw (a) -- (x);
  \draw (b) -- (y);
  \draw (b) -- (z);

  \draw[trace] (a.south east)
    to[out=-18,in=200] (b.south west);
    \draw[trace] (x.south east)
    .. controls +(0.95,-0.40) and +(-0.5,-2.5) .. (c.south);
  \draw[trace] (x.north east) -- (b.south west);
  \draw[trace] (x.east) -- (y.west);
  \draw[trace] (x.south east)
    .. controls +(0.95,-0.40) and +(-0.85,-0.45) ..
    (z.south west);
\end{tikzpicture}
    \caption{In this part of the genealogical tree $\cT$ picture, suppose we have $z\prec_\cT y\prec_\cT x\prec_\cT c\prec_\cT b\prec_\cT a$, so that $a$ is the oldest and $z$ is the youngest vertex. Then, $\cl_\cT(\{x,y\})=\{\{x,y\},\{x,b\},\{a,b\}\}, \cl_\cT(\{x,z\})=\{\{x,z\},\{x,b\},\{a,b\}\}$ and $\cl_\cT(\{x,c\})=\{\{x,c\},\{x,b\},\{a,b\}\}$.}
    \label{fig:trace}
\end{figure}
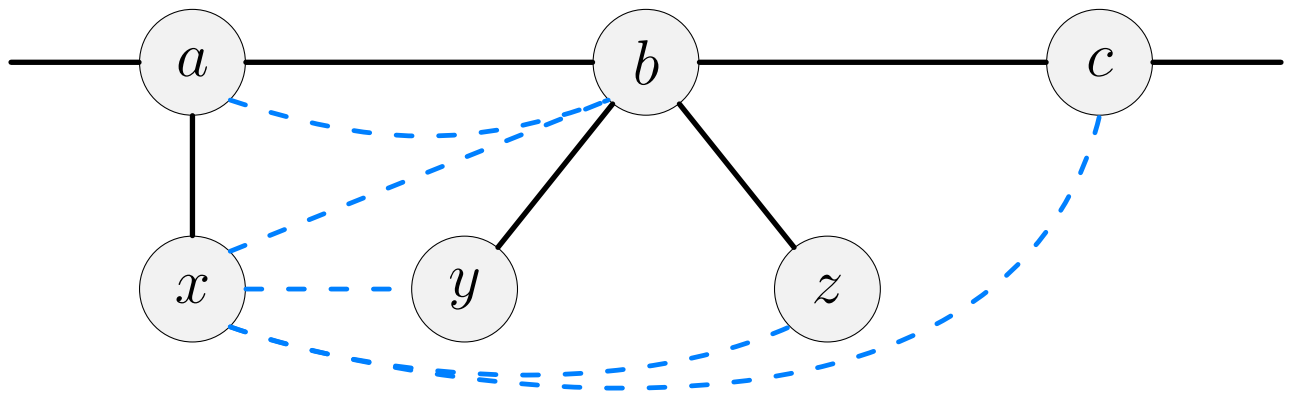

We think of there being a blue edge between vertex pairs $\{x,y\}$ with $x,y$ non-adjacent in $\cT$ but $J_{x,y}=1$.  Moreover, a blue edge $e=\{x,y\}, d_\cT(x,y)\ge 2$, is present exactly when all marks in 
$\cl_\cT(e)$ are one. See Figure~\ref{fig:trace} for an illustration.
We claim that the graph $\cG_o$ is the limiting random graph in the local sense. We prove the following result. 
\begin{thm}\label{thm:spine-local-law}
Almost surely, $\cG$ is locally finite and $(\cG_o,o)\sim\Pi.$ For distinct pair of vertices $x,y$,
\begin{align}\label{eq:spine-edge-marginal}
 \pr(J_{xy}=1\mid\cT)=rp^{d_{\cT}(x,y)-1}.
\end{align}
In general, for a finite set $F$ of unordered pairs of distinct vertices from $\cT$, we have
\begin{align}\label{eq:spine-joint-law}
 \pr(F\subseteq E(\cG)\mid\cT) =p^{|\cC(F)|}r^{|\cR(F)|},
\end{align}
where $\cR(F),\cC(F)$ are as defined in~\eqref{def:CRF}. Inclusion--exclusion determines every finite conditional adjacency law.
\end{thm}
The proof is presented in Section~\ref{ssec:pf-spine-local-law}. Moreover, we have the following results regarding the structure of a sample from $\Pi$.

\begin{prop}\label{prop:limit-component-dichotomy} Assume that $p<p_c$, where $p_c$ is as given in~\eqref{def:pc}. For $(G_\infty,o)\sim\Pi$, almost surely we have
\[
 \begin{cases}
 \abs{V(G_\infty)}<\infty,&\text{ if }0<r<1,\\
 G_\infty\text{ is infinite and has exactly one end},&\text{ if } r=1.
 \end{cases}
\]
\end{prop}
Here one end means that deleting any finite vertex set leaves exactly one infinite component.
A size-biased exponential-age description is also available when $r<1$.

\begin{prop}\label{prop:component-exp-law}
Assume $0<r<1$. Let $C_t$ be the ordinary unrooted component-growth process started from a singleton: each vertex has a rate-one duplication clock, and unattached newborns are discarded. Let $E\sim\operatorname{Exp}(1)$ be independent of this process. Then 
\begin{align}\label{eq:component-exp-normalizer}
	M:=\E\abs{V(C_E)}=1/((1+r)\pi_0)<\infty.
\end{align}
For every bounded rooted-graph functional $F$,
\begin{align}\label{eq:component-exp-law}
	 \Pi F:=\E F(G_\infty,o)=\frac1M\E \sum_{v\in V(C_E)}F(C_E,v),
\end{align}
where $(G_\infty,o)\sim \Pi$. Equivalently, bias the law of $C_E$ by its size and choose a uniform root. In particular, conditional on the unrooted component, the root is uniform.
\end{prop}
The proofs are presented in Section~\ref{sec:limit-law}.

\subsection{Fixed vertices and the maximum degree}\label{ssec:maxdegree}

The results so far concern a uniformly chosen vertex.  The degree of a
\emph{fixed} vertex grows like $t^p$, for every $p\in(0,1)$, and so does
the maximum degree, although for $p<p_c$ the empirical degree
distribution has the power-law tail $y^{-\gk}$ with $\gk<1/p$.  Define
$$D_{\max}(t):=\max_{v\leq t}\deg_t(v).$$ The variables $(\xi_v)_{v\ge1}$ below are defined on the probability space of the single graph process, as the simultaneous almost-sure limits of its normalized degrees. 

\begin{thm}[Fixed vertices and the maximum degree]\label{thm:max-degree}
Let $0<p<1$ and $0<r\leq1$.
\begin{enumeratea}
\item For every vertex $v\geq1$,
\[
t^{-p}\cdot \deg_t(v)\To \xi_v
 \text{ almost surely and in }L^2,
\]
where $\xi_v\in(0,\infty)$ almost surely, and the convergence holds
jointly for every finite set of vertices. Moreover, we have $\E\xi_1=\E\xi_2=r/(p\gC(1+p))$.
\item Almost surely,
\[
t^{-p}\cdot D_{\max}(t)\To\Xi:=\sup_{v\geq1}\xi_v\in(0,\infty),
\]
and the supremum is a maximum. For every $\eps>0$, the set
$\{v\mid \xi_v>(p+\eps)\Xi\}$ is almost surely finite. So
$\Xi=\max_{v\leq V_0}\xi_v$ for a finite random $V_0$, and for all
large $t$ the maximum degree is attained by a vertex $v\leq V_0$.  More
generally, if the $k$-th largest value $\xi_{(k)}$ of $(\xi_v)_{v\geq1}$
exceeds $p \Xi$, then the $k$-th largest degree at time $t$, divided by $t^p$, converges almost surely to $\xi_{(k)}$. 
\item The law of each $\xi_v$ has no atoms, and $\pr(\xi_v=\xi_w)=0$
for $v\neq w$.  Consequently the maximum $\Xi$ is attained at a unique
random vertex $V^*$, and there is a finite random time after which
$V^*$ is the unique vertex of maximum degree, \ie\ \emph{the hub persists}.
More generally, all order statistics above $p \Xi$ are distinct, and
the vertices carrying them are eventually fixed.
\end{enumeratea}
\end{thm}
The proofs of parts a), b),~and c) are presented in Sections~\ref{ssec:pf-max-degree-a},~\ref{ssec:pf-max-degree-b} and~\ref{ssec:pf-max-degree-c}. 
The normalization of a fixed degree follows from its affine conditional drift, as in martingale analyses of preferential growth. Here a newborn vertex can already inherit a substantial fraction of its parent's degree due to duplication. The edge-count scale $t^{2p}$ is therefore compatible with a maximum degree of order $t^p$. Edge growth involves the descendants of many edges, rather than the degree of a single hub.

The cutoff $p\Xi$ in Theorem~\ref{thm:max-degree} is a separation threshold. Its proof controls late vertices uniformly above this level; it does not assert convergence of order statistics at or below it. When $p<p_c$, the stationary power-law degree distribution and the maximum-degree scale also concern different limits. Indeed, $\kappa<1/p$, since $p^{1/p}+1>1$ in~\eqref{def:kappa}. The theorem about a typical vertex takes $t\to\infty$ at a fixed degree threshold and does not provide a uniform tail approximation at thresholds growing with $t$. Persistence of the hub is known for preferential attachment model~\cite{DM09,banerjee2021persistence},
where it is proved through the continuous-time embedding of the degree process.  Here the persistence comes from the factor $p<1$ lost at each copy, and it holds for every $p\in(0,1)$.

\subsection{Outline of the proofs}\label{sec:outline}
The proof separates the stationary degree chain from its connection to the graph.
\begin{enumeraten}
\item \textbf{The stationary chain (Section~\ref{sec:chain}).}
A newborn copy of a degree-$d$ parent has degree $\ber(r)+\bin(d,p)$, while an existing degree grows at rate $(r+pd)/t$. On logarithmic time these operations give the generator $\cL$ in~\eqref{eq:degree-generator}. Foster--Lyapunov estimates yield positive recurrence for $p<p_c$ and transience for $p>p_c$. A coupling with the multiplier $M=pU^{-p}$ permits an application of Goldie's implicit renewal theorem~\cite{goldie91}; the equation $\E M^\kappa=1$ gives~\eqref{def:kappa}. A stopped excursion argument proves positivity of the tail constant.
\item \textbf{Empirical convergence and local limit properties (Sections~\ref{sec:gen},\ref{sec:local} and \ref{sec:limit-law}).}
The empirical degree measure satisfies the exact conditional recursion with $K_t=I+\cL/(t+1)$. Euler approximation, tightness and a summable localized noise bound give convergence through Proposition~\ref{prop:empirical-convergence}. A backward martingale and coupling contraction additionally give the variance upper bounds in Theorem~\ref{thm:edegree-lln}. For rooted balls, the main additional ingredients are a finite-influence bound and regeneration of the graph-valued chain. These give all fixed radii throughout $p<p_c$, without requiring a finite $\ell$-th degree moment. The results on the properties of the local limit are in Section \ref{sec:limit-law}.
\item \textbf{Fixed and maximum degrees (Section~\ref{sec:max-degree}).}
An affine degree martingale gives the fixed-vertex limits. Uniform control of newborn degrees and their later fluctuations confines the largest degrees to finitely many old vertices. An extra-edge coupling rules out atoms and ties in their limiting degrees, yielding a persistent unique hub.
\end{enumeraten}

\subsection{Notation}
Throughout the paper, all asymptotic notation is understood as $n\to\infty$, unless stated otherwise. For two sequences $(a_t)$ and $(b_t)$, with $b_t>0$ for all sufficiently large $t$, we write $a_t=O(b_t)$, if there exist constants $C<\infty$ and $t_0\in\dN$ such that
$
|a_t|\leq Cb_t\text{ for all }t\geq t_0.
$ We write $a_t=o(b_t)$ if
$
\lim_{t\to\infty}a_t/{b_t}=0.
$
Finally, we write $a_t=\Theta(b_t)$ if there exist constants $0<c<\infty$ and $t_0\in\dN$ such that
$
c^{-1}b_t\leq |a_t|\leq cb_t\text{ for all }t\geq t_0.
$
Equivalently, $a_t=\Theta(b_t)$ means that $a_t=O(b_t)$ and $b_t=O(a_t)$.

\section{Related literature}\label{sec:literature}
In this section we review literature on structural properties of the graph sequence $(G_t(p,r))_{t\geq 1}$.
The \emph{partial duplication model} with $f\equiv 1$ and $r=0$ was introduced in~\cite{Bhan2002} and has since been studied extensively; see, for example,~\cite{chung2003, Li2013, Hermann2016, Jordan2018, Jacquet2020}. Closely related duplication-divergence models, motivated largely by biological network formation, have also been considered in~\cite{Ispolatov2005, Bebek2006, PSS03, FTS24}. At the other extreme, when $f\equiv 1$ and $r=1$, each new vertex first attaches to a uniformly chosen parent and then independently connects to each neighbor of that parent with probability $p$. This \textit{copying} mechanism has been studied in the network-science and physics literature; see~\cite{KR05, Bhat_16, LPBR16, tws24}. 

\subsection*{Partial duplication model with $r=0$.} 
The classical partial duplication model ($r=0$ and $p\in(0,1)$) is well studied in the literature. A feature of this model is that the degree distribution of the full graph, including isolated vertices, has no nontrivial power-law limit. Hermann and Pfaffelhuber~\cite{Hermann2016} identified a phase transition at $p_c \approx 0.567143$, the unique solution of $p_ce^{p_c}=1$. For $p\le p_c$, the proportion of isolated vertices converges almost surely to one, and hence the empirical degree distribution converges to $\gd_0$. For $p>p_c$, the limiting expected degree frequencies are defective, \ie\, positive mass escapes to degrees tending to infinity, so they do not define a probability distribution on the nonnegative integers. In~\cite{hermann2019partial}, the same authors studied the effect of partial duplication with \emph{edge-deletion} and gave results on the degree-distribution and certain subgraph counts. 

Li, Choi, and Wu~\cite{Li2013} studied expected degree-frequency recurrences and established a phase transition for the proportion of isolated vertices, thereby answering a question from~\cite{Bebek2006}.  They showed that the expected degree tends to zero for $p<1/2$, stays of constant order at $p=1/2$, and diverges for $p>1/2$.  Thus, for $1/2<p\leq p_c$, almost all vertices are isolated even though the expected average degree diverges; a vanishing fraction of high-degree vertices carries the excess degree.

Although the degree distribution of the full graph is degenerate in the subcritical regime $(p<p_c)$, when started from a single edge, its unique non-isolated connected component has a nontrivial degree distribution. For $0<p<e^{-1}$, Jordan~\cite{Jordan2018} proved that the degree distribution within this component converges in probability to a limiting distribution $(\pi_k)_{k\ge1}$. Moreover, this distribution has an approximate power-law behavior with exponent $\gk_0=\gk_0(p)$, \ie~$\pi_k = k^{-\gk_0 + o(1)}$, where $\gk_0>2$ is the nontrivial solution of
$
    p^{\gk_0-2}=3-\gk_0,
$
confirming predictions of Ispolatov, Krapivsky and Yuryev~\cite{Ispolatov2005}. The heuristics from~\cite{Ispolatov2005} further predicts the exponent $\gk_0=2$ for $e^{-1}\le p<1/2$, a regime not covered by Jordan's theorem. Jacquet, Turowski and Szpankowski~\cite{Jacquet2020} later sharpened Jordan's result for $0<p<e^{-1}$, proving the exact power-law asymptotic
$
    \pi_k = (C_{p}+ o_k(1))\cdot k^{-\gk_0}
$
for an explicit constant $C_p>0$. 
Barbour and Lo~\cite[Theorem 6.5]{BL22} proved a central limit theorem for the
logarithm of the degree of a randomly chosen vertex for a general duplication-divergence model.
The Gaussian limit describes the mass escaping to large degrees and the complementary mass is absorbed at degree zero.

Hermann and Pfaffelhuber~\cite{Hermann2016} obtained precise asymptotics for several structural statistics. In particular, if $I_t$ denotes the number of edges and the initial graph contains at least one edge, then,
$
    t^{-2p}\cdot I_t \implies I_\infty
$
almost surely and in $L^2$, to a nontrivial limit $I_\infty$. The average degree is therefore of order $t^{2p-1}$, yielding a separate phase transition at $p=1/2$. They also obtained a detailed description of the large-scale behavior of the partial duplication model. 
For every initially non-isolated vertex $k$, they proved that $t^{-p}\deg_t(k)$ converges almost surely and in $L^q$, for every $q\ge1$, to a nondegenerate limit. Thus, the degree of a fixed vertex grows on the scale $t^p$. Together, these results provide a precise description of the large-scale geometry for the $r=0$ model and identify the natural growth scales of several fundamental graph statistics.

A related model was considered by Kim, Krapivsky, Kahng and Redner in~\cite{kim2002infinite}, where $r=0$ and $p>0$, and in addition, an incoming vertex at step $t$ joins \emph{all other} existing vertices with probability $\beta/t$. The model can be thought of as combining the effect of duplication with \emph{uniformly growing random graphs}, see~\cite[Chapter 3]{RGCN_2}. Using non-rigorous techniques, the authors predict infinite-order phase transitions for giant components in this model when $p$ is sufficiently small, and that the fluctuation of the size of the giant undergoes phase transitions. Some rigorous results on the tail of the degree distribution of this model is also available in~\cite{Bebek2006}, asymptotic results on the degree of a fixed vertex in~\cite{turowski2021towards}, large-deviation results on the degree of a fixed or uniform vertex in ~\cite{frieze2020degree} and concentration of the maximum degree in~\cite{FTS24}. Recently, compression algorithms for partial duplication models have been studied in~\cite{Turowski2020Compression}.

\subsection*{Duplication-divergence models with parent links $r\in (0,1]$} Introducing parent-links (\ie\ $r>0$), however, changes the process qualitatively. They continually seed new edges, prevent the edgeless state from being absorbing, and create structures that need not be present initially. The resulting recurrences contain both an \emph{immigration} term generated by parent-links and a \emph{reinforcement} term generated by duplication. We summarize literature for this case here.

When $r=1$, every new vertex attaches to its parent and then to each neighbor of the parent with probability $p$.  This is the duplication-divergence model of Krapivsky and Redner~\cite{KR05}, and the ``friend of a friend'' model of~\cite{tws24}. Bhat, Krapivsky, Lambiotte and Redner~\cite{Bhat_16,LPBR16} predicted, by non-rigorous rate equations, that the number of edges grows linearly for $p<1/2$ and like $t^{2p}$ for $p>1/2$, and that the degree distribution has a power-law tail whose exponent solves $p^\gk+p\gk=1$. {This falls into the scope of our Theorem~\ref{thm:degree-limit}}. A similar power-law prediction also appears in a work of Raval~\cite{raval2003some}, where he makes further predictions on the rate of convergence of the degree distribution to its limit, and conjectures on the possibility of undertaking \emph{network archaeology} in these models. Duplication with a uniformly chosen parent also appears in the models for the web graph of~\cite{kleinberg1999web,kumar2000stochastic}, where a new page copies the out-links of a random earlier page. 

For $r>0$,~\cite{ispolatov2005cliques} studies clique-counts in these graphs using non-rigorous methods, and predicts the abundance of $4$ and $5$-cliques. The effect of the \emph{seed} graph, i.e., the initial graph at time $t=0$ from which the dynamics begin, on the evolution of the network was studied in~\cite{hormozdiari2007not}. Knudsen and Wiuf~\cite{Knudsen2008} developed a Markov chain approach to the degree dynamics, including a continuous-time approximation. For $r>0$, they proved recurrence of the associated jump chain for $p<p_c$, ergodicity for $p<1/2$, and transience for $p>(\sqrt{5}-1)/2$. However, they did not explicitly identify the stationary degree distribution. This paper establishes the stationary degree tail, proving the predictions of~\cite{Bhat_16} for the degree law, empirical degree convergence, induced local limits, and persistent-hub behavior for the uniform-parent model with $r>0$. 

Backhauz and Móri~\cite{backhausz2015asymptotic, backhausz2016further} study the effect of duplication together with \emph{vertex deletion}, and show that under some parameter ranges the degree distribution is \emph{stretched-exponential} instead of a power-law and give bounds on the maximum degree. More recently,~\cite{lo2025isolated} examines the effect of edge-deletion together with partial duplication, and gives results on isolated vertices. For a recent review article on the model see~\cite{borrelli2025duplicationdivergencegrowinggraphmodels}.

\subsection*{Comparison with preferential attachment}
In the affine preferential attachment model~\cite{barabasi1999emergence, BRST2001,hofstad17}, the degree of a fixed vertex grows like $t^{1/(2+\delta)}$, the degree distribution has a power law with exponent $(3+\delta)$, and the local weak limit is a random tree~\cite{BBCS14}.  The hub of maximum degree  persists~\cite{Mori05,DM09,banerjee2021persistence}. The duplication-divergence model shares the persistent hub and the power law, but the mechanisms differ. Duplication is a \emph{local} form of reinforcement -- a vertex gains edges only when one of its neighbors is chosen as a parent, at rate $p\deg(v)/t$, so the reinforcement is preferential in effect while the parent choice is uniform.  The degree of a fixed vertex grows like $t^p$ for every $p\in(0,1)$ (Theorem~\ref{thm:max-degree}), and the tail exponent $\gk=\gk_p$ of the typical degree can be any positive number, with $\gk<1/p$.  Consequently the maximum degree, of order $t^p$, is much smaller than the maximum of $t$ independent samples from $\pi$, which is of order $t^{1/\gk}$.  The local limit is not a tree, since copied edges create triangles, and the degree of the root has a power-law tail with a $p$-dependent exponent. 

Cohen, Jordan and Voliotis~\cite{cohen2010preferential} study a perferential version of the duplication model, where the parent is being chosen with probability proportional to its current degree (unlike the uniform case that we study), and show that certain descendants of initital vertices dominate the graph evolution, a phenomenon they call `trapping', using a comparison with vertex-reinforced random walks.

\subsection*{A Markovian approach to local limits}
Theorem~\ref{thm:l-ball-limit} is a local weak limit in the sense of Benjamini and Schramm~\cite{BS01}, and Aldous and Steele~\cite{AS04} (for a nice introduction to the topic we also refer the reader to~\cite[Chapter 2]{RGCN_2}), in the strong form of convergence in probability of the empirical neighborhood distribution. In fact, as the result shows, the convergence to the limiting $\ell$-neighborhoods is in the almost sure sense. The limit is described as the stationary law of a graph-valued reset-and-growth process, rather than through an explicit construction of the limiting rooted graph and coupling it with graph exploration in the pre-limit. It is worth mentioning that the latter strategy has been the most familiar technique to establish local limits (see the local limit proofs in~\cite{RGCN_2, maitra2026multipartite, rath2026local}) in the literature; in this paper we take a difference perspective. Namely, the local limit is established relying heavily on the graph dynamics itself, viewing the evolution of the neighborhood of a typical vertex as a Markov process on the space of rooted graphs, and invoking Markov chain convergence theorem to obtain the limit. As a result, the convergence is also almost sure, instead of just in probability. 
The explicit construction of the limiting rooted graph is in Section \ref{sec:lim_graph_constn}.

\section{The limiting reset-and-growth chain}\label{sec:chain}
We now isolate the one-dimensional Markov chain that governs the degree of a single vertex, and prove the stationarity and tail statements behind Theorem~\ref{thm:degree-limit}.  Section~\ref{sec:local} then shows that this chain indeed gives the local limit of the graph.
The starting point is the birth rate of a single degree.  If a vertex has degree $d$ at time $s-1$, its degree increases at time $s$ precisely when it is selected as the parent and the parent edge is present, or when one of its $d$ neighbors is selected and the corresponding edge is copied.  These alternatives are disjoint, so 
\begin{align}\label{eq:fixed-vertex-birth-rate}
	\pr\bigl(\deg_s(v)=d+1\mid\cF_{s-1},\deg_{s-1}(v)=d\bigr) =\frac{r+pd}{s-1}.
\end{align}

\subsection{The chain and its transition law}
For $x\in\dZ_+$, let $H_x$ be supported on
$\{x,x+1,\ldots\}$ with
\begin{align}\label{eq:H-law}
\pr(H_x\geq k)=\prod_{i=x}^{k-1}\frac{r+pi}{1+r+pi},
\qquad k\geq x.
\end{align}
To interpret this law, start a continuous-time pure-birth process at $x$, give state $i$ birth rate $r+pi$, and kill the process at an independent exponential time of rate $1$. Then $H_x$ is the terminal size of the process.  Indeed, to reach $k$, the process must win the birth-versus-killing race in every state $i=x,\ldots,k-1$, which gives~\eqref{eq:H-law}.  Equivalently, if $\cS_x$ counts the number of steps through the corresponding discrete chain, including the absorbing step, then $H_x\equald x+\cS_x-1$.  The recursion
\[
	\E\cS_x=1+\frac{r+px}{1+r+px}\E\cS_{x+1}
\]
has the solution
\begin{align}\label{eq:H-mean}
	\E\cS_x=\frac{1+r+p(x-1)}{1-p}, \quad \E H_x=\frac{r+x}{1-p}.
\end{align}

\begin{figure}[ht]
\centering
\begin{tikzpicture}[>=stealth,font=\small,x=1.6cm]
  \tikzset{st/.style={circle,draw,minimum size=7mm,inner sep=0pt}}
  \foreach \j in {0,...,5} {\node[st] (s\j) at (\j,0) {$\j$};}
  \node at (5.85,0) {$\cdots$};
  \foreach \j [evaluate=\j as \n using int(\j+1)] in {0,...,4}
    {\draw[->,thick] (s\j) to[bend left=35] node[above,font=\scriptsize] {$r+\j p$} (s\n);}
  \draw[->,thick] (s5) to[bend left=35] node[above,font=\scriptsize] {$r+5p$} (5.75,0.15);
  \draw[->,red!70!black,thick,dashed] (s4) to[bend left=30] (s1);
  \draw[->,red!70!black,thick,dashed] (s4) to[bend left=35] node[below,font=\scriptsize,pos=0.5] {rate $1$: $4\mapsto\ber(r)+\bin(4,p)$} (s0);
  \draw[->,red!70!black,thick,dashed] (s4) to[bend left=20] (s2);
  \draw[->,red!70!black,thick,dashed] (s4) to[bend left=20] (s3);
   \draw[->,red!70!black,thick,dashed] (s4) to [loop below](s4);
  \draw[->,red!70!black,thick,dashed] (s4) to[bend right=40] (s5);
\end{tikzpicture}
\caption{The reset-and-growth process with generator $\cL$.  From state
$j$ the degree grows to $j+1$ at rate $r+pj$ (solid arrows); at rate one
it is reset to $\ber(r)+\bin(j,p)$ (dashed arrows, drawn from state $4$).
The state seen at reset times is the chain $P$ in~\eqref{eq:XY-chain}.}
\label{fig:reset-growth}
\end{figure}
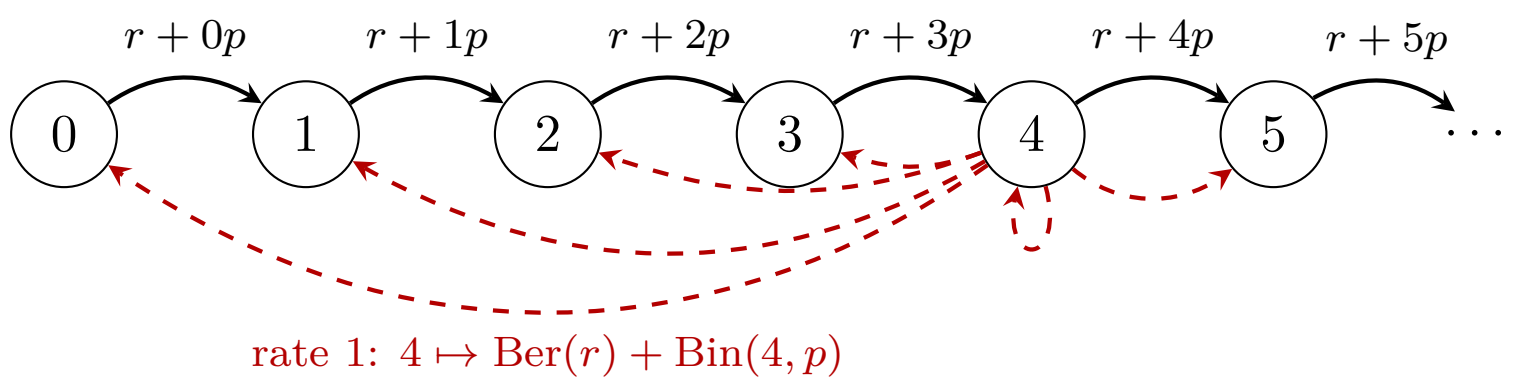

Define a homogeneous Markov chain $\mvY=(Y_n)_{n\geq0}$ on $\dZ_+$ by
\begin{align*}
	(X_{n+1}\mid Y_n=y)\sim\ber(r)+\bin(y,p), \quad
	(Y_{n+1}\mid X_{n+1}=x)\sim H_x.
\end{align*}
Write $P$ for its transition kernel.  When $r<1$, this chain is irreducible and aperiodic on $\dZ_+$: state $0$ has a self-loop, every state can reach $0$, and every state is accessible from $0$.  When $r=1$, the same statements hold on $\dN$, with a self-loop at $1$.
The chain $\mvY$ is the reset-time skeleton of the continuous-time Markov process $\mvZ:=(Z_s)_{s\geq0}$ with generator
\begin{align}\label{eq:degree-generator}
	(\cL\varphi)(j)=\E\varphi\bigl(\ber(r)+\bin(j,p)\bigr)-\varphi(j)+(r+pj)\bigl(\varphi(j+1)-\varphi(j)\bigr).
\end{align}
Indeed, a rate-one clock resets $j$ to $\ber(r)+\bin(j,p)$; between resets, the state grows at rate $r+pj$ when at state $j$. The state immediately before the next reset has the law in~\eqref{eq:XY-chain}.  Because the reset clock has a state-independent exponential rate, a stationary distribution for $P$ is also stationary for $\mvZ$, and conversely.  Figure~\ref{fig:reset-growth} summarizes the two mechanisms.  In terms of the graph, $Y_n$ is to be thought of as the degree of a uniformly sampled vertex in $G_n$; the reset corresponds to the uniformly chosen vertex being replaced by a newborn copy, and the growth corresponds to the vertex gaining edges from later arrivals. The correspondence is made exact in Section~\ref{ssec:ball-limits}.

We analyze the homogeneous Markov chain $(X_k,Y_k)_{k\geq0}$ defined in~\eqref{eq:XY-chain}.
Fix $n\ge x$. Using $1+t/p\le(1+t)^{1/p}$ for $t>-1$ and $0<p\le1$, we obtain
\begin{align*}
   \left(\frac{px+r}{pn+r} \right)^{1/p} 
   \le  \prod_{i=x}^{n-1} \frac{1}{1+1/(pi+r)} 
   &=\pr(H_x \ge n)\\
 &  =\prod_{i=x}^{n-1} \left( 1-\frac{1}{pi+1+r} \right)\le \left(\frac{px+r+1-p}{pn+r+1-p} \right)^{1/p}.
\end{align*}
Consequently, in the stochastic order, we have
\begin{align}\label{eq:H-stochastic-bounds}
    \lfloor(x+r/p)\cdot U^{-p}-r/p\rfloor 
    \preccurlyeq H_{x} 
    \preccurlyeq \lfloor (x+(r+1-p)/p)\cdot U^{-p}-(r+1-p)/p\rfloor,
\end{align}
where $U\sim\unif(0,1)$.

\subsection{Recurrence, transience and the stationary tail}
The following lemmas allow us to apply Foster--Lyapunov drift criteria to $\mvY$. Their proofs follow Lemma~\ref{lem:kappa} below. For $r=1$, recurrence statements refer to the closed class $\dN$, an initial state $0$ enters that class at the first reset. 

\begin{lem}\label{lem:positive-recurrence}
If $p<p_c$, there exist $\ga,A,c>0$ and $\gr\in(0,1)$ such that
\[
\E(1+Y_1^\ga \mid Y_0 = y)
\le \gr \cdot (1+y^\ga)+c\ind_{\{y\le A\}},
\qquad y\ge0.
\]
In particular, the chain $\mvY$ is positive recurrent, aperiodic and converges to its unique stationary distribution $\pi$.
\end{lem}

\begin{lem}\label{lem:transience}
If $p>p_c$, the chain $\mvY$ is transient. Consequently, because its state space is $\dZ_+$, $Y_k\to\infty$ almost surely as $k\to\infty$.
\end{lem}

We also need the following technical lemma.
\begin{lem}\label{lem:kappa}
    Let $p_c$ be the unique solution of $pe^p=1$. For $p\neq p_c$, there exists a unique nonzero $\gk=\gk_{p}$ such that $p^\gk=1-\gk p$. Moreover,
    \begin{enumeraten}
        \item If $p<p_c$, then $\gk>0$ and $p^x<1-px$ for every $x\in(0,\gk)$.
        \item If $p>p_c$, then $\gk<0$ and $p^x<1-px$ for every $x\in(\gk,0)$.
    \end{enumeraten}
\end{lem}
\begin{proof}
Fix $p\in(0,1)$ and define
\[
h_p(x):=p^x-1+px=\exp(x\log p)-1+px,
\qquad x\in\dR.
\]
The function $h_p$ is strictly convex, $h_p(0)=0$, and $h_p(x)\to\infty$ as $x\to\pm\infty$. Hence it has at most one nonzero root. Furthermore,
$
h_p'(0)=\log p+p=\log(pe^p).
$
If $p<p_c$, $h_p'(0)<0$ and strict convexity yields a unique root $\gk>0$ and $h_p(x)<0$ for $x\in(0,\gk)$. If $p>p_c$, $h_p'(0)>0$, and there is a unique root $\gk<0$ and $h_p(x)<0$ for $x\in(\gk,0)$.
\end{proof}

Before proving Lemma~\ref{lem:positive-recurrence}, we record a standard sufficient drift criterion from~\cite{MeyneTweedie_MC_book}. For the statement below, for a time-inhomogeneous Markov chain, a \emph{petite} set is one on which a mixture of transition kernels has a common nonzero lower bound; every finite set in an irreducible countable-state chain is petite (which is the case for us).

\begin{thm}[{\cite[Theorem 15.0.1]{MeyneTweedie_MC_book}}]\label{thm:meyn-tweedie-recurrence}
    Suppose that $(Y_k)_{k\ge0}$ is an irreducible, aperiodic Markov chain on $\cX$. The chain is positive recurrent if there exist a petite set $C\in\cB(\cX)$, finite constants $b,c_1>0$, and a function $V:\cX\to[1,\infty]$ that is finite at some point and satisfies
    \[
    \gD V(y):=\E\bigl(V(Y_1)-V(Y_0)\mid Y_0=y\bigr)
    \le-c_1V(y)+b\ind_C(y),
    \qquad y\in\cX.
    \]
\end{thm}

\begin{proof}[Proof of Lemma~\ref{lem:positive-recurrence}]
For $p<p_c$, Lemma~\ref{lem:kappa} allows us to choose $\ga\in(0,1]$ such that $p^\ga<1-\ga p$. Let $\E_y(\cdot):=\E(\cdot\mid Y_0=y)$ and $U\sim\unif(0,1)$. By~\eqref{eq:H-stochastic-bounds},
\begin{align*}
    \E_y(1+Y_1^\ga)
    &\le 1+\E_y\left((X_1^\ga+(r+1-p)^\ga/p^\ga)\cdot U^{-p\ga}\right)\\
    &\le 1+\frac{p^\ga y^\ga+r^\ga+(r+1-p)^\ga/p^\ga}{1-p\ga}
    =\gr_0(1+y^\ga)+c,
\end{align*}
where $\gr_0:=p^\ga/(1-p\ga)<1$ and
\[
c:=\frac{r^\ga+(r+1-p)^\ga/p^\ga}{1-p\ga}+1-\gr_{0}>0.
\]
Choose $\gr\in(\gr_0,1)$ and $A>0$ large so that $c\le(\gr-\gr_0)A^\ga$. Then
\[
\E_y(1+Y_1^\ga)\le\gr (1+y^\ga)+c\ind_{\{y\le A\}},
\]
which proves the lemma. Theorem~\ref{thm:meyn-tweedie-recurrence}, applied to $V(y)=1+y^\ga$ and the finite set $C=\{y\in\dZ_+:y\le A\}$, shows that the chain is positive recurrent. Irreducibility then gives a unique stationary distribution, and aperiodicity gives convergence to it.
\end{proof}

We use the following criterion, specialized to a discrete-state-space chain, to prove transience when $p>p_c$.
\begin{thm}[{\cite[Theorem 8.0.2]{MeyneTweedie_MC_book}}]\label{thm:meyn-tweedie-transience}
    Suppose that $(Y_k)_{k\ge0}$ is irreducible on a discrete state space $\cX$. The chain is transient if and only if there exist a bounded nonnegative function $V$ and a nonempty set $C\subseteq\cX$ such that
    \[
    \gD V(x)\ge0,
    \qquad x\in C^c,
    \]
    and $D:=\{x:V(x)>\sup_{y\in C}V(y)\}$ is nonempty.
\end{thm}

\begin{proof}[Proof of Lemma~\ref{lem:transience}]
    Fix $p>p_c$. By Lemma~\ref{lem:kappa}, we can choose $\ga\in(0,1]$ such that $p^{-\ga}<1+\ga p$. Set
    \[
    V(x):=1-(1+x)^{-\ga}
    \]
    and let $\E_y(\cdot):=\E(\cdot\mid Y_0=y)$. By~\eqref{eq:H-stochastic-bounds}, we may couple the variables so that
    \[
    1+Y_1\ge\max\{X_1U^{-p},1\},
    \]
    where $U\sim\unif(0,1)$ is independent of $X_1$. Hence
    \[
    (1+Y_1)^{-\ga}\le\min\{X_1^{-\ga}U^{\ga p},1\},
    \]
    with the convention $0^{-\ga}=\infty$. For any $\eps\in(0,1)$ and $y\ge1$,
     \begin{align*}
         \E_y\left((1+Y_1)^{-\ga}\right)
         &\le \E_y\left(\min\{X_1^{-\ga}U^{\ga p},1\}\right)
         \le \frac{((1-\eps)p)^{-\ga}}{1+p\ga}\,y^{-\ga} + \pr_y(X_1\le (1-\eps)py).
     \end{align*}
     Given $Y_0=y$, the variable $X_1$ stochastically dominates $\bin(y,p)$. Therefore, the Chernoff bound gives
     \begin{align*}
         \pr_y(X_1\le (1-\eps)py)
         \le \pr(\bin(y,p)\le (1-\eps)py)
         \le \exp\left(-\eps^2py/2\right).
     \end{align*}
      It follows that
    \[
         \E_y(1+Y_1)^{-\ga}\le h_\eps(y)(1+y)^{-\ga},
    \]
    where
    \begin{equation*}
        h_\eps(y):=\frac{(1-\eps)^{-\ga}p^{-\ga}}{1+\ga p}(1+1/y)^{\ga}+(1+y)^\ga e^{-\eps^2py/2}.
    \end{equation*}
    Choose $\eps>0$ small enough that $(1-\eps)^{-\ga}p^{-\ga}<1+\ga p$, and then choose an integer $R$ large enough that $\gr:=\sup_{y>R}h_\eps(y)<1$. With $C:=\{0,1,\ldots,R\}$, for every $y\in C^c$ we have
    \begin{align*}
        \gD V(y)
        = (1+y)^{-\ga} - \E_y(1+Y_{1})^{-\ga}
        \ge (1-\gr) \cdot (1+y)^{-\ga}.
    \end{align*}
    Because $V$ is strictly increasing and bounded, $D:=\{x:V(x)>\sup_{y\in C}V(y)\}=C^c$. Theorem~\ref{thm:meyn-tweedie-transience} therefore implies that $(Y_k)_{k\ge0}$ is transient. Finally, a transient chain on $\dZ_+$ visits each finite set only finitely often, so $Y_k\to\infty$ almost surely.
\end{proof}

\subsection{Moments of the stationary law}

From now on, assume that $p<p_c$, where $p_c$ is the unique solution of $pe^p=1$. Following~\cite{DiaconisFreedman98}, we realize a transition as a random map $\Psi(\cdot,\theta)$ on $\dZ_+$. Here the innovation space is
\[
 \cE=[0,1]^{\dZ_+}\times(0,1)\times[0,1]^{\dZ_+},
\]
with its product Borel sigma-field and product uniform law. Its three coordinates are the reset uniforms $(U_j)_{j\ge0}$, the growth uniform $U$, and the correction uniforms $(U'_j)_{j\ge0}$ constructed below. Independent innovations $\theta_k$ give $Y_{k+1}=\Psi(Y_k,\theta_{k+1})$. We suppress the subscript on $\theta$ or $\theta$ itself when only one transition is under consideration. If $Y\sim\pi$ is independent of the random map $\Psi$, then stationarity implies that 
\[Y\equald\Psi(Y,\theta).\] 
First, we couple $\Psi(Y,\theta)$ with $pU^{-p}Y$, where $U\sim\unif(0,1)$ is independent of $Y$ but is  a component of the innovation $\theta$. Then, we use Goldie's implicit renewal theorem~\cite{goldie91} to study the tail of the stationary distribution.  In the following general theorem, $M$ is coupled with the innovation map and independent of $Y$. Our application takes $M=pU^{-p}$ using the growth uniform in that innovation.

\begin{thm}[{\cite[Corollary 2.4]{goldie91}}]\label{thm:goldie}
Let $Y\ge0$ satisfy $Y\equald\Psi(Y)$, where $Y$ is independent of the pair $(\Psi,M)$ and $M>0$ almost surely. Suppose that, for some $\gk>0$,
\begin{enumeraten}
\item $\E M^\gk=1$;
\item $\E(M^\gk\log^+M)<\infty$;
\item the distribution of $\log M$ is nonarithmetic.
\end{enumeraten}
Let $m:=\E(M^\gk\log M)$. If
\[
\E\abs{\Psi(Y)^\gk-(MY)^\gk}<\infty,
\]
then $m\in(0,\infty)$ and
\[
\lim_{x\to\infty}x^\gk\pr(Y>x)=C,
\qquad
C:=\frac{1}{\gk m}\E\left(\Psi(Y)^\gk-(MY)^\gk\right)\ge0.
\]
If $C>0$, then $\pr(Y>x)\sim Cx^{-\gk}$.
\end{thm}

To apply this theorem, we construct all relevant variables on a common probability space using~\eqref{eq:H-stochastic-bounds}. We first establish a moment bound.
\begin{lem}[Moment bounds]\label{lem:moment-bounds}
Let $p<p_c$ and $\ga\in(0,\gk)$. There exist constants $C_\ga>0$ and $\gr_\ga\in(0,1)$ such that
\[
\E Y_t^\ga\le\gr_\ga^t\E Y_0^\ga+\frac{C_\ga}{1-\gr_\ga},
\qquad t\ge0.
\]
In particular, if $Y\sim\pi$, then $\E Y^\ga<\infty$ for every $\ga\in(0,\gk)$.
\end{lem}
\begin{proof}
    The upper bound in~\eqref{eq:H-stochastic-bounds} implies the convenient, slightly weaker estimate
    \[
    (Y_1\mid X_1=x)\preccurlyeq\left(x+\frac{r+1}{p}\right)U^{-p},
    \]
    where $U\sim\unif(0,1)$. Consequently, for $\ga<1/p$,
    \begin{align*}
        \E(Y_1^\ga\mid Y_0=y)
        \le\frac{1}{1-\ga p}
        \E\left(\left(X_1+\frac{r+1}{p}\right)^\ga\middle|Y_0=y\right).
    \end{align*}
    First suppose that $\ga\in(0,1]\cap(0,\gk)$. Jensen's inequality and subadditivity gives
    \begin{align*}
        \E\left(\left(X_1+\frac{r+1}{p}\right)^\ga\middle|Y_0=y\right)
        &\le\left(\E(X_1\mid Y_0=y)+\frac{r+1}{p}\right)^\ga\\
        &=\left(py+\frac{1+r+rp}{p}\right)^\ga
        \le p^\ga y^\ga+\left(\frac{1+r+rp}{p}\right)^\ga.
    \end{align*}
    With
    \[
    C_\ga:=\frac{(1+r+rp)^\ga}{p^\ga(1-\ga p)},
    \qquad
    \gr_\ga:=\frac{p^\ga}{1-\ga p}<1,
    \]
    where the last inequality follows from Lemma~\ref{lem:kappa}, we obtain
    \[
    \E Y_1^\ga\le\gr_\ga\E Y_0^\ga+C_\ga.
    \]
    Iteration yields the asserted bound.

    Now suppose that $\ga\in(1,\gk)$. A standard centered-binomial moment bound gives a constant $c_\ga>0$ such that
    \[
    \E\bigl(\abs{X_1-py}^\ga\mid Y_0=y\bigr)
    \le c_\ga(1+y^{\ga/2}).
    \]
    For $\eps>0$, $a>1$, and $x,y\ge0$,
    \[
    (x+y)^a\le(1+\eps)^{a-1}x^a+(1+1/\eps)^{a-1}y^a.
    \]
    It follows that, for some $C_{\ga,\eps}>0$,
    \begin{align*}
        \E\left(\left(X_1+(r+1)/{p}\right)^\ga\mid Y_0=y\right)
        &\le(1+\eps)^{\ga-1}p^\ga y^\ga+C_{\ga,\eps}(1+y^{\ga/2})\\
        &\le(1+\eps)^{\ga-1}p^\ga y^\ga
        +C_{\ga,\eps}\left(1+\eps'y^\ga+\frac{1}{4\eps'}\right).
    \end{align*}
    Because $p^\ga/(1-\ga p)<1$, choose $\eps>0$ and then $\eps'>0$ small enough that
    \[
    \gr_\ga:=\frac{(1+\eps)^{\ga-1}p^\ga+C_{\ga,\eps}\eps'}{1-\ga p}<1.
    \]
    Thus, after enlarging $C_\ga$ if necessary,
    \[
    \E(Y_1^\ga\mid Y_0=y)\le\gr_\ga y^\ga+C_\ga,
    \]
    and iteration again proves the first assertion.

    For the stationary moment claim, start the chain from $Y_0=0$ and let $g_L(x):=\min\{|x|^\ga,L\}$. Since $g_L$ is bounded and the law of $Y_t$ converges to $\pi$,
    \[
    \E g_L(Y)=\lim_{t\to\infty}\E g_L(Y_t)
    \le\limsup_{t\to\infty}\E Y_t^\ga
    \le\frac{C_\ga}{1-\gr_\ga}.
    \]
    Letting $L\to\infty$ and applying monotone convergence proves the claim.
\end{proof}

\subsection{The uniform--threshold coupling}\label{ssec:coupling}

Let $Y\sim\pi$, independently of all the uniform random variables introduced below. We construct a random map $\Psi=\Psi(\cdot,\theta)$ for which $\Psi(Y)\equald Y$ and couple it with the lower bound in~\eqref{eq:H-stochastic-bounds}.
\begin{enumeraten}
    \item Let $(U_j)_{j\ge0}$ be i.i.d.~$\unif(0,1)$. Given $Y=y$, define
    \[
    X:=\ind_{\{U_0\le r\}}+\sum_{j=1}^y\ind_{\{U_j\le p\}}.
    \]
    Then $(X\mid Y=y)\sim\ber(r)+\bin(y,p)$.
    \item Let $U\sim\unif(0,1)$ be independent of $(U_j)_{j\ge0}$. Given $X=x$, define
    \begin{align}\label{eq:Z-coupling}
    Z:=\lfloor(x+r/p)U^{-p}-r/p\rfloor\ge x.
    \end{align}
    For $k\ge x$,
    \begin{align*}
        \pr(Z\ge k\mid X=x)
        &=\pr\left(U^p\le\frac{px+r}{pk+r}\right)
        =\left(\frac{px+r}{pk+r}\right)^{1/p}
        =\prod_{i=x}^{k-1}a_i,
    \end{align*}
    where
    \[
    a_i:=\left(\frac{pi+r}{pi+p+r}\right)^{1/p},
    \qquad i\ge0.
    \]
    Set $b_i:=(r+pi)/(1+r+pi)$. For $i\ge0$, we have
    \begin{align}\label{eq:ai-bi}
    a_i=\frac{1}{(1+p/(pi+r))^{1/p}}
    \le\frac{1}{1+1/(pi+r)}
    =b_i.
    \end{align} 
    \item Let $(U'_j)_{j\ge0}$ be another independent i.i.d.~sequence of $\unif(0,1)$ random variables. Since $0\le a_i\le b_i<1$, both $(1-b_i)/(1-a_i)$ and $(b_i-a_i)/(1-a_i)$ lie in $[0,1]$. Given $Z=k$, define
    \begin{align}\label{eq:coupling-variable}
        W:=k+\ind_{\{U'_0\le(b_k-a_k)/(1-a_k)\}}
        \left(1+\sum_{j\ge1}\prod_{t=1}^j\ind_{\{U'_t\le b_{k+t}\}}\right).
    \end{align}
    By construction, $W\ge Z$ almost surely. Moreover,
    \[
    \pr(W=k\mid Z=k)=\frac{1-b_k}{1-a_k},
    \]
    and, for $m>k$,
    \[
    \pr(W=m\mid Z=k)
    =\frac{b_k-a_k}{1-a_k}\prod_{j=k+1}^{m-1}b_j\cdot (1-b_m).
    \]
    Moreover,
    \[
    (W\mid Z=k, W\neq k) \equald H_{k+1}.
    \]
\end{enumeraten}
The next lemma verifies that $W$ realizes one transition of the original chain.
\begin{lem}\label{lem:coupling}
Let $W$ be defined by~\eqref{eq:coupling-variable}. Then
\[
(W\mid X=x)\equald H_x
\]
and
\[
W\ge\lfloor(X+r/p)\cdot U^{-p}-r/p\rfloor \text{ almost surely}.
\]
Thus we may define $\Psi(Y,\theta):=W$; this map has the transition law in~\eqref{eq:XY-chain}.
\end{lem}
\begin{proof}
    Since
    \[
    \pr(Z=j\mid X=x)=(1-a_j)\prod_{i=x}^{j-1}a_i,
    \qquad j\ge x,
    \]
    the definition of $W$ gives, for $k\ge x$,
    \begin{align*}
        \pr(W=k\mid X=x)
        &=\frac{1-b_k}{1-a_k}\pr(Z=k\mid X=x)\\
        &\quad +(1-b_k)\sum_{j=x}^{k-1}
        \frac{b_j-a_j}{1-a_j}
        \prod_{t=j+1}^{k-1}b_t\,\pr(Z=j\mid X=x)\\
        &=(1-b_k)\left(
        \prod_{t=x}^{k-1}a_t
        +\sum_{j=x}^{k-1}(b_j-a_j)
        \prod_{t=x}^{j-1}a_t
        \prod_{t=j+1}^{k-1}b_t
        \right)\\
        &=(1-b_k)\prod_{t=x}^{k-1}b_t
        =\pr(H_x=k).
    \end{align*}
    The penultimate equality is the elementary telescoping identity obtained by replacing the factors $a_t$ by $b_t$ one at a time. The almost-sure inequality follows  from $W\ge Z$ and~\eqref{eq:Z-coupling}.
\end{proof}

\begin{prop}\label{prop:stationary-goldie-condition}
Let $Y\sim\pi$ and let $W,U$ be defined by the uniform-threshold coupling in Section~\ref{ssec:coupling}. Then
    \begin{equation*}
         \E\abs{W^\gk-(pU^{-p}Y)^\gk}<\infty,
    \end{equation*}
where $\gk>0$ is the nonzero solution of $p^\gk=1-\gk p$ and the same uniform variable $U$ is used in~\eqref{eq:Z-coupling}.
\end{prop}

\begin{proof}
Let $Z$ be defined by~\eqref{eq:Z-coupling}. By the triangle inequality,
\begin{align}\label{eq:goldie-triangle}
\abs{W^\gk-(pU^{-p}Y)^\gk}
&\le W^\gk-Z^\gk +\abs{Z^\gk-(pU^{-p}Y)^\gk},
\end{align}
because $W\ge Z$ almost surely.
We first bound the coupling correction. Given $Z=z$, either $W=z$, or $W$ has the law of $H_{z+1}$. Hence
\begin{align*}
\E(W^\gk-Z^\gk\mid Z=z)
&=\frac{b_z-a_z}{1-a_z}\cdot 
\E \left(H_{z+1}^\gk-z^\gk\right).
\end{align*}
The upper stochastic bound in~\eqref{eq:H-stochastic-bounds} and the fact that $\gk p<1$ imply
$
\E H_{z+1}^\gk\le C\cdot(z+1)^\gk
$
for a constant $C=C(r,p)$. To control the prefactor, put
\[
q_z:=\frac{1}{z+1+r/p}.
\]
Then $a_z=(1-q_z)^{1/p}$. Bernoulli's inequality and $a_z\le1-q_z$ give
$
a_z\ge1-q_z/p,1-a_z\ge q_z.
$
Moreover,
\[
b_z-a_z
\le b_z-(1-q_z/p)
=\frac{(1-p)q_z^2}{p^2(1+(1/p-1)q_z)}
\le\frac{(1-p)q_z^2}{p^2}.
\]
Consequently,
\[
\frac{b_z-a_z}{1-a_z}\le\frac{(1-p)/p^{2}}{z+1},
\]
and therefore
\[
\E(W^\gk-Z^\gk\mid Z=z)
\le K\cdot (z+1)^{(\gk-1)_+},
\]
where $(u)_+:=\max\{u,0\}$. From~\eqref{eq:Z-coupling} and $X\le Y+1$, we have
\[
Z+1\le K U^{-p}\cdot (Y+1).
\]
Thus
\begin{align}\label{eq:goldie-correction-one}
\E(W^\gk-Z^\gk)
&\le K\E\left((Z+1)^{(\gk-1)_+}\right) 
\le K\E\left(U^{-p(\gk-1)_+}(Y+1)^{(\gk-1)_+}\right)
<\infty,
\end{align}
because $p(\gk-1)_+<1$ and $(\gk-1)_+<\gk$.

For the second term in~\eqref{eq:goldie-triangle}, set $s:=\min\{1,\gk\}$. For nonnegative $u,v$,
\[
\abs{u^\gk-v^\gk}
\le K\abs{u-v}^s(u+v)^{\gk-s}.
\]
The definition of $Z$ gives the deterministic bounds
\begin{align*}
\abs{Z-pU^{-p}Y}
&\le K U^{-p}(1+\abs{X-pY}),\\
Z+pU^{-p}Y
&\le K U^{-p}(1+Y).
\end{align*}
Since $U$ is independent of $(X,Y)$ and $\gk p<1$,
\begin{align*}
\E\abs{Z^\gk-(pU^{-p}Y)^\gk}
&\le K\E\left((1+\abs{X-pY})^s(1+Y)^{\gk-s}\right).
\end{align*}
Conditionally on $Y=y$, we have $X\sim\ber(r)+\bin(y,p)$. Lyapunov's inequality yields
\[
\E(\abs{X-pY}^s\mid Y)
\le r^s+\left({Y}/{4}\right)^{s/2}.
\]
It follows that
\begin{align}\label{eq:goldie-correction-two}
\E\abs{Z^\gk-(pU^{-p}Y)^\gk}
&\le K\E(1+Y)^{\gk-s/2}<\infty,
\end{align}
because $\gk-s/2<\gk$ and Lemma~\ref{lem:moment-bounds} applies. Combining~\eqref{eq:goldie-triangle}--\eqref{eq:goldie-correction-two} proves the claim.
\end{proof}

\subsection{Positivity of the tail constant}\label{sec: proof of C>0 p<pc}

Fix $p\in(0,p_c)$ and $r\in(0,1]$.  Let $\gk>0$ be the positive solution
of $p^\gk=1-p\gk$, and let $\pi$ be the stationary distribution of the
chain in~\eqref{eq:XY-chain}.  We write $\pr_y$ for the law started from
$Y_0=y$ and $\pr_\pi$ for the stationary law.
\begin{lem}\label{lem:positive-tail-constant}
There exists $c=c(p,r)>0$ such that
\[
\liminf_{y\to\infty}y^\gk\pr_\pi(Y_0>y)\geq c.
\]
\end{lem}
\begin{rem}
    As a consequence, we have $\E_\pi Y^\gk=\infty$, complementing Lemma~\ref{lem:moment-bounds}: $\E_\pi Y^\ga<\infty$ for $\ga<\gk$. The proof used only three features of the model: $(X_1\mid Y_0=y)\succcurlyeq\bin(y,p)$; the lower bound $H_x\succcurlyeq\floor{xU^{-p}}$ from~\eqref{eq:H-stochastic-bounds}; and positive recurrence with an accessible atom. 
\end{rem}

\begin{proof}
We follow the classical argument by Kesten~\cite{kesten73}, by reducing the tail probability to the probability that a single excursion goes above $y$ and estimate the probability by an exponential change of measure. 
We consider a countable collection of i.i.d.~$\unif(0,1)$ r.v.s 
\[
(U_k, (U_{k,j}, U'_{k,j})_{j\ge 0})_{k\ge 1}
\]
and extend the uniform-threshold coupling (see~\ref{ssec:coupling}) to the full Markov Chain. Given $Y_0$, we sequentially define for $k=1,2,\ldots$,
\begin{align*}
    X_k &:= \ind_{U_{k,0}\le r} + \sum_{j=1}^{Y_{k-1}} \ind_{U_{k,j}\le p},\\
    Z_k &:=\lfloor (X_k + r/p)\cdot U_k^{-p} - r/p\rfloor,\text{ and }\\
    Y_k &:= Z_k + \ind_{U'_{k,0}\le \gamma(Z_k)}\cdot \left(1+\sum_{j\ge 1} \prod_{t=1}^j \ind_{U'_{k,t}\le b_{t+Z_k}} \right)
\end{align*}
where $\gamma(z):=(b_z-a_z)/(1-a_z)$. It follows from our previous computations that $(Y_k)_{k\ge 0}$ is the MC as defined in~\eqref{eq:XY-chain}. We define $S_0=0$, and for $t\ge 1$,
\begin{align*}
    M_t:=pU_t^{-p},\quad  
    S_t=\sum_{i=1}^t \log M_i,  \quad 
    \cF_t=\gs(Y_0, (U_k, \bigl((U_{k,j})_{j\ge 0}, (U'_{k,j})_{j\ge 0}\bigr)_{1\le k\le t}).
\end{align*}
Note that $\E M_t^\gk = 1$ for all $t\ge 1$.

We will work with a simplified collection of random variables defined on the same probability space as follows: $\wY_0=Y_0$, and sequentially define for $k=1,2,\ldots$,
\begin{align*}
    \wX_k :=  \sum_{j=1}^{\wY_{k-1}} \ind_{U_{k,j}\le p} \le X_k,\qquad
    \wY_k=\wZ_k :=\lfloor (\wX_k + r/p)\cdot U_k^{-p} - r/p\rfloor \le Z_k \le Y_k.
\end{align*}

Now, consider the tilted probability measure $\wpr$ defined as
\begin{align*}
    \frac{\dd\wpr}{\dd\!\pr}\bigg|_{\cF_t} = e^{\gk S_t} \text{ for } t=0,1,2,\ldots.
\end{align*}
Under the $\wpr$ measure the distribution of $(U_{k,j},U'_{k,j})_{j\ge 0,k\ge 1}$ are i.i.d.~$\unif(0,1)$ and $(U_k)_{k\ge 1}$ are i.i.d.~with density $p^\gk x^{-p\gk}\ind_{x\in(0,1)}dx$. Define
\[
\gl := p^\gk = 1-p\gk.
\]
Thus, under $\wpr$, $\log M_t\sim \log p + {\mathrm{Exp}(\lambda/p)}$ with mean $\log p +p/\gl=\log p +p^{1-\gk}=m>0$. 

In particular, under $\wpr$, we have $S_t/t\to m$ a.s.~as $t\to\infty$. Thus we have for any $\theta\in (0,m)$ and $u\in (0,\infty)$
\[
\wpr(\inf\{S_t - \theta t\mid t\ge 0\}>-\infty)
=\wpr(\tau_{u+}:=\inf\{t\ge 0\mid S_t>u\}<\infty)=1.
\]
Moreover, by the memory-less property of exponential random variables, the overshoot satisfies
\[
S_{\tau_{u+}}-u\sim{\mathrm{Exp}(\lambda/p)}.
\]

We first transfer the problem to an excursion estimate. For any $y>y_0>1$, by positive recurrence of the MC, we have
\begin{align*}
    \pr_\pi(Y_0>y)
    &= \pi_1\cdot  \E_1\left(\sum_{t=0}^{\tau_1^Y-1} \ind_{Y_t>y}\right)\\
    &\ge \pi_1\pr_1(Y_1=y_0)\cdot \pr_{y_0}(\exists\ n\ge 1 \text{ such that }  Y_n>y, Y_t>1, 0\le t<n)\\
    &\ge \pi_1\pr_1(Y_1=y_0)\cdot \pr_{y_0}(\exists\ n\ge 1 \text{ such that } \wY_n>y, \wY_t>1, 0\le t<n).
\end{align*}
Let
\[
G_y:=\{\exists\ n\ge 1 \text{ such that } \wY_n>y, \wY_t>1, 0\le t<n\}.
\]
Heuristically, $\log\wY_n$ is concentrated around $S_n+\log y_0$ by our coupling, so the first time $\wY_n$ crosses $y$ is roughly the first time $S_n$ crosses $\log(y/y_0)$. We make this argument rigorous below in five steps.

\begin{enumeratea}
    \item Fix $\theta\in (0,m)$. Choose $a=a_\theta>0$ large enough so that $\wpr(A_a)\ge 3/4$, where 
\begin{align*}
    A_a:=\{\inf(S_t-\theta t\mid t\ge 0)\ge -a\}.
\end{align*}

\item Given $y_0$, define $\ell_0=y_0e^{-a}/2$ and
\begin{align*}
    \ell_n &:= \ell_0 e^{\theta n} \text{ and }
    \eps_n := \ell_n^{-1/3} \text{ for }n\ge 0.
\end{align*}
We take $y_0$ big enough so that $p\eps_n\ell_n=p\ell_n^{1/3}\ge p\ell_0^{1/3}\ge 1$ or $\ell_0\ge 1/p^3$ and 
\begin{align*}
    \sum_{n=0}^\infty \eps_n = \ell_0^{-1/3} \sum_{n=0}^\infty e^{-\theta n/3} = \frac{\ell_0^{-1/3}}{1-e^{-\theta/3}} \le 1/4
    \text{ or } \ell_0 \ge 4^3/(1-e^{-\theta/3})^3.
\end{align*}
\item Consider the stopping time
\begin{align*}
    \gs:=\inf\{n\ge 1\mid \wX_{n}< (1-\eps_{n-1})\cdot p\wY_{n-1}\}.
   \end{align*}
   Note that $\wY_0=y_0\ge \ell_0$. When $A_a$ holds and $\gs>n$, by induction, we have 
   \[
   \wY_t \ge \left(1-2\sum_{i=0}^{t-1} \eps_i\right) y_0 e^{S_t} \ge  \ell_t \text{ for all } t=0,1,2,\ldots,n.
   \]
   Here, we used the fact that with $A_a, \wY_n \ge \left(1-2\sum_{i=0}^{n-1} \eps_i\right) y_0 e^{S_n} \ge  \ell_n$ and $\wX_{n+1} \ge (1-\eps_{n})\wY_n$ we have
\begin{align*}
    \wY_{n+1} \ge \wX_{n+1}U_{n+1}^{-p}-1 
    &\ge (1-\eps_n)p \wY_n U_{n+1}^{-p} - \widehat  Y_n U_{n+1}^{-p}/\ell_n\\
    & = (1-\eps_n) \wY_n M_{n+1} - \eps_n\widehat  Y_n  M_{n+1}/(p\eps_n\ell_n)\\
    &\ge (1-2\eps_n)  \widehat  Y_n  M_{n+1}
    \ge \left(1-2\sum_{i=0}^{n} \eps_i\right) y_0 e^{S_{n+1}} \ge \ell_{n+1}.
\end{align*}

   \item Choose $u:=\log (2y/y_0)$, so that when $S_n>u$ we have $\wY_n\ge y_0/2\cdot e^{S_n}>y$. Hence,
    \[
        \{Y_0=y_0, \gs \ge \tau_{u+}=n\}\cap A_a \subseteq \{\wY_n>y, \wY_t>1, 0\le t<n\}.
    \]
    \item  In particular, for any $K>1$, we have
\begin{align*}
    \pr_{y_0}(G_y) 
    &\ge \pr_{y_0}(A_a\cap \{\tau_{u+}<\infty, \gs\ge \tau_{u+}, S_{\tau_{u+}}-u\le \log K\}) \\
&=\wE_{y_0}\left(e^{-\gk S_{\tau_{u+}}}; A_a\cap \{\tau_{u+}<\infty, \gs\ge \tau_{u+}, S_{\tau_{u+}}-u\le \log K\} \right)\\
&\ge  e^{-\gk(u+K)}\cdot  \wpr_{y_0}\left( A_a\cap \{\gs \ge \tau_{u+}, S_{\tau_{u+}}-u\le \log K\} \right)\\
&\ge (2K/y_0)^{-\gk}\cdot y^{-\gk} \cdot \left(\wpr_{y_0}(A_a) - \wpr_{y_0}(S_{\tau_{u+}}-u> \log K) - \sum_{n=1}^\infty \wpr_{y_0}(A_a\cap\{\gs=n\})\right).
\end{align*}
Note that, $\wpr_{y_0}(A_a)\ge 3/4, \wpr_{y_0}(S_{\tau_{u+}}-u> \log K) = K^{-\gl/p}<1/4$ for large $K$, and by Chernoff's inequality for binomial r.v., we have
\begin{align*}
    \sum_{n=1}^\infty \wpr_{y_0}(A_a\cap\{\gs=n\}) 
    &\le \sum_{n=1}^\infty \wpr_{y_0}(\wX_{n}< (1-\eps_{n-1})\cdot p\wY_{n-1}, \wY_{n-1}\ge \ell_{n-1})\\
    &\le \sum_{n=1}^\infty \exp(-p\eps_{n-1}^2\ell_{n-1}/2) = \sum_{n=0}^\infty \exp(-p\ell_{n}^{1/3}/2)<1/4
\end{align*}
for $y_0$ large. Thus, we can choose $y_0>1$ such that for all $y>y_0$ we have 
\[
\pr_{y_0}(G_y)  \ge \left(\frac{2K}{y_0}\right)^{-\gk}\cdot y^{-\gk}\cdot \left(\frac34 -\frac14-\frac14\right).
\]
\end{enumeratea}
This completes the proof.
\end{proof}
\begin{rem}
Goldie's implicit renewal theorem yields the existence of the limit $C=\lim_yy^\gk\pr(Y>y)$ but not its positivity. In this theory, positivity is always established separately. 

See~\cite{kesten73,grincevicius75,grey94}, and~\cite[\S2.4]{BDM16} for a textbook account. 
For our chain, the effective multiplier $(1+(X_1-pY_0)/{pY_0})\cdot pU^{-p}$ is not independent of $Y_0$ and the recursion is only asymptotically multiplicative. So the Kesten--Grincevi\v{c}ius--Goldie theorem does not apply verbatim. General results for Lipschitz maps are available in~\cite{mirek11} for the continuous state space $\dR^d$.
The direct proof above follows Kesten's original scheme, the level-dependent tolerance $\eps_n$ keeps the exponent equal to $\gk$.
\end{rem}

\subsection{Proof of Theorem~\ref{thm:limit-chain}}\label{ssec:pf-limit-chain}
\begin{enumeraten}
\item Suppose $p<p_c$. Positive recurrence and uniqueness of the stationary distribution follow from Lemma~\ref{lem:positive-recurrence}. Let $Y\sim\pi$, independently of the innovation variables in the uniform-threshold coupling, and set $M:=pU^{-p}$. The equation
\[
\E M^\gk=\frac{p^\gk}{1-p\gk}=1
\]
is equivalent to $p^\gk=1-p\gk$. Furthermore, $\E(M^\gk\log^+M)<\infty$, and $\log M=\log p-p\log U$ has a continuous, hence nonarithmetic, distribution. Finally,
\begin{align*}
m:=\E(M^\gk\log M)
&=p^\gk\E\left(U^{-p\gk}(\log p-p\log U)\right)
=\log p+p^{1-\gk}>0.
\end{align*}
Proposition~\ref{prop:stationary-goldie-condition} verifies the remaining hypothesis of Theorem~\ref{thm:goldie}. Therefore,
\[
\lim_{y\to\infty}y^\gk\pr(Y>y)=C,
\]
where
\[
C=\frac{1}{\gk m}\E\left(W^\gk-(pU^{-p}Y)^\gk\right)\ge0.
\]
Lemma~\ref{lem:positive-tail-constant} proves that $C>0$.
\item If $p>p_c$, the result follows from Lemma~\ref{lem:transience}.\qed
\end{enumeraten}

\section{A general criterion for empirical convergence}\label{sec:gen}
\subsection{Euler approximations}\label{ssec:euler}
The next lemma is the Euler approximation for a general countable-state Markov chain. It is used for the degree and the rooted-ball chains. The applications below use the exact recursions and moment bounds proved in Section~\ref{sec:local}. Those preliminary calculations do not depend on the empirical-convergence criterion of this section. Let $Q$ be the generator of a continuous-time Markov chain on a countable state space $S$, with total jump rate $\gl_{y}<\infty$ out of $y\in S$, and let $(\cP_s)_{s\geq0}$ be its semigroup. Fix $\theta>0$.  Define $h_j:=1/(j+\theta), j\ge 1$ and $K_j:=I+h_jQ$.  We assume that
\begin{enumerate}[label=(\textup{E}\arabic*),ref=\textup{(E\arabic*)}]
\item\label{E1} there exists $V\colon S\to[1,\infty)$ with finite sublevel sets such that $QV\leq c_VV$ for a constant $c_V>0$,
\item\label{E2} for every finite $F\subseteq S$ there is $m_F$ such that, for $m\geq m_F$, the kernels $(K_j)_{j\geq m}$, define a Markov chain started in $F$ at time $m$, that is, $h_j\gl_{y}\leq1$ for every state $y$ that the chain can occupy at time $j$.
\end{enumerate}

We prove the following lemma.
\begin{lem}[General Euler approximation]\label{lem:euler-approximation}
Assume~\ref{E1} and~\ref{E2}.  Let $F\subseteq S$ be finite, $s>0$, and let $m_t\leq t$ be integers with $m_t\to\infty$ and $\gL_t:=\sum_{j=m_t}^{t-1}h_j\to s$.  Then
\[
 \max_{y\in F}\Ptv{K_{m_t}K_{m_t+1}\cdots K_{t-1}(y,\cdot)-\cP_s(y,\cdot)}\to0 .
\]
\end{lem}

The lemma says that the products of the $K_j$ form an Euler scheme for $\partial_su=Qu$ with variable step sizes $h_j$ whose mesh tends to zero. The only difficulty is that $S$ is infinite and the rates are unbounded, and~\ref{E1} handles it by localization.  The proof uses only $m_t\to\infty$ and $\gL_t\to s$, so it applies without change to step sizes $h_j=1/(j+\theta)$ for any fixed $\theta>0$.

\begin{proof}
Let $\bigl(Y^{(t)}_j\bigl)_{m_t\leq j\leq t}$ be the Euler chain and $(Y_u)_{u\geq0}$ the chain with generator $Q$, both started at $y\in F$.  By~\ref{E1}, $e^{-c_Vu}V(Y_u)$ is a nonnegative supermartingale, and by~\ref{E1}--\ref{E2}, so is 
$$
\prod_{i=m_t}^{j-1}(1+c_Vh_i)^{-1}V(Y^{(t)}_j).
$$ 
Doob's maximal inequality gives, for every $N$,
\begin{align}\label{eq:euler-exit-bound}
 \pr_y\Bigl(\sup_{u\leq s}V(Y_u)>N\Bigr)\leq {e^{c_Vs}V(y)}/N,
 \quad
 \pr_y\Bigl(\max_{m_t\leq j\leq t}V(Y^{(t)}_j)>N\Bigr)\leq {e^{c_V\gL_t}V(y)}/N.
\end{align}
In particular the continuous-time chain is nonexplosive.  Given $\eps>0$, choose $N$ so large that both right sides are at most $\eps$ for all $y\in F$ and all $t$, using $\sup_t\gL_t<\infty$.

Let $A_N:=\{V\leq N\}$, which is  a finite set. We add a cemetery state $\partial$.  Let $Q^N$ be the generator on $A_N\cup\{\partial\}$ obtained from $Q$ by sending every jump out of $A_N$ to $\partial$ and making $\partial$ absorbing.  It is a finite matrix with operator norm $\norm{Q^N}\leq2\max_{A_N}\gl$ on bounded functions.  Let $K^N_j:=\mathrm{Id}+h_jQ^N$, which is a stochastic matrix once $h_j\max_{A_N}\gl\leq1$, hence for all $j\geq m_t$ when $t$ is large.  The chain with kernels $K^N_j$ coincides with the Euler chain until the first exit from $A_N$, and the chain with generator $Q^N$ coincides with $Y$ until its
first exit from $A_N$.  Hence, for every $f$ with $\norm{f}_\infty\leq1$ and $f(\partial):=0$, and every $y\in F$,~\eqref{eq:euler-exit-bound} gives
\[
 \bigl|K_{m_t}\cdots K_{t-1}f(y)-K^N_{m_t}\cdots K^N_{t-1}f(y)\bigr|\leq2\eps,
 \qquad
 \bigl|\cP_sf(y)-e^{sQ^N}f(y)\bigr|\leq2\eps .
\]

On $A_N\cup\{\partial\}$ the matrices $K^N_j$ and $e^{hQ^N}$ are stochastic, hence of operator norm one, and a telescoping sum gives
\begin{align*}
 \norm{\prod_{j=m_t}^{t-1}K^N_j-\prod_{j=m_t}^{t-1}e^{h_jQ^N}}
 &\leq\sum_{j=m_t}^{t-1}\norm{I+h_jQ^N-e^{h_jQ^N}}\\
 &\leq\sum_{j=m_t}^{t-1}h_j^2 \norm{Q^N}^2 e^{h_j\norm{Q^N}}
 \leq e^{\norm{Q^N}}\norm{Q^N}^2\gL_t\max_{j\geq m_t}h_j ,
\end{align*}
where we used $\norm{e^{X}-I-X} \leq \norm{X}^2e^{\norm{X}}$.  The right side tends to zero because $\max_{j\geq m_t}h_j=1/(m_t+\theta)\to0$.  Moreover
\[
\prod_je^{h_j Q^N}=e^{\gL_t Q^N}\text{ and }\norm{e^{\gL_tQ^N}-e^{sQ^N} }\leq \abs{\gL_t-s}\cdot \norm{Q^N} \to0.
\]
Combining for all large $t$ and all $y\in F$, we have $\sup_{\norm{f}_\infty\leq1} \abs{K_{m_t}\cdots K_{t-1}f(y)-\cP_sf(y)} \leq 4\eps+o(1)$. This completes the proof.
\end{proof}

\subsection{Empirical convergence}\label{ssec:empirical}
We now prove a general statement that covers both the empirical degree distribution and the empirical distribution of rooted balls.  Let $S$, $Q$, $\gl$, $V$ and $\cP_s$ be as in Section~\ref{ssec:euler}, satisfying \ref{E1} and~\ref{E2}. Define the finite set 
\[
E_K:=\{y\in S:V(y)\leq K\}.
\]
Let
$(\gn_n)_{n\geq1}$ be random probability measures on $S$ with $\gn_n$ being $\cF_n$-measurable, and define the martingale difference sequence
\[
\eta_{n+1}(F):=\gn_{n+1}F-\E(\gn_{n+1}F\mid\cF_n).
\]
 We assume the following.
\begin{enumerate}[label=(\textup{C}\arabic*),ref=\textup{(C\arabic*)}]
\item\label{C1} \emph{recursion:} $\E(\gn_{n+1}F\mid\cF_n)=\gn_n (I+h_nQ)F$ for every bounded $F$,
\item\label{C2} \emph{tightness:} $\gd_K:=\sup_{n\geq1}\E(\gn_n(E_K^c))\to0$ as $K\to\infty$,
\item\label{C3}  \emph{noise bound:} for every $K$ there is $C_K\in (0,\infty)$ such that, for every $F$ supported on $E_K$, $$\E\abs{\eta_{n+1}(F)}^2 \leq C_K\norm{F}_\infty^2h_n^2m_2(n),$$ where $\sum_{j\geq s}h_j^2m_2(j)\to0$ as $s\to\infty$,
\item\label{C4}  \emph{ergodicity:} the chain with generator $Q$ has a unique stationary law $\Pi$, and $\cP_s(y,\cdot)\to\Pi$ in total variation for every $y$, and 
\item\label{C5}  \emph{entrance rate:} for every $K$ and $y\notin E_K$, the total rate of jumps of $Q$ from $y$ into $E_K$ is at most $1$.
\end{enumerate}

\begin{prop}[Empirical convergence]\label{prop:empirical-convergence}
Assume~\ref{E1}--\ref{E2},~\ref{C1}--\ref{C5}. Then 
\[
\E \Ptv{\gn_n-\Pi} \to0\text{ as }n\to\infty.
\]
\end{prop}

The idea is simple.  By~\ref{C1}, $\gn_n$ follows the Euler approximation scheme up to a martingale noise whose variance is controlled by~\ref{C3}, and the scheme converges to $\cP_s$ by Lemma~\ref{lem:euler-approximation}.  Running the scheme for a long logarithmic time $s$ brings any initial measure close to $\Pi$, by~\ref{C4}, and the assumption~\ref{C2} makes everything uniform over the random initial measure $\gn_m$.

\begin{proof}
Fix $K$, a function $F$ supported on $E_K$, and $s>0$. Define $m=m_n:=\lfloor ne^{-s}\rfloor$.  Let $Q_K$ be the restriction of $Q$ to $E_K$ with killing on exit. Note that $Q_K$ is a finite matrix.  Restricting~\ref{C1} to sets $A\subseteq E_K$ gives, for $j\geq1$,
\begin{align}\label{eq:truncated-euler}
 \gn_{j+1}\big|_{E_K} =\gn_j\big|_{E_K}\bigl(I+h_jQ_K\bigr)+h_jR_j+\eta_{j+1}\big|_{E_K},
 \text{ where }
 R_j(A):=\sum_{B\notin E_K}\gn_j(B)Q(B,A),
\end{align}
 and $R_j$ is the mass entering $E_K$ from outside.  By~\ref{C5}, $R_j\geq0$ and $\norm{R_j}_1\leq\gn_j(E_K^c)$.  For $m$ large the matrices $(I+h_jQ_K)_{j\geq m}$, are substochastic. So the functions $F_j:=(I+h_jQ_K)\cdots(I+h_{n-1}Q_K)F$, $m\leq j\leq n$, with $F_n:=F$, satisfy $\norm{F_j}_\infty\leq\norm{F}_\infty$.  Testing \eqref{eq:truncated-euler} against $F_{j+1}$ and summing over
$m\leq j<n$, we get
\begin{align}\label{eq:main-rec}
 \gn_nF=\gn_mF_m+\sum_{j=m}^{n-1}h_jR_jF_{j+1}+\sum_{j=m}^{n-1}\eta_{j+1}(F_{j+1}).
\end{align}
The second term is at most $\norm{F}_\infty\sum_{j=m}^{n-1}h_j\gn_j(E_K^c)$ in absolute value, whose expectation is at most $\norm{F}_\infty\gd_K \sum_{j=m}^{n-1}h_j$ by~\ref{C2}.  The third term is a sum of martingale differences. Since $F_{j+1}$ is a deterministic function bounded by $\norm{F}_{\infty}$, by~\ref{C3} its second moment is at most 
$$
C_K\norm{F}_\infty^2\sum_{j\geq m}h_j^2m_2(j)\to0.
$$
Using the proof of Lemma~\ref{lem:euler-approximation}, applied to the finite matrix $Q_K$, we have $\|F_m-e^{sQ_K}F\|_\infty\to0$ as $n\to\infty$.  The killed semigroup $e^{sQ_K}$ differs from $\cP_s$ only on paths that leave $E_K$ before time $s$, so for $J\leq K$,
\[
 \E\bigl|\gn_me^{sQ_K}F-\gn_m\cP_sF\bigr|
 \leq2\norm{F}_\infty\Bigl\{\gd_J+\max_{y\in E_J}\pr_y\Bigl(\sup_{u\leq s}V(Y_u)>K\Bigr)\Bigr\},
\]
and the last probability tends to zero as $K\to\infty$ by \eqref{eq:euler-exit-bound}.  Combining with~\eqref{eq:main-rec}, we have
\[
 \limsup_{n\to\infty}\E\bigl|\gn_nF-\gn_{m_n}\cP_sF\bigr|
 \leq\norm{F}_\infty\Bigl\{s\gd_K+2\gd_J+2\max_{y\in E_J}\pr_y\Bigl(\sup_{u\leq s}V(Y_u)>K\Bigr)\Bigr\},
\]
for every $J\leq K$ with $K$ at least the support level of $F$. By~\ref{C4}, convergence of $\cP_sF(y)\to\Pi F$ is uniform over the finite set $E_J$, so
\[
 \sup_m\E\bigl|\gn_m\cP_sF-\Pi F\bigr|
 \leq\max_{y\in E_J}|\cP_sF(y)-\Pi F|+2\norm{F}_\infty\gd_J .
\]
Given $\eps>0$, choose $J$ with $\gd_J<\eps$, then $s$ so large that the maximum is less than $\eps$, then $K\geq J$ so large that $s\gd_K<\eps$ and the exit probability from $E_{K}$ is less than $\eps$.  Combining we get $\limsup_n\E|\gn_nF-\Pi F|\leq8\norm{F}_\infty\eps$ for every $F$ supported on $E_J$. Hence $\gn_nF\to\Pi F$ in $L^1$ for every finitely supported $F$.

For probability measures $\gn,\Pi$ on $S$ and finite $E\subseteq S$, we have 
$$
\E\Ptv{\gn-\Pi}\leq2\sum_{y\in E}\E|\gn(y)-\Pi(y)|+2\Pi(E^c),
$$ 
as
$\gn(E^c)-\Pi(E^c)=\Pi(E)-\gn(E)$.  We apply this with $E=E_J$, let $n\to\infty$  for the indicators of the points of $E_J$, and finally let $J\to\infty$ to complete the proof. 
\end{proof}

\subsection{Application to empirical degree}
Take $S=\dZ_+, Q=\cL, h_n=(n+1)^{-1}, V(j)=1+j$, and $\nu_n=\edeg_n$. The preceding degree admissibility check gives~\ref{E1}--\ref{E2}, and Lemma~\ref{lem:one_step} gives~\ref{C1}. Choose $0<\alpha<\kappa$. The degree-moment bound gives
\[
 \sup_n\E\edeg_n(\{j:1+j>K\})\le CK^{-\alpha},
\]
which is~\ref{C2}. For bounded $F$, a single graph step changes $\sum_vF(d_n(v))$ by at most $C\norm{F}_\infty(1+d_n(u_{n+1}))$. Centering this increment, and accounting for the normalization, gives
\[
 \E|\eta_{n+1}(F)|^2\le C\|F\|_\infty^2(n+1)^{-2}m_2(n),
 \qquad m_2(n):=\E\edeg_n((1+x)^2).
\]
The moment bounds following~\eqref{eq:quenched-variance-bound} show that $\sum_n(n+1)^{-2}m_2(n)<\infty$ for $p<p_c$, proving~\ref{C3}. Positive recurrence of the degree process proves~\ref{C4}, including convergence from the transient state $0$ when $r=1$. Finally, only resets can decrease $j$, and their rate is one, which proves~\ref{C5}. Proposition~\ref{prop:empirical-convergence} consequently gives $$\E\norm{\edeg_n-\pi}_{\tv}\to0.$$

\subsection{Application to $\ell$-balls}
Take $S=\ballsp_\ell, Q=\sL_\ell, V(B)=1+N(B)=1+|V(B)|$, $\gn_n=\gn_n^{(\ell)}$, and write $(\bB_s^{(\ell)})_{s\geq0}$ for the corresponding continuous-time Markov chain. This chain records the Markovian evolution of the entire $\ell$-neighborhood of a typical vertex, seen as a Markov chain with state space $S$. Every transition  increases $N$ by at most one and the total rate is at most
$1+N(B)$, so~\ref{E1} holds with $c_V=2$.  A ball at time $j$ has at most $j$ vertices, so $h_j\gl\leq1$ on the reachable states and~\ref{E2} holds. Assumption \ref{C5} holds because growth cannot decrease $N$, so only the rate-one reset can enter $E_K$ from outside. Assumption~\ref{C1} is Lemma~\ref{lem:l-ball-closure}.  The remaining conditions are the following three lemmas.

\begin{lem}[Uniform tightness -- Condition~\ref{C2} for $\ell$-ball]\label{lem:ball-uniform-tightness}
Under the assumptions of Lemma~\ref{lem:ball-regeneration}, we have
\[
\gd_K=\sup_{n\geq1}\E \gn_n^{(\ell)}(E_K^c) \to0\text{ as }K\to\infty.
\]
\end{lem}
\begin{proof}
Embed the graph process in continuous time by giving each vertex an independent rate-one reproduction clock. The continuous time graph $G(s)$ at time $s$ has a Yule population $\cN(s)$, and the sequence of graphs at the jump times is $(G_n)$, independently of the holding times.  The generator of the continuous-time process acts on $\Phi_F(s):=\sum_{v\in G(s)}F(\ball_\ell(G(s),v))$ as $I+\sL_\ell$. Thus, we have 
\begin{align}\label{eq:ball-yule-many-to-one}
 \E\bigl(e^{-s}\Phi_F(s)\bigr)=P_s^{(\ell)}F(B_\circ), 
\end{align}
for bounded $F$, where $B_\circ$ denotes the one-vertex rooted graph $\ball_\ell(G_1,1)$.  Let $A_K(n):=n\gn_n^{(\ell)}(E_K^c)$, which is nondecreasing in $n$ because the ball of an existing root cannot lose vertices.  At $s=\log n$,~\eqref{eq:ball-yule-many-to-one} and the independence of the holding times give
\[
 P_{\log n}^{(\ell)}(B_\circ,E_K^c) =n^{-1}\E(A_K(\cN(\log n))) \geq\pr\bigl(\cN(\log n)\geq n\bigr)\E(\gn_n^{(\ell)}(E_K^c)),
\]
and $\pr(\cN(\log n)\geq n)=(1-1/n)^{n-1}\geq e^{-1}$ for $n\geq2$.  By Lemma~\ref{lem:ball-regeneration} and nonexplosion, the family $\left(P_s^{(\ell)}(B_\circ,\cdot)\right)_{s\geq0}$ is tight, which proves the claim.
\end{proof}

\begin{lem}[Noise bound -- Condition~\ref{C3} for $\ell$-ball]\label{lem:ball-finite-influence}
If $F$ is supported on $E_K$ and $u$ is the parent at step $n+1$, then
\begin{align}\label{eq:ball-finite-influence}
 \Bigl|\sum_vF(\ball_\ell(G_{n+1},v))-\sum_vF(\ball_\ell(G_n,v))\Bigr| \leq C_{\ell,K}\norm{F}_\infty(1+\deg_n(u)),
\end{align}
and consequently $\E\abs{\eta_{n+1}(F)}^2\leq C_{\ell,K}\norm{F}_\infty^2m_2(n)/(n+1)^2$ with $m_2(n):=\E\edeg_n((1+x)^2)$.
\end{lem}
\begin{proof}
An old-root summand can change only if its ball has at most $K$ vertices before or after the step, hence before the step by monotonicity, and its root $o$ is within distance $\ell$ of $u$.  On a shortest path from $u$ to such an $o$, every vertex except possibly $u$ is within distance $(\ell-1)$ of $o$, so all of its neighbors belong to the old ball and its degree is at most $(K-1)$.  Hence the number of such roots is at most
$$
1+\deg_n(u)\sum_{j=0}^{\ell-1}(K-1)^j\leq C_{\ell,K}(1+\deg_n(u)).
$$
Adding the newborn-root summand proves~\eqref{eq:ball-finite-influence}.
\end{proof}

Recall $(\bB_s^{(\ell)})_{s \geq 0}$ is the Markov chain on $S=\ballsp_\ell$ that records the evolution of the $\ell$-neighborhood of a uniform vertex.

\begin{lem}[Ergodicity  -- Condition~\ref{C4} for $\ell$-ball]\label{lem:ball-regeneration} For every $\ell\geq1$, $0<p<p_c$ and $0<r\le1$, the rooted-ball chain has a unique stationary law $\Pi_\ell$, and from every state its law converges to $\Pi_\ell$ in total variation.  A positive recurrent atom $\bullet$ is the one-vertex ball when $r<1$, and the path of length $\ell$ rooted at an endpoint when $r=1$.  Consequently, if $\tau_\bullet^+$ is the first return time to the atom after leaving it,
\begin{align}\label{eq:k-ball-regenerative-law}
 \Pi_\ell F= \frac{\E_\bullet\bigl(\int_0^{\tau_\bullet^+}F(\bB_s^{(\ell)})\,ds\bigr)}
      {\E_\bullet(\tau_\bullet^+)} .
\end{align}
\end{lem}
\begin{proof}
Nonexplosion follows from Lemma~\ref{lem:l-ball-euler}. The root degree $D_s=d(\bB_s^{(\ell)})$ is a Markov chain with generator $\cL$ -- a reset gives $\ber(r)+\bin(D_s,p)$, while ordinary updates increase the degree at total rate $r+pD_s$.

If $r<1$, a connected rooted ball has degree zero precisely when it is the singleton. The degree chain is irreducible and positive recurrent. So it hits zero in finite mean time from every degree. Starting from the singleton, the ball chain leaves it at rate $2r$, to a rooted edge, and returns after an additional mean time $\E_1H_0<\infty$. Thus the singleton is a positive recurrent atom accessible from every state.

For $r=1$, degree zero is inaccessible after the first edge, and a separate argument is needed. First, we choose $0<\alpha\le1$ with $c=1-p^\alpha-\alpha p>0$. Define $W(B):=(1+d(B))^\alpha$. Jensen's inequality and concavity gives
\[
 \sL_\ell W\le-cW+C,\qquad
 \cP_1^{(\ell)}W\le aW+b,\quad a=e^{-c}<1.
\]
Choose $m$ large enough so that $aW+b\le a'W$ when $d(B)>m$, for some $a'<1$. For the chain sampled at integer times, the hitting time $H_C$ of
$C=C_m:=\{B\mid d(B)\le m\}$ therefore satisfies
\[
 \E_B H_C\le W(B)/(1-a').
\]
This follows by stopping the drift inequality and using $W\ge1$ outside $C$. Note that the set $C$ need not be finite.

We next show that for the endpoint-rooted path $A_\ell$ of length $\ell$ there exists $\delta=\delta(p,\ell,m)>0$ with
\begin{align}\label{eq:ball-path-minorization}
 \inf_{B\in C}\cP_1^{(\ell)}(B,\{A_\ell\})\ge\delta.
\end{align}
This requires exactly $\ell$ reset-clock rings during $[0,1]$, and requires every copying mark at these resets to fail. Until the first reset, suppress ordinary updates that add an edge at the initial root. Subsequently, protect all the newly created roots by suppressing ordinary updates that add edges at any of them. The first root has degree at most $m$ before its reset. Each subsequent current root has degree one. Each new root is joined only to the previous root, since $r=1$. The protected new vertices therefore form a path and have degrees at most two. At the final reset, the initial root is at distance $\ell$ and all other initial vertices lie outside the radius-$\ell$ ball. Keep suppressing edges at the protected new vertices until time one. The final ball is exactly $A_\ell$.

The reset clock is independent and has rate one. Given $d(B)\le m$, the total rate of forbidden ordinary updates is bounded by $1+pm+\ell(1+2p)$. The rate of adding an edge at a vertex of degree $d$ is $(1+pd)$, and a union bound suffices for several protected vertices. The reset copying conditions have conditional probability at least $(1-p)^{m+\ell-1}$. Thus the specified event has probability at least
\[
 \delta:=\frac{e^{-1}}{\ell!}\cdot e^{-(1+pm+\ell(1+2p))}\cdot (1-p)^{m+\ell-1}>0.
\]
Other updates, including updates in an arbitrarily large outer neighborhood, need not be suppressed. This proves~\eqref{eq:ball-path-minorization}.

For completeness, the return-time bound and this minorization imply finite mean time to $A_\ell$. At each visit to $C$, we can try the next one-unit interval. Its success probability is at least $\delta$. Moreover,
\[
 \sup_{B\in C}\E_B W(\bB_1^{(\ell)})\le a(1+m)^\alpha+b=:D.
\]
On failure, the expected additional time to $C$, multiplied by the failure indicator, is at most $D/(1-a')$. Each trial thus has expected cost at most $1+D/(1-a')$, and the probability of needing $j$ trials is at most $(1-\delta)^{j-1}$. Summing proves a finite mean hitting time from every state and a finite mean sampled return time from $A_\ell$. The sampled chain has a self-loop at $A_\ell$, since no effective update over a unit interval has positive probability. Hence it has one accessible positive recurrent aperiodic class and converges in total variation to its unique invariant law. This law is also invariant for every $\cP_s^{(\ell)}$, by commutation and uniqueness. Contraction of total variation under $\cP_s^{(\ell)}$ then gives continuous-time convergence as well.

In both $r<1$ and $r=1$ cases, the atom has a positive finite holding time and finite mean return time after departure. The standard cycle occupation identity gives~\eqref{eq:k-ball-regenerative-law}, with the cycle including the initial holding time. Accessibility from every initial ball gives uniqueness of the stationary law and total-variation convergence from every state.
\end{proof}

\section{The empirical degree distribution and local limits}\label{sec:local}
In Section~\ref{sec:chain}, we analyzed the limiting reset-and-growth chain $P$ separately without any graph structure. We now show that the annealed law of the degree of a uniformly sampled vertex, is an Euler discretization of the CTMC generated by $\cL$ on the logarithmic scale. We use a backward martingale representation and a monotone coupling estimate to improve the result to the random empirical degree distribution. Finally, we generalize the method to find limiting behavior of rooted balls of any fixed radius. 

\subsection{The empirical degree recursion}

Recall that $\edeg_t$ is the empirical degree distribution. For $t\ge1$ and a function $\varphi:\dZ_+\to\dR$, define
\[
\cZ_\varphi(t):=\frac1t\sum_{v=1}^{t}\varphi\bigl(\deg_t(v)\bigr) = \edeg_t \varphi,
\]
and the operators
\begin{align}
(\cL\varphi)(j) &:=\E\varphi\bigl(B_r+B_{j,p}\bigr)-\varphi(j)+(r+pj) \gD\varphi(j),\label{def:L}\\
K_t &:= I+\cL/(t+1).\label{def:K}
\end{align}
where  $B_r\sim \ber(r)$ and $B_{j,p}\sim \bin(j,p)$ are independent random variables. Here and below, the Bernoulli and binomial random variables in a reset are independent, and $\nb_t(u)$ denotes the set of neighbors of $u$ in $G_t$. 
 We write 
\[
\norm{\varphi}:=\sup_j\abs{\varphi(j)},\qquad \gD\varphi(j):=\varphi(j+1)-\varphi(j),
\]
and $d_t(v):=\deg_t(v)$. Note that $\cL$ is the generator of a CTMC on $\dZ_+$ and $K_t$ is a Markov kernel on the feasible degree states $0\le j< t$ for all $t\ge 0$. Moreover, $K_t$'s commute with each other for all $t\ge 0$. First, we compute the conditional mean under one time step of the duplication-divergence model dynamics.

\begin{lem}\label{lem:one_step}
For every $\varphi:\dZ_+\to\dR$, and every $t\ge1$, we have
$
\E\bigl(\cZ_\varphi(t+1)\mid\cF_t\bigr)
=\cZ_{K_t\varphi}(t).
$
\end{lem}
\begin{proof}
Let $u:=u_{t+1}$ be the parent, $\gx:=\ind_{\{\{t+1,u\}\in\cI_{t+1}\}}\sim\ber(r)$ and, for $w\sim u$ in $G_t$, $\gh_w:=\ind_{\{\{t+1,w\}\in\cI_{t+1}\}}\sim\ber(p)$, all independent given $\Fplus_t$. Then
\[
d_{t+1}(t+1)=\gx+\sum_{w\sim u}\gh_w, \quad 
d_{t+1}(u)=d_t(u)+\gx,\quad
d_{t+1}(w)=d_t(w)+\gh_w 
\text{ for } w\sim u,
\]
and $d_{t+1}(v)=d_t(v)$ otherwise.
Hence
\begin{align*}
    &(t+1)\E\bigl(\cZ_\varphi(t+1)\mid\Fplus_t\bigr)\\
&\qquad =t\cZ_\varphi(t)+\E \varphi\bigl(\ber(r)+\bin(d_t(u),p)\bigr)+r \gD\varphi(d_t(u))+p\sum_{w\sim u}\gD\varphi(d_t(w)).
\end{align*}
Averaging over $u\sim\unif[t]$ and using the identity
$$
\sum_{u\le t}\sum_{w\sim u}\gD\varphi(d_t(w))=\sum_{w\le t}d_t(w)\gD\varphi(d_t(w))
$$
we get
\[
(t+1)\E\bigl(\cZ_\varphi(t+1)\mid\cF_t\bigr)=(t+1)\cZ_\varphi(t) + \cZ_{\cL\varphi}(t) = (t+1)\cZ_{K_t\varphi}(t).
\]
This proves the claim.
\end{proof}

In particular,  $\cL\vone=0$, and $\cL\;\mathrm{id}=(2p-1)\;\mathrm{id}+2r$, since
$\E(\ber(r)+\bin(j,p))=r+pj$. Noting $I_t=\frac12t \cZ_{\mathrm{id}}(t)$ recovers
$\E(I_{t+1}\mid\cF_t)=(1+2p/t)I_t+r$.

\begin{lem}\label{lem: spine}
Let $K_s$ be as defined in~\eqref{def:K}. Then, $K_s$ is a Markov kernel from $\{0,1,\ldots,s-1\}$ to $\{0,1,\ldots,s\}$ and for $0\le j\le s-1$, $K_s(j,\cdot)$ is the probability measure on $\dZ_+$ supported on $\{0,1,\dots,j+1\}$ that puts mass
\[
\frac{s-r-pj}{s+1}\text{ on } j,
\quad
\frac{r+pj}{s+1} \text{ on } j+1,
\quad
\frac{1}{s+1}\pr\bigl(\ber(r)+\bin(j,p)=u \bigr)\text{ on any } u\in \{0,\dots,j+1\}.
\] 
Writing 
\[
\cK_{t,T}:=\prod_{s=t}^{{T-1}}K_{s},\qquad \cK_{T,T}:=\mathrm{Id},
\]
the kernel $\cK_{t,T}$ is the law at time $T$ of the time-inhomogeneous Markov chain $(J_{s})_{t\le s\le T}$  with transitions $K_s$.
\end{lem}

\begin{proof}
The three masses are nonnegative for $0\le j\le s-1$, because $r+pj\le 1+p(s-1)=ps+(1-p)\le s$, and they sum to one. For any $\varphi$, we have
\[
\varphi(j)+\frac1{s+1}(\cL\varphi)(j)=\frac{s-r-pj}{s+1}\varphi(j)+\frac{r+pj}{s+1}\varphi(j+1)+\frac{1}{s+1}\E\varphi\bigl(\ber(r)+\bin(j,p)\bigr).
\]
Finally, if $J_s\le s-1$ then $J_{s+1}\le s$, so the chain is well defined for all $s\ge t$ once $J_t\le t-1$, which holds since degrees in $G_t$ are at most $t-1$.
\end{proof}

\begin{lem}[Degree moment bounds]\label{lem:uniform-empirical-moment}
    For every $\ga>0$, there exists $C=C(p,r,\ga)<\infty$ such that, for all $t\ge1$,
    \[
    \E\edeg_t\bigl((1+x)^\ga\bigr)\le C\cdot 
    \begin{cases} 
    1, & \text{if } p^\ga+\ga p<1,\\
    1+\log t, & \text{if } p^\ga+\ga p=1,\\
    t^{p^\ga+\ga p-1}, & \text{if } p^\ga+\ga p>1.
    \end{cases}
    \]
    In particular, if $p<p_c$, then for every $0<\ga<\gk$, we have
    $
\sup_{t\ge1}\E\edeg_t\bigl((1+x)^\ga\bigr)<\infty.
    $
\end{lem}
\begin{proof}
    We first estimate the growth and reset terms in $\cL$. For every $\ga>0$, Taylor expansion gives
    \[
    (r+pj)\left((j+2)^\ga-(j+1)^\ga\right)
    =\ga p(1+j)^\ga+O\left((1+j)^{\ga-1}\right).
    \]
    Note that for the reset term $(B_j\sim \ber(r)+\bin(j,p))$ we have,
    \begin{align*}
        \E(1+B_j)=1+r+jp\ge p(1+j), \qquad \var(B_j)=r(1-r)+j\,p(1-p)
    \end{align*}
    Since, $0\le \frac{1+B_j}{1+r+pj}\le \frac{j+2}{p(j+1)}\le 2/p$, we have by Taylor's theorem, 
    \begin{align*}
        \abs{\E(1+B_j)^\ga-(1+r+pj)^\ga}&\le C(1+r+pj)^{\ga-2}\var(B_j)\\
        &\le C(1+j)^{\ga-1}.
    \end{align*}
    Here we used 
    \[\abs{z^\ga-1-\ga(z-1)}\le C (z-1)^2, \qquad0\le z\le 2/p.\]
    Combining with the fact that $(1+r+pj)^\ga =p^\ga (1+j)^\ga+\cO\left((1+j)^{\ga-1}\right)$, we get 
    \[
    \cL\left((1+x)^\ga\right)(j)=(p^\ga+\ga p-1)(1+j)^\ga+\cO\left((1+j)^{\ga-1}\right).
    \]
    Taking expectations in Lemma~\ref{lem:one_step}, and using$(1+j)^{\ga-1}\le(1+j)^{(\ga-1)_+}$, we get
    \begin{align*}
    \E\edeg_{t+1}\left((1+x)^\ga\right)&\le\left(1+\frac{p^\ga+\ga p-1}{t+1}\right)\E\edeg_t\left((1+x)^\ga\right)+\frac{C}{t+1}\E\edeg_t\left((1+x)^{(\ga-1)_+}\right).
    \end{align*}
    This iteration gives,
    \begin{align}\label{eq:degree-moment-iteration}
    \E\edeg_t\left((1+x)^\ga\right)\le Ct^{p^\ga+\ga p-1}\left(1+\sum_{s=1}^{t-1}s^{-p^\ga-\ga p}\E\edeg_s\left((1+x)^{(\ga-1)_+}\right)\right).
    \end{align}
    For $0<\ga\le 1$, we get $\E\edeg_s\bigl((1+x)^{(\ga-1)_+}\bigr)=\E\edeg_s(1)=1$, which gives the result.

    Now, let us consider the case $\ga >1$. Assume that the statement has already been proved with $\ga$ replaced by $\ga-1$ (for induction). By convexity of $x\mapsto p^x+xp-1$, we have that $p^{\ga-1}+(\ga-1)p-1<\frac{\ga-1}{\ga}(p^\ga+\ga p-1)$. If $p^\ga+\ga p\le 1$, the second term in \eqref{eq:degree-moment-iteration}, becomes bounded. Thus we get a constant bound when $p^\ga+\ga p<1$, and a bound of order $1+\log t$ when $p^\ga+\ga p=1$. 
    
    Finally, when $p^\ga+\ga p>1$, we have $(p^{\ga-1}+(\ga-1)p-1)_+<p^\ga+\ga p-1$. Thus, induction again gives that the second term in \eqref{eq:degree-moment-iteration} is finite. Hence we get the final bound.
\end{proof}

\subsection{Convergence of the inhomogeneous MC}

\begin{lem}\label{lem:euler-degree-approximation}
Let $(\cP_s)_{s\geq0}$ be the semigroup generated by $\cL$.  For every fixed $s>0$ with $m_t=\lfloor te^{-s}\rfloor$, we have
\[
\sup_{y\in F}\norm{K_{m_t}K_{m_t+1}\cdots K_{t-1}(y,\cdot)
-\cP_s(y,\cdot)}_{\tv}\to 0 \text{ as } t\to\infty,
\]
for any fixed finite set $F\subset \dZ_+$.  More generally, the same conclusion holds for any integers $m_t\leq t$ with
$m_t\to\infty$ and $\sum_{j=m_t}^{t-1}(j+1)^{-1}\to s$.
\end{lem}

This is the degree specialization of Lemma~\ref{lem:euler-approximation}. Take $V(y)=1+y$, for which $\cL V(j)=2r+(2p-1)j\le2V(j)$. For a finite starting set $F$, begin at $m\ge1+\max F$. Then a state at time $s\ge m$ is at most $\max F+s-m\le s-1$, so all Euler transitions are feasible. This verifies both hypotheses of that lemma.

\subsection{Proof of Theorem~\ref{thm:degree-limit}}
\label{ssec:pf-degree-limit}
Fix $p\in (0,p_c)$.
Let $\mu_t:=\E\edeg_t$ be the annealed law of the degree of a uniformly sampled vertex of $G_t$. Taking expectations in Lemma~\ref{lem:one_step}, we get $\mu_{t+1}=\mu_tK_t$. Thus, for every $1\le m\le t$, we have $\mu_t=\mu_m\cK_{m,t}$. By Lemma~\ref{lem:uniform-empirical-moment}, there exists $\ga>0$ such that
$
C:=\sup_{t\ge1}\mu_t\bigl((1+\cdot)^\ga\bigr)<\infty.
$
Consequently, for every integer $N\ge1$, we have
\begin{align}\label{eq:degree-uniform-tightness}
\sup_{t\ge1}\mu_t\bigl([N,\infty)\bigr)
\le CN^{-\ga}.
\end{align}
By Lemma~\ref{lem:positive-recurrence}, the continuous time Markov chain converges to the stationary distribution $\pi$. Hence, by the ergodic theorem, we get
\[
\lim_{s\to\infty}\norm{\cP_s(j,\cdot)-\pi}_{\tv}=0,
\text{ for all } j\in\dZ_+.
\]
In particular, for every fixed $N$, we have
\begin{align}\label{eq:degree-semigroup-convergence}
\lim_{s\to\infty}
\max_{0\le j\le N}\norm{\cP_s(j,\cdot)-\pi}_{\tv}=0.
\end{align}
Fix $s>0$, for all sufficiently large $t$, we have, with $m_t:=\lfloor te^{-s}\rfloor$,
\begin{align*}
    \norm{\mu_t-\pi}_{\tv} =\norm{\mu_{m_t}\cK_{m_t,t}-\pi}_{\tv}\le\norm{\mu_{m_t}\cK_{m_t,t}-\mu_{m_t}\cP_s}_{\tv}+\norm{\mu_{m_t}\cP_s-\pi}_{\tv}.
\end{align*}
Splitting each term in the last inequality above into $j\ge N$, and $j<N$, and using~\eqref{eq:degree-uniform-tightness}, we get
\begin{align*}
    \norm{\mu_t-\pi}_{\tv}\le \sup_{0\le j\le N}\norm{\cK_{m_t,t}(j,\cdot)-\cP_s(j,\cdot)}_{\tv}+\sup_{0\le j\le N}\norm{\cP_s(j,\cdot)-\pi}_{\tv}+2CN^{-\ga}.
\end{align*}
All states $j\le N$, are admissible as long as $m_t\ge N$ for large enough $t.$ Hence, by Lemma~\ref{lem:euler-degree-approximation}, we get,
\begin{align*}
    \limsup_{t\to\infty}\norm{\mu_t-\pi}_{\tv}\le \sup_{0\le j\le N}\norm{\cP_s(j,\cdot)-\pi}_{\tv}+2CN^{-\ga}.
\end{align*}
Letting, $s\to\infty$ and  using~\eqref{eq:degree-semigroup-convergence}, we have
$
\limsup_{t\to\infty}\norm{\mu_t-\pi}_{\tv}
\le 2CN^{-\ga}.
$
Letting $N\to\infty$, we obtain the result that $
\limsup_{t\to\infty}\norm{\mu_t-\pi}_{\tv}=0.
$
Finally, we note that $\E f(\deg_T(U_T))=\mu_T f$ for all bounded function $f:\dZ_+\to\dR$.

Moreover, when $p_c<p<1$, Lemma~\ref{lem:euler-degree-approximation} gives that $\deg_T(U_T)\to \infty$ almost surely, as $T\to\infty$. Let $h(x):= \log(1+x)$, and $D_t:= \deg_t(U_t)$, then we have by~\eqref{def:L}
\[
    (\cL h )(k)= \E \log\left(\frac{1+B_r+B_{k,p}}{1+k}\right)+(r+pk)\log \left(\frac{2+k}{1+k} \right)
\]
this is $\log p + p + O((1+k)^{-1})$. Using the definition of the Markov kernel $K_t$ in~\eqref{def:K}, we get
\begin{align*}
h(D_T)-h(D_1)&= \sum_{t=1}^{T-1} \log \left(\frac{1+D_{j+1}}{1+D_j} \right)= \sum_{t=1}^{T-1}\frac{\cL h(D_t)}{t+1} + M_T
\end{align*}
where $M_T:= \sum_{t=1}^{T-1}\left[\log \left(\frac{1+D_{t+1}}{1+D_t} \right)-\frac{\cL h(D_t)}{t+1}\right]$, $M_1=0$. $M_T$ is a mean 0 martingale and $\E M_T^2\le C \log T$, and hence $M_T/\log T\to 0$, in $L_2$. And since $D_T\to \infty$, we have that, $\log (1+D_T)/\log T \probc \log(pe^p)$, as $T\to\infty.$
\qed


\subsection{The empirical degree distribution}\label{ssec:quenched-degree}

\begin{lem}\label{lem: general martingal bdd phi}
Fix $T>1$ and a bounded $\varphi$. Define
$
g_{t,T}:=\cK_{t,T}\varphi,\text{ for } 1\le t\le T,
$
so that $g_{T,T}=\varphi$ and $\norm{g_{t,T}}\le\norm{\varphi}$. Then
\[
(M_t:=\cZ_{g_{t,T}}(t))_{1\le t\le T},
\]
is an $(\cF_t)$-martingale with $M_T=\cZ_\varphi(T)$.
\end{lem}

\begin{proof}
By Lemma~\ref{lem: spine}, $\cK_{t,T}= K_t\cK_{t+1,T}$, so $g_{t,T}=K_t g_{t+1,T}$. Lemma~\ref{lem:one_step} then gives $\E(M_{t+1}\mid\cF_t)=\cZ_{K_t g_{t+1,T}}(t)=\cZ_{g_{t,T}}(t)=M_t$.
\end{proof}

The next estimate quantifies how fast $\cK_{t,T}$ forgets its initial state.  It is the reason that $p_c$ is the threshold for the law of large numbers.

\begin{lem}[Degree kernel contraction]\label{lem:degree-kernel-contraction}
Suppose $0<p<p_c$, and let $\ga\in(0,\gk\wedge1]$ with $p^\ga+\ga p<1$. Put $c_\ga:=1-p^\ga-\ga p>0$.  For admissible states $j,j'$ and $1\leq t\leq T$,
\begin{align}\label{eq:degree-kernel-tv-contraction-general}
 \norm{\cK_{t,T}(j,\cdot)-\cK_{t,T}(j',\cdot)}_{\tv}
 \leq \,\abs{j-j'}^\ga\left(\frac{t+1}{T+1}\right)^{c_\ga}.
\end{align}
Consequently, for bounded $\varphi$ and $g_{t,T}$ as in Lemma~\ref{lem: general martingal bdd phi}, we have
\begin{align}\label{eq:degree-kernel-difference-bound}
 \left|g_{t,T}(j)-g_{t,T}(j')\right|
 \leq 2\norm{\varphi}_\infty\,\abs{j-j'}^\ga \left(\frac{t+1}{T+1}\right)^{c_\ga},\text{ and }
 \left|\gD g_{t,T}(j)\right|
 \leq 2\norm{\varphi}_\infty \left(\frac{t+1}{T+1}\right)^{c_\ga}
\end{align}
where $\gD g(j):=g(j+1)-g(j)$.
\end{lem}

\begin{proof}
Couple one step of $K_s$ from $j'\leq j$.  With probability $1/(s+1)$, resample both states using the same Bernoulli parent-link and the same Bernoulli sequence for the copied edges.  Otherwise use one uniform variable for the two birth decisions.  The order is preserved, and the difference $D_s$ of the two states is itself a Markov chain. Given $D_s=D$, it becomes $\bin(D,p)$ with probability $1/(s+1)$, becomes $D+1$ with probability $pD/(s+1)$, and is unchanged otherwise.  Once the states meet they remain together, so
\begin{align}\label{eq:coupling-tv-bound}
 \norm{\cK_{t,T}(j,\cdot)-\cK_{t,T}(j',\cdot)}_{\tv}
 \leq\pr(D_T>0\mid D_t=j-j')\leq\E( D_T^\ga \mid D_t=j-j'),
\end{align}
because $D_T$ is a nonnegative integer.  For $\ga\leq1$, Jensen's inequality gives $\E\bin(D,p)^\ga\leq(pD)^\ga$, and concavity gives $(D+1)^\ga-D^\ga\leq\ga D^{\ga-1}$ for $D\geq1$.  Hence
\[
 \E(D_{s+1}^\ga\mid D_s=D)  \leq D^\ga+\frac{(p^\ga-1)D^\ga+p\ga D^\ga}{s+1} =\Bigl(1-\frac{c_\ga}{s+1}\Bigr)D^\ga,
\]
and iteration from $D_t=\,\abs{j-j'}$ yields
\[
 \E(D_T^\ga\mid D_{t}=j-j')\leq\,\abs{j-j'}^\ga\prod_{s=t}^{T-1}\Bigl(1-\frac{c_\ga}{s+1}\Bigr) \leq\,\abs{j-j'}^\ga\Bigl(\frac{t+1}{T+1}\Bigr)^{c_\ga}.
\]
Such an $\ga$ exists by Lemma~\ref{lem:kappa}.  The bounds in~\eqref{eq:degree-kernel-difference-bound} follow from $|\cK_{t,T}\varphi(j)-\cK_{t,T}\varphi(j')| $ $\leq 2\norm{\varphi}_\infty\cdot$ $\norm{\cK_{t,T}(j,\cdot)-\cK_{t,T}(j',\cdot)}_{\tv}.$
\end{proof}

\subsection{Proof of Theorem~\ref{thm:edegree-lln}}\label{ssec:pf-edegree-lln}
Fix a bounded $\varphi$ and $T>1$.  By Lemma~\ref{lem: general martingal bdd phi}, we have
\[
 \cZ_\varphi(T)-\E\cZ_\varphi(T)=\sum_{t=1}^{T-1}\gD M_{t+1},
 \text{ where }
 \gD M_{t+1}:=\cZ_{g_{t+1,T}}(t+1)-\cZ_{g_{t,T}}(t).
\]
Write $g=g_{t+1,T}$. Let $u=u_{t+1}$ be the parent,  $\chi\sim\ber(r)$ indicate the parent edge, and $\eta_w\sim\ber(p)$, $w\in\nb_t(u)$, indicate the copied edges.  The new vertex has degree $\chi+\sum_{w\sim u}\eta_w$, and
\[
 (t+1)\cZ_g(t+1)-t\cZ_g(t)-g(0)
 =g(\chi+\sum_{w\sim u}\eta_w)-g(0)+\chi\,\gD g(\deg_t(u))+\sum_{w\sim u}\eta_w\,\gD g(\deg_t(w)).
\]
Fix $\ga \in (0,1]$ and $c=c_\ga=1-p^\ga-\ga p$ as in Lemma~\ref{lem:degree-kernel-contraction}.  By
\eqref{eq:degree-kernel-difference-bound}, we have, for some $C>0$,
\begin{align*}
    \abs{g(d)-g(0)} &\le C\norm{\varphi}_\infty T^{-c}t^c(1+d)^\ga \text{ and }
    \abs{\gD g(d)} \le C\norm{\varphi}_\infty T^{-c}t^c.
\end{align*}
Hence
$|\gD M_{t+1}|\leq X_{t+1}+\E(X_{t+1}\mid\cF_t)$, where
\[
 X_{t+1}:=\frac1{t+1}\bigl|(t+1)\cZ_g(t+1)-t\cZ_g(t)-g(0)\bigr|
 \leq C\norm{\varphi}_\infty T^{-c}t^{c-1}
 \bigl(1+\deg_t(u_{t+1})\bigr)^{\ga}.
\]
Therefore, by martingale orthogonality and $\E(\gD M_{t+1})^2\leq4\E X_{t+1}^2$, we have
\begin{align}\label{eq:quenched-variance-bound}
    \var\bigl(\cZ_\varphi(T)\bigr)
    &=\sum_{t<T}\E(\gD M_{t+1})^2 \leq C\norm{\varphi}_\infty^2 T^{-2c}\sum_{t<T} t^{2c-2}\cdot m_{2\ga}(t)\notag\\
    &\leq C\norm{\varphi}_\infty^2 T^{2(p^\ga+\ga p)-2}\sum_{t<T} t^{(p^{2\ga}+2\ga p -1)_+ - 2(p^\ga+\ga p)} \log(1+t),
\end{align}
where $m_{2\ga}(t):=\E\frac1t\sum_{v\leq t}\bigl(1+\deg_t(v)\bigr)^{2\ga}$ and by Lemma~\ref{lem:uniform-empirical-moment}, it satisfies
$
m_{2\ga}(t)\leq C t^{(p^{2\ga}+2\ga p -1)_+}(1+\ind_{p^{2\ga}+2\ga p =1}\log t).
$
For $p\le 1/e<\sqrt2-1$, we can take $\ga=1$, to get the upper bound 
\begin{align*}
    \var\bigl(\cZ_\varphi(T)\bigr) &\le 
C\norm{\varphi}_\infty^2 T^{4p-2+(1-4p)_+}(1+\ind_{p=1/4}\log T)\\
&\le C\norm{\varphi}_\infty^2 ( \ind_{p<1/4} T^{-1} + \ind_{p=1/4}T^{-1}\log T +\ind_{1/4<p\le 1/e} T^{4p-2}).
\end{align*}

For $p\in(1/e, p_c)$, we choose $\ga=\ga_*(p):=1+ \log\abs{\log p}/\abs{\log p}$, to maximize $c_\ga=1-p^\ga -\ga p\in (0,1]$. Using  $\gb(p)=p^{\ga}+\ga p = p+p\log\abs{e\log p}/\abs{\log p}$ and $p^{2\ga}+2\ga p\ge 1$, we get the upper bound
\begin{align*}
    \var\bigl(\cZ_\varphi(T)\bigr) &\le C\norm{\varphi}_\infty^2 T^{2\gb(p)-2}.
\end{align*}
Since $\E\cZ_\varphi(T)\to\pi\varphi$ by Theorem~\ref{thm:degree-limit}, this proves that $\cZ_\varphi(T)\to\pi\varphi$ in $L^2$, and $\varphi=\ind_{\{j\}}$ gives $N_j(T)/T\to\pi_j$. For the total-variation convergence, fix $j_0$ and use 
\[
\norm{\edeg_T-\pi}_{\tv}\leq\sum_{j\leq j_0}|N_j(T)/T-\pi_j|+\pi(\{j\mid j>j_0\}),
\]
which follows from $\sum_j(\pi_j-N_j(T)/T)=0$.  The finite sum tends to zero in probability, and $\pi(\{j\mid j>j_0\})$ is arbitrarily small for $j_0$ large.\qed

\subsection{Limits of Rooted $\ell$-balls}\label{ssec:ball-limits}
Assume $0<p<p_c$, and $0<r\le 1$ throughout the section. The degree law of Theorem~\ref{thm:degree-limit} does not record degree correlations across edges or wedges, and those correlations are needed for fluctuation control.  We therefore retain the whole rooted neighborhood of a given radius and prove Theorem~\ref{thm:l-ball-limit}.  The argument is the same for every radius, so we carry a general $\ell$ throughout. 
Fix an integer $\ell\geq1$.  For a finite graph $G$ and a vertex $o$, let
\[
 \ball_\ell(G,o):=G[\{v:\dist_G(o,v)\leq\ell\}]
\]
be the induced radius-$\ell$ ball rooted at $o$, and let
$\ballsp_\ell$ be the countable set of rooted isomorphism classes of finite
rooted graphs of radius at most $\ell$.  Recall from
\eqref{def:empirical-l-ball} the random empirical measure
\[
\gn_t^{(\ell)}=\frac1t\sum_{o\leq t}\gd_{\ball_\ell(G_t,o)}.
\]
If $U_t$ is a uniformly chosen vertex of $G_t$, let
\[
\mu^{(\ell)}_t:= \mathrm{L}(\ball_\ell(G_t,U_t))=\E\gn_t^{(\ell)}
\]
be the annealed law of the empirical $\ell$-ball measure.

We first define the dynamics of these balls. For $B=(V,E,o)\in\ballsp_\ell$ and $j\le\ell$, write $\tau_jB:=B_j(B,o)$. We retain the successive balls through the vector
\[
    \sB^{(\ell)}(G,o):=(B_1(G,o),\dots,B_\ell(G,o)).
\]
Its state space is $\cS_\ell:=\{(\tau_1B,\dots,\tau_\ell B):B\in\ballsp_\ell\}.$
The map
\[
    \mvgt_\ell(B):=(\tau_1B,\dots,\tau_\ell B)
\]
is a bijection from $\ballsp_\ell$ to $\cS_\ell$, with inverse given by the last coordinate.

For $B=(V,E,o)\in\ballsp_\ell$ and $u\in V$, perform one duplication
step inside $B$. Take independent variables
\[
 \chi\sim\ber(r),
 \qquad \eta_z\sim\ber(p)\text{ for } z\in\nb_B(u),
\]
add a new vertex $*$, join it to $u$ when $\chi=1$, and join it to each $z\in\nb_B(u)$ when $\eta_z=1$.  Let $\mathsf G_u^{(\ell)}(B)$ be the induced radius-$\ell$ ball about the old root $o$ after this update.  When the parent is the root, let $\mathsf R_\ell(B)$ instead be the induced radius-$\ell$ ball about the new vertex $*$.  The randomness in these maps is fresh at every use.  For bounded $F\colon\ballsp_\ell\to\dR$, define
\begin{align}\label{def:l-ball-generator}
 (\sL_\ell F)(B)
 :=\E\bigr(F(\mathsf R_\ell(B))-F(B)\bigl)
 +\sum_{u\in V}\E\bigr(F(\mathsf G_u^{(\ell)}(B))-F(B)\bigl).
\end{align}
Write $(\bB_s^{(\ell)})_{s\geq0}$ for the corresponding continuous-time Markov chain.  Thus, the first term in \eqref{def:l-ball-generator} is a rate-one reset to a newborn root, while each vertex visible in the ball acts as a potential parent at rate one.

\begin{lem}[$\ell$-ball recurrence]\label{lem:l-ball-closure}
The maps in~\eqref{def:l-ball-generator} depend only on the rooted ball $B$. Moreover, for every bounded $F$,
\begin{align}\label{eq:k-ball-euler}
 \E\bigl(\gn_{t+1}^{(\ell)}F\mid\cF_t\bigr) =\gn_t^{(\ell)} \left(I+\frac1{t+1}\sL_\ell\right) F.
\end{align}
\end{lem}

\begin{proof}
Adding the newborn cannot shorten distances between old vertices. Any path segment through it can be replaced by a path of no greater length through its parent, since its old neighbors lie in the parent's closed neighborhood. Thus the old vertices of a rooted ball do not change. Suppose a potential parent $u$ is at distance $j\leq\ell$ from the old root.   If $j<\ell$, every neighbor of $u$ already belongs to the ball. If $j=\ell$, a neighbor outside the ball is at distance $\ell+1$, and duplicating an edge to it cannot bring it into the updated radius-$\ell$ ball. Thus every copied edge that can affect the old rooted ball is visible in $B$.  For a newborn root, a path of length at most $\ell$ starts through the parent or one of its neighbors. Every old vertex on such a path is within distance $\ell$ of the parent.  Hence the reset map is also determined by $B$.

Conditional on a parent $u$, the newborn-root contribution is
\[
 F\bigl(\mathsf R_\ell(\ball_\ell(G_t,u))\bigr).
\]
An old rooted ball can change only if its root lies in $\ball_\ell(G_t,u)$.  Averaging over the uniform parent and interchanging the sums over roots and parents gives
\[
 \E\bigl((t+1)\nu_{t+1}^{(\ell)}(F)\mid\cF_t\bigr) = t\nu_{t}^{(\ell)}(F)+\gn_t^{(\ell)} (F+\sL_\ell F).
\]
Dividing by $t+1$ proves~\eqref{eq:k-ball-euler}.
\end{proof}

Assume for now that $0<p<p_c$ and $0<r\le 1$. For $B\in\ballsp_\ell$, write
\[
    N(B):=\abs{V(B)},    \qquad d(B):=\deg_B(o).
\]
For $j\le\ell$, write $\tau_jB:=B_j(B,o)$. Define, on states with at most $t$ vertices,
\[
    K_t^{(\ell)}:=I+\frac{\sL_\ell}{t+1},\qquad\cK_{m,T}^{(\ell)}:=K_m^{(\ell)}\cdots K_{T-1}^{(\ell)}.
\]
Since the total update rate is $N(B)+1$, $K_t^{(\ell)}$ is stochastic on these states. Taking expectations in the $\ell$-ball recurrence gives 
\begin{align}\label{eq:l-ball-mean-recursion}
    \mu_{t+1}^{(\ell)}=\mu_t^{(\ell)}K_t^{(\ell)},\qquad\mu_T^{(\ell)}=\mu_m^{(\ell)}\cK_{m,T}^{(\ell)}.
\end{align}

\begin{lem}[Finite-time approximation]\label{lem:l-ball-euler}
The chain with generator $\sL_\ell$ is nonexplosive. Write $(\cP_s^{(\ell)})_{s\ge0}$ for its semigroup. For every fixed $s>0$, with $m_T=\lfloor Te^{-s}\rfloor$, and every finite $A\subset\ballsp_\ell$, 
\begin{align}\label{eq:l-ball-euler-convergence}
    \sup_{B\in A}\norm{\cK_{m_T,T}^{(\ell)}(B,\cdot)-\cP_s^{(\ell)}(B,\cdot)}_{\tv}\longrightarrow0.
\end{align}
The same approximation holds for the degree projection, including when both degree processes are stopped on first hitting $0$.
\end{lem}
\begin{proof}
    We apply Lemma~\ref{lem:euler-approximation} with $Q=\sL_\ell$ and $h_t=1/(t+1)$. Take $V(B):=1+N(B)$. The sublevel sets of $V$ are finite. Every update adds at most one vertex before truncation, and the total update rate is $N(B)+1$. Therefore,
    \[
    \sL_\ell V(B)\le N(B)+1=V(B),
    \]
    which verifies \ref{E1}. Fix a finite $A\subset\ballsp_\ell$, and choose $M$ such that $N(B)\le M$ for every $B\in A$. If the Euler chain starts from $A$ at time $m\ge M$, then at time $t\ge m$, $N(B_t)\le M+t-m\le t$. Hence, its total jump rate satisfies
    \[
    h_t\lambda_{B_t}\le \frac{N(B_t)+1}{t+1}\le1.
    \]
    Thus, \ref{E2} holds. And hence the nonexplosion argument in Lemma~\ref{lem:euler-approximation} applies. Taking $m_T=\floor{Te^{-s}}$, gives the result. 
    Moreover, the same argument holds for degree projection, including when both degree processes are stopped on first hitting $0$.
\end{proof}

\subsection{Proof of Theorem~\ref{thm:l-ball-limit}}\label{ssec:pf-l-ball-limit}
Fix $\ell\ge1$, $0<p<p_c$, and $0<r\le1$. Lemma~\ref{lem:l-ball-euler} proves \ref{E1}, \ref{E2}. Lemma~\ref{lem:l-ball-closure} proves \ref{C1}. Lemmas~\ref{lem:ball-uniform-tightness}, \ref{lem:ball-finite-influence}, and~\ref{lem:ball-regeneration} prove \ref{C2}, \ref{C3}, and~\ref{C4} respectively. Finally, \ref{C5} holds because growth cannot decrease the size of the ball. Hence, Proposition~\ref{prop:empirical-convergence} finishes the proof.

Finally, let $1\le j\le\ell$. For every bounded $F$, we have
\[
    \nu_t^{(j)}F=\nu_t^{(\ell)}(F\circ\tau_j).
\]
Applying the convergence just proved at radii $j$ and $\ell$, the two sides converge in probability to $\Pi_jF$ and $\Pi_\ell(F\circ\tau_j)$, respectively. Since these limits are deterministic, we get
\[
    \Pi_jF=\Pi_\ell(F\circ\tau_j).
\]
Thus, truncating a $\Pi_\ell$-distributed ball at radius $j$ gives $\Pi_j$, proving projective consistency. Similarly, for every bounded $\varphi$, we have
\[
    \nu_t^{(\ell)}(\varphi\circ d)=\nu_t^{\deg}\varphi.
\]
By the convergence just proved and Theorem~\ref{thm:edegree-lln}, the two sides converge in probability to $\Pi_\ell(\varphi\circ d)$ and $\pi\varphi$, respectively. Hence,
\[
    \Pi_\ell(\varphi\circ d)=\pi\varphi.
\]
Therefore, the root degree under $\Pi_\ell$ has distribution $\pi$.\qed

\begin{proof}[Proof of Corollary~\ref{cor:local-weak-limit}]
    By Theorem~\ref{thm:l-ball-limit}, the family $(\Pi_\ell)_{\ell\ge1}$ is projectively consistent. Hence, by the Kolmogorov extension theorem, there exists a random sequence $(B_\ell)_{\ell\ge 1}$ such that $B_\ell\sim \Pi_\ell$, $\tau_jB_{\ell}=B_j$, for $1\le j\le \ell$ almost surely.
    Choose representatives of these balls with a common root, so that $B_\ell$ is the induced radius-$\ell$ ball of $B_{\ell+1}$, and take their union. This gives a rooted connected graph $(G_\infty,o_\infty)$ satisfying
    \[
    B_\ell(G_\infty,o_\infty)=B_\ell,\qquad \ell\ge1.
    \]
    This is locally finite: if a vertex is at distance $j$ from the root, all its neighbors belong to the finite ball $B_{j+1}$. Let $\Pi$ be its law. The events specifying a finite rooted ball generate the $\gs$-algebra of the local topology, so the prescribed ball determines $\Pi$ uniquely.

    We now prove convergence of the empirical measures. We use the usual local metric under two rooted graphs agreeing at radius $\ell$, having distance at most $2^{-\ell}$ (see \cite[Section 1.2]{BS01}). 
    Let $f$ be bounded by one and $1$-Lipschitz, and put $F_\ell(B):=f(B)$, viewing $B$ as a finite rooted connected graph. Then
    \[
    \abs{f(G,o)-F_\ell(B_\ell(G,o))}\le 2^{-\ell}.
    \]
    Consequently,
    \begin{align*}
        \abs{\nu_t^{\mathrm{loc}}f-\Pi f}&\le 2^{1-\ell}+\abs{\nu_t^{(\ell)}F_\ell-\Pi_\ell F_\ell}\le 2^{1-\ell}+2\norm{\nu_t^{(\ell)}-\Pi_\ell}_{\tv}.
    \end{align*}
    This bound is uniform over such $f$. By Theorem~\ref{thm:l-ball-limit}, first letting $t\to\infty$ for fixed $\ell$, and then letting $\ell\to\infty$, prove that $\gn_t^{\loc} $ converges weakly in probability to $\Pi$. Since $U_t$ is conditionally uniform given $G_t$, we have $\E f(G_t,U_t)=\E [\gn_t^{\loc}f].$ Taking expectations in the preceding bound, and using 
    \[
    \E \norm{\gn_t^{(\ell)}-\Pi_\ell}_{\tv} \to 0,
    \]
    which follows from boundedness and convergence in probability, gives $(G_t,U_t)\implies (G_\infty,o_\infty)$ in the local topology. Finally, the root degree is determined by the radius-one ball. Since this ball has law $\Pi_1$, Theorem~\ref{thm:l-ball-limit} shows that $\deg_{G_\infty}(o_\infty)$ has distribution $\pi$.
\end{proof}

\section{Proofs for the direct construction of the local limiting law}\label{sec:limit-law}

\subsection{Proof of Theorem~\ref{thm:spine-local-law}}\label{ssec:pf-spine-local-law}
By the mark construction in~\ref{def:local-lim B}, an edge is present exactly when every mark along its trace succeeds. At each nonterminal step of the trace, replacing the younger endpoint by its parent decreases the genealogical distance by one. The trace of $e=\{x,y\}$ therefore uses $d_{\cT}(x,y)-1$ distinct copying marks and one parent-link mark. Recall that, 
\[
J_{x,y}=\prod_{e\in \cl_\cT(\{x,y\})} M_e
\]
where $M_{\{u,v\}}=P_{\{u,v\}}$ if $d_\cT(u,v)>1$ and $M_{\{u,v\}}=R_{\{u,v\}}$ if $d_\cT(u,v)=1$. For a finite set $F$, all its edges are present exactly when every distinct mark in $\cl_\cT(F)$ succeeds. Conditional independence of these marks proves~\eqref{eq:spine-edge-marginal} and~\eqref{eq:spine-joint-law}.

For local finiteness, let $\cB_k$ be the vertex set of the Yule bush $\mvt_k$ constructed as in~\ref{def:local-lim A}, including $a_k$. These bushes partition $\cT$. Define
\begin{equation}\label{def:Z_k}
 Z_k=\sum_{v\in\cB_k}p^{d_{\cT}(a_k,v)},
 \quad
 Z=\sum_{v\in\cT}p^{d_{\cT}(o,v)}
   =\sum_{k\ge0}p^kZ_k.
\end{equation}

 We use the following standard result.
\begin{lem}\label{lem:yule distance}
    Let $(\cN_t)_{t \geq 0}$ be a Yule process, where we denote the root of the underlying genealogical tree by $o'$, and by $d(\cdot,\cdot)$ the graph distance in it. Define $W_t:=\sum_{u \in \cN_t}p^{d(u,o')}.$ Then we have for any $t \geq 0$, $\E(W_t)=e^{pt}.$
\end{lem}
\begin{proof}
    Define $m(t):=\E (W_t)$, then we have
    \begin{align*}
        m(t+s)-m(t)&=\E\left(\ind\{\text{only one new vertex $\nu$ born in }[t,t+s)\}\cdot p^{d(\nu,o')}\right)+o(s)\\
        &=\E \left(\sum_{u\in \cN_t} p^{1+d(u,o')} \right) \cdot s+o(s)= p\cdot m(t)\cdot s+o(s)
    \end{align*}
    as $s\to 0$. Thus $m(t)$ solves $m'(t)=p\cdot m(t)$, with $m(0)=1$. The result follows.
\end{proof}

Recalling $Z_k$ from~\eqref{def:Z_k}, note that by Lemma~\ref{lem:yule distance}, $\E\bigl(Z\mid(A_j)_{j\ge0}\bigr)  =\sum_{k\ge0}p^ke^{pA_k}<\infty$ almost surely, as $A_k/k\to1$ and $\log p+p<0$. In particular, $Z<\infty$ almost surely. For every vertex $v$, the triangle inequality gives
\[
 \sum_{w\in\cT}p^{d_{\cT}(v,w)} \le p^{-d_{\cT}(v,o)}Z<\infty.
\]
By~\eqref{eq:spine-edge-marginal},
\[
 \E\bigl(\deg_{\cG}(v)\mid\cT\bigr) =\frac rp\left(\sum_{w\in\cT}p^{d_{\cT}(v,w)}-1\right)<\infty.
\]
Countability gives finite degrees simultaneously at every vertex, hence finite balls of every radius.

We now identify the law of the rooted component $(\cG_o,o)$ defined in~\ref{def:local-lim C}. 
At time $-T$, start from the oldest spine vertex born by $-T$ and retain every subsequent birth whose parent is retained, discarding all ordinary descendants already born before the cutoff. Let $H_T$ be the full finite graph on the retained vertices, with adjacencies inherited from $\cG$. During the construction, distinguish the current spine vertex as its root. The points $(-A_k)_{k\ge0}$ form a rate-one Poisson process on $(-\infty,0)$, independent of the ordinary reproduction processes. Its points in $(-T,0]$ are therefore the reset times, and they move the distinguished root to its prescribed newborn spine child. By independent increments, the retained vertices have fresh rate-one ordinary reproduction clocks after the cutoff at $-T$. All duplication and parent-link marks also have their original independent laws. Since a newborn can only attach within its parent's component, ordinary births outside the current root component cannot affect its rooted balls. Thus, for every $\ell$, the radius-$\ell$ ball about the distinguished root follows the reset-and-growth chain with generator $\sL_\ell$, started from a singleton and run for time $T$. At time zero, the distinguished root is $o=a_0$.

The retained genealogical trees increase to $\cT$ and contain the paths between their vertices. The younger-endpoint reduction of every retained pair stays within its genealogical path, so its trace also stays inside the retained tree. Hence the $H_T$ are consistent induced subgraphs of $\cG$. Their union contains every vertex and edge of $\cG$. 

Fix $\ell\ge 1$. Since $B_\ell(\cG_o,o)$ is finite, all of its vertices eventually belong to $H_T$. Induction then gives eventual equality of every fixed-radius rooted ball \ie\,$B_{\ell}(H_T,o)=B_\ell(\cG_o,o)$ for all sufficiently large $T$, almost surely. Together with the established ergodicity of the radius-$\ell$ reset-and-growth chain, this yields
\[
 \operatorname{Law}\bigl(B_\ell(\cG_o,o)\bigr)
 =\lim_{T\to\infty}\operatorname{Law}\bigl(B_\ell(H_T,o)\bigr)
 =\Pi_\ell.
\]
These ball laws determine $\Pi$, and hence $(\cG_o,o)\sim \Pi$. \qed

\begin{proof}[Proof of Proposition~\ref{prop:limit-component-dichotomy}]
Suppose $r<1$. Run reset-and-growth on finite connected rooted graphs, retaining only the root component after each update. This restriction is closed under the dynamics, since a newborn can only attach within its parent's component and distinct components never merge. An ordinary birth adds at most one vertex; a reset either adds one vertex or leaves a singleton. At size $n$, the update rate is at most $n+1$, so the chain is nonexplosive.

Its root degree follows the positive recurrent degree chain, and degree zero is equivalent to the singleton. The singleton is left at rate $2r$ to a rooted edge. Starting from this edge, the first hitting time of the singleton is exactly the first hitting time of zero by the root degree, which has finite mean by positive recurrence of the degree chain. Hence the full component chain has a positive recurrent singleton state and a stationary law supported on finite graphs. Its radius projections have generators $\sL_\ell$, by the same closure argument as in Lemma~\ref{lem:l-ball-closure}, so their stationary laws are $\Pi_\ell$. These identify the full stationary law as $\Pi$.

When $r=1$, every genealogical edge is present, so $\cG_o=\cG$ and $\cT$ is a spanning tree. Recall that $\cT$ consists of an infinite ancestral ray with finite side bushes. Fix a finite vertex set $S$. Beyond the finitely many bushes meeting $S$, the ancestral ray and its attached bushes remain connected and contain all but finitely many vertices of $\cT$. Hence $\cT\setminus S$ has exactly one infinite component, with only finitely many vertices outside it. The graph $\cG\setminus S$ is obtained by adding edges to $\cT\setminus S$. These edges can only merge its components, and the finitely many vertices outside the infinite tree component cannot form another infinite component. Thus $\cG\setminus S$ also has exactly one infinite component. Since $S$ was arbitrary, $\cG$ has exactly one end.
\end{proof}

\begin{proof}[Proof of Proposition~\ref{prop:component-exp-law}]
Write $K_1$ for the singleton. For a finite connected unrooted graph $H$,
define
\[
 f_H(u)=(1-r)(1-p)^{d_H(u)}, \quad a(H)=\sum_{u\in V(H)}\bigl(1-f_H(u)\bigr).
\]
Here $f_H(u)$ is the probability that a newborn with parent $u$ is unattached, and $a(H)$ is the total rate at which growth updates add a vertex to the component. Let $b(H,K)$ be the sum over parents $u\in V(H)$ of the probabilities that a copying birth attaches and produces a graph isomorphic to $K$. Thus $b(H,K)=0$ unless $|V(K)|=|V(H)|+1$, and $\sum_K b(H,K)=a(H)$.
Let $L$ be the full component generator, and write
\[
 S_F(H)=\sum_{o\in V(H)}F(H,o),\quad n=|V(H)|.
\]
Averaging over the root gives
\[
 \frac1n\sum_{o\in V(H)}LF(H,o)
 =\frac1n\left(
   \sum_K b(H,K)S_F(K)    -(1+a(H))S_F(H)    +\sum_{u\in V(H)}f_H(u)F(K_1)
 \right).
\]
Indeed, suppose a duplication at $u$ attaches and produces $K$. The growth updates retain each of the old vertices as a possible root, while the reset with root $u$ gives the newborn as root. Summing these contributions gives $S_F(K)$. If the newborn is unattached, a growth update leaves the component unchanged, while a reset gives the singleton.

We claim that, if the root is initially uniform conditional on the unrooted graph, this remains true at every time. The displayed identity depends on $F$ only through its averages over the possible roots of each graph. First, kill the process when its size exceeds $N$. The state space is then finite, and the forward equations show that this property is preserved. Letting $N\to\infty$, nonexplosion gives the same conclusion for the full process. Starting from the singleton and passing to the stationary law from Proposition~\ref{prop:limit-component-dichotomy}, we conclude that the root under $\Pi$ is uniform conditional on the unrooted graph.

For each unrooted shape $H$, define
\[
 \mu_H=\pr_\Pi(G_\infty\simeq H),  \quad  c_H={\mu_H}/{|V(H)|}.
\]
Each shape indicator has a bounded generator: growth connects only adjacent size levels, and the reset rate is at most one. Stationary balance therefore gives, without assuming a finite mean component size,
\begin{align}\label{eq:component-triangular-law}
 (1+a(H))c_H  = q\mathbf1_{\{H=K_1\}}   +\sum_{|V(K)|=|V(H)|-1}c_Kb(K,H),
\end{align}
where
\[
 q=\E_\pi\bigl((1-r)(1-p)^D\bigr)>0,  \quad D\sim\pi.
\]
Here $q$ is the total stationary rate of resets producing a singleton, including resets from the singleton which leave its shape unchanged. The singleton case gives $q=(1+r)\pi_0$, since $a(K_1)=r$ and $c_{K_1}=\pi_0$.
Now let $w_H=\pr(C_E\simeq H)$. Integrating the pure-growth forward equations against $e^{-t}$ yields
\[
 (1+a(H))w_H
 =\mathbf1_{\{H=K_1\}}
  +\sum_{|V(K)|=|V(H)|-1}w_Kb(K,H).
\]
Both recursions are triangular in graph order, so $c_H=qw_H$. Consequently,
\[
 1=\sum_H|V(H)|c_H  =q\sum_H|V(H)|w_H=qM.
\]
Since $q>0$, $M=1/q\le\infty$. This proves~\eqref{eq:component-exp-normalizer} and uniform rooting then gives~\eqref{eq:component-exp-law}, for every bounded $F$.
\end{proof}

\section{Fixed vertices and the maximum degree}\label{sec:max-degree}

We prove Theorem~\ref{thm:max-degree} in three stages: convergence of
individual degrees, control of the maximum degree, and uniqueness of the
eventual vertex with maximal degree. Throughout, $0<p<1$ and $0<r\leq1$, and all random
variables are defined on the same probability space carrying the graph process.
Vertex $v\geq2$ is born at time $v$ and chooses a parent $u_v$ uniformly
from $[v-1]$. Conditional on $\cF_{v-1}$ and $u_v$, its initial degree
$d_v:=\deg_v(v)$ has distribution
\[
 d_v\equald\ber(r)+\bin(\deg_{v-1}(u_v),p),
\]
with independent summands. Vertex $1$ has $d_{1}=\deg_1(1)=0$. Define
\[
 R_t:=\frac{\gC(t+p)}{\gC(t)},
 \qquad
 M_v(t):=\frac{\deg_t(v)+r/p}{R_t},
 \qquad t\geq v.
\]
Then (by standard gamma function properties and estimates) $R_{t+1}/R_t=1+p/t$, $R_t$ is increasing, and
$c_pt^p\leq R_t\leq C_pt^p$ for constants $0<c_p\leq C_p<\infty$.
The shift $r/p$ absorbs the contribution of parent links, making $M_v$
a martingale. Since $R_t/t^p\to1$, convergence of $M_v$ gives the
required normalization of the degree.

\begin{lem}[The fixed-vertex martingale]\label{lem:fixed-vertex-martingale}
For every $v$ and $t\geq v$, we have
 $\E(M_v(t+1)\mid\cF_t)=M_v(t)$,
\[
 |M_v(t+1)-M_v(t)|\leq\frac1{R_{t+1}},
 \qquad
 \E\bigl((M_v(t+1)-M_v(t))^2\bigm|\cF_t\bigr)\leq\frac{p M_v(t)}{t R_{t+1}} .
\]
\end{lem}

\begin{proof}
At step $t+1$ the degree of $v$ increases by one if the new vertex chooses $v$ as parent and attaches its parent link (probability $r/t$), or chooses a neighbor of $v$ and copies the edge to $v$ (probability $p\deg_t(v)/t$). These alternatives are exclusive.  Hence, we have $$\E(\deg_{t+1}(v)+r/p\mid\cF_t)=(\deg_t(v)+r/p)(1+p/t),$$ for all $t\ge v$. This is the martingale property, and $\deg_{t+1}(v)-\deg_t(v)\in\{0,1\}$ gives the increment bound. Write $\gD\deg_t(v):=\deg_{t+1}(v)-\deg_t(v)$. Finally $$M_v(t+1)-M_v(t)=R_{t+1}^{-1}\bigl[\gD\deg_t(v) -\E(\gD\deg_t(v)\mid\cF_t)\bigr],$$ whose conditional variance is at most 
$$
R_{t+1}^{-2}\pr(\gD\deg_t(v)=1\mid\cF_t)=R_{t+1}^{-2}(r+p\deg_t(v))/t=pR_tM_v(t)/(tR_{t+1}^2)\leq pM_v(t)/(tR_{t+1}).
$$
This completes the proof.
\end{proof}

\subsection{Proof of Theorem~\ref{thm:max-degree}\textup{(a)}}\label{ssec:pf-max-degree-a}
Fix $v$.  By Lemma~\ref{lem:fixed-vertex-martingale} and $\E M_v(t)=\E M_v(v)\leq(v+r/p)/R_v<\infty$,
\[
 \sum_{t\geq v}\E\bigl(M_v(t+1)-M_v(t)\bigr)^2 \leq\sum_{t\geq v}\frac{p \E M_v(v)}{tR_{t+1}}\leq C\sum_{t\geq1}t^{-1-p}<\infty,
\]
so $M_v$ is bounded in $L^2$ and converges almost surely and in $L^2$ to a limit $\xi_v\geq0$. Since $t^{-p}R_t\to1$, $t^{-p}\deg_t(v)=t^{-p}R_tM_v(t)-t^{-p}r/p\to\xi_v$ in the same senses. Taking a countable intersection, these almost-sure limits exist simultaneously for every vertex. In particular, convergence holds jointly for any finite set of vertices, both almost surely and in $L^2$. The $L^2$ convergence also preserves expectations, so
\[
 \E\xi_1=M_1(1)=\frac{r}{p\gC(1+p)},  \qquad  \E\xi_2=\E M_2(2)  =\frac{r+r/p}{\gC(2+p)}  =\frac{r}{p\gC(1+p)}.
\]
Vertices $1$ and $2$ have the same degree at time $2$, and the graph process is invariant in law under interchanging their labels thereafter.
Thus their limits have the same distribution.

It remains to show that the limit is strictly positive. Fix $t\geq v$.  Conditionally on $\cF_t$, the martingale $(M_v(u))_{u\geq t}$ has mean $M_v(t)$ and, by Lemma~\ref{lem:fixed-vertex-martingale} and the martingale property,
\begin{align*}
 \E(\xi_v^2\mid\cF_t)&=M_v(t)^2+\sum_{u\geq t}\E\bigl((M_v(u+1)-M_v(u))^2\bigm|\cF_t\bigr)\\
 &\leq M_v(t)^2+pM_v(t)\sum_{u\geq t}\frac1{uR_{u+1}}
 \leq M_v(t)^2+CM_v(t)t^{-p}.
\end{align*}
The conditional Paley--Zygmund inequality gives
\[
 \pr(\xi_v>0\mid\cF_t)\geq\frac{M_v(t)^2}{\E(\xi_v^2\mid\cF_t)}
 \geq\frac1{1+C/(t^pM_v(t))}\geq c_0>0,
\]
since 
$
t^pM_v(t)\geq t^p(r/p)/R_t\geq r/(pC_p).
$
By L\'evy's zero--one law, $\pr(\xi_v>0\mid\cF_t)\to\ind_{\{\xi_v>0\}}$ almost surely, so $\ind_{\{\xi_v>0\}}\geq c_0$, \ie\  $\xi_v>0$ almost surely.\qed

\subsection{Control of the maximum degree}
Convergence for each fixed vertex does not by itself give convergence of the maximum, because new vertices arrive at every step. We therefore control late vertices uniformly over their entire future evolution. A late vertex starts at approximately $p$ times its parent's normalized degree, hence at most $p$ times the current maximum (Lemma~\ref{lem:newborn-bound}). Its subsequent upward fluctuations are small relative to the running maximum (Lemma~\ref{lem:post-birth-fluctuations}).

Since the normalization $R_t$ is common to all vertices,
\begin{align}\label{eq:max-degree-identity}
 D_{\max}(t)+\frac rp=R_t \gL(t),
 \qquad
 \gL(t):=\max_{v\leq t}M_v(t),
\end{align}
and $t^{-p}R_t\to1$, so part (b) amounts to showing that $\gL(t)\to\Xi$ almost surely.  Note that $M_v(u)\leq\gL(u)$ for all $u\geq v$.  Fix a rational $\gd\in(0,\min\{p,(1-p)/8\})$ and a rational $\eta>0$. This is needed for taking union bounds over countably many null events later. 
Below, constants denoted by $C$ may change from line to line and depend only on $p,r,\gd$. We use the convention that a maximum over an empty set is zero. All the almost-sure bounds below may be taken to hold on a single event of probability one for every such pair $(\gd,\eta)$, since there are only countably many pairs.

\begin{lem}[Newborn vertices]\label{lem:newborn-bound}
Almost surely there is $V_1<\infty$ such that for all $v\geq V_1$,
\[
 M_v(v)\leq(p+\gd)\gL(v-1)+\frac{C\log v}{R_v}.
\]
\end{lem}

\begin{proof}
Let $d:=\deg_{v-1}(u_v)$, so that $d+r/p=R_{v-1}M_{u_v}(v-1)\leq
R_{v-1}\gL(v-1)$.  The multiplicative Chernoff bound
$$
\pr(\bin(d,p)\ge(1+\theta)pd)\leq e^{-\theta^2pd/3}, 0<\theta\leq1,
$$
with $\theta=\gd/p$, gives 
$$
\pr(\bin(d,p)\ge(p+\gd)d)\leq
e^{-\gd^2d/(3p)}\leq v^{-2}\text{ when } d\geq(6p/\gd^2)\log v,
$$
while
$\bin(d,p)\leq d<(6p/\gd^2)\log v$ otherwise.  Hence, for every value
of $d$, 
$$
\pr\bigl(d_v\ge1+(p+\gd)d+(6p/\gd^2)\log v\bigm|\cF_{v-1}\bigr)\leq v^{-2}.
$$
By the Borel--Cantelli lemma, the complementary inequality holds eventually for all
$v\geq V_1$.  Then
$$
M_v(v)=(d_v+r/p)/R_v\leq(p+\gd)R_{v-1}\gL(v-1)/R_v+(1+r/p+(6p/\gd^2)\log v)/R_v,
$$
and $R_{v-1}\leq R_v$.
\end{proof}

\begin{lem}\label{lem:post-birth-fluctuations}
Almost surely there is $V_2<\infty$ such that for all $v\geq V_2$ and all
$t\geq v$,
\[
 M_v(t)-M_v(v)\leq2\gd \max\Bigl\{\eta,\ \max_{v\leq u<t}\gL(u)\Bigr\}.
\]
\end{lem}

\begin{proof}
Let $Q_v(t):=\sum_{u=v}^{t-1}\E((M_v(u+1)-M_v(u))^2\mid\cF_u)$ be the predictable quadratic variation of $M_v$ after its birth. Since $R_{u+1}/R_u=1+p/u$, we have $p/(uR_{u+1})=1/R_u-1/R_{u+1}$. Thus by Lemma~\ref{lem:fixed-vertex-martingale} and $M_v(u)\leq\gL(u)$, 
\begin{align}\label{eq:post-birth-variation}
 Q_v(t)\leq p\max_{v\leq u<t}\gL(u)\sum_{u\geq v}\frac1{uR_{u+1}}
 \leq C_1v^{-p}\max_{v\leq u<t}\gL(u),
 \qquad C_1:=p\sup_{v\geq1}v^p\sum_{u\geq v}\frac1{uR_{u+1}}<\infty .
\end{align}
For $A>0$, let $G_{v,A}$ be the event that $M_v(t)-M_v(v)\geq\gd A$ and $Q_v(t)\leq2C_1v^{-p}A$ for some $t\geq v$.  The increments of $M_v$ are bounded by $1/R_{v+1}\leq Cv^{-p}$, so Freedman's inequality
\cite[Theorem~1.6]{Freedman1975} gives
\[
 \pr(G_{v,A})\leq\exp\Bigl(-\frac{\gd^2A^2}{2\bigl(2C_1v^{-p}A+Cv^{-p}\gd A\bigr)}\Bigr)
 =\exp\bigl(-c \gd^2Av^{p}\bigr),
 \qquad c:=\frac{1}{4C_1+2C\gd}.
\]
For fixed $\eta>0$, the exponential bound is summable over $v\geq1$ and the dyadic values $A=\eta2^j$, $j\geq0$. Indeed, $2^jv^p\geq(2^j+v^p)/2$ separates the two sums. Thus, by the Borel--Cantelli lemma there is $V_2$ such that $G_{v,\eta2^j}$ fails for all $v\geq V_2$ and all $j\geq0$.  Now let $v\geq V_2$, $t\geq v$, and $\wh\gL:=\max\{\eta,\max_{v\leq u<t}\gL(u)\}$. Choose $j\geq0$ with $A:=\eta2^j\in[\wh\gL,2\wh\gL)$.  By~\eqref{eq:post-birth-variation},
$$
Q_v(t)\leq C_1v^{-p}\wh\gL\leq2C_1v^{-p}A,
$$
so the failure of $G_{v,A}$ forces $M_v(t)-M_v(v)<\gd A\leq2\gd\wh\gL$.
\end{proof}

Combining the two lemmas with 
$$
\gL(v-1)\leq\max\{\eta,\max_{v-1\leq u<t}\gL(u)\},
$$
we obtain that almost surely, with $V:=\max\{V_1,V_2\}$, for all $v\geq V$ and $t\geq v$,
\begin{align}\label{eq:late-vertex-bound}
 M_v(t)\leq(p+3\gd)\max\Bigl\{\eta,\max_{v-1\leq u<t}\gL(u)\Bigr\}+\frac{C\log v}{R_v}.
\end{align}

\subsection{Proof of Theorem~\ref{thm:max-degree}\textup{(b)}}\label{ssec:pf-max-degree-b}
\emph{Step 1: the maximum is finite.}  Let
$B$ be the random variable $\max_{w\leq V}\sup_{t\geq w}M_w(t)$, which is finite because each $M_w$ converges. Let 
\[
C':=\sup_{v\geq1}C\log v/R_v<\infty.
\]
Define
\[
L:=\max\{B,\eta,2C'/(1-p-3\gd)\}.
\]
Suppose $\gL(t)>L$ for some $t$, and let $t$ be the first such time, with $w$ a vertex attaining $M_w(t)=\gL(t)$.  Since $\gL(t)>B$, $w>V$, and \eqref{eq:late-vertex-bound} together with $\max_{u<t}\gL(u)\leq L$ gives $\gL(t)\leq(p+3\gd)L+C'\leq L$, a contradiction.  Hence $\sup_t\gL(t)\leq L<\infty$ almost surely.

\emph{Step 2: the limit.}  Since $M_v(t)\to\xi_v$ for every $v$, $\liminf_t\gL(t)\geq\sup_v\xi_v=\Xi$.  Let $\bar L:=\limsup_t\gL(t)\leq L$. For $T\geq V$,~\eqref{eq:late-vertex-bound} bounds every vertex $v\geq T$ by 
$$
\sup_{t\geq v}M_v(t)\leq(p+3\gd)\max\{\eta,\sup_{u\geq T-1}\gL(u)\}+C\sup_{v\geq T}\log v/R_v.
$$
Hence
\[
 \bar L\leq\max\Bigl\{\max_{v<T}\xi_v,\ (p+3\gd)\max\{\eta,\sup_{u\geq T-1}\gL(u)\} +C\sup_{v\geq T}\frac{\log v}{R_v}\Bigr\}.
\]
Letting $T\to\infty$ gives $\bar L\leq\max\{\Xi,(p+3\gd)\max\{\eta,\bar L\}\}$. If $\bar L>\Xi$, then $\bar L\leq(p+3\gd)\max\{\eta,\bar L\}$, which forces $\bar L\leq(p+3\gd)\eta$ because $(p+3\gd)\bar L<\bar L$.  Thus $\bar L\leq\max\{\Xi,(p+3\gd)\eta\}$ almost surely for every rational $\eta>0$, so $\bar L\leq\Xi$, and $\gL(t)\to\Xi$.  Moreover $\Xi\leq L<\infty$ and $\Xi\geq\xi_1>0$.  By \eqref{eq:max-degree-identity}, $t^{-p}D_{\max}(t)\to\Xi$.

\emph{Step 3: late vertices stay below the largest degrees.}
First fix $0<\eps<1-p$, so that $(p+\eps)\Xi<\Xi$. Choose a rational $\gd$ satisfying the earlier restriction and $\gd\leq\eps/8$, and a rational $0<\eta<\Xi$.  Choose $T\geq V$ so large that $\sup_{u\geq T-1}\gL(u)\leq(1+\gd)\Xi$ and $C\sup_{v\geq T}\log v/R_v\leq\gd\Xi$. Then~\eqref{eq:late-vertex-bound} gives, for every $v\geq T$ and $t\geq v$, 
$$
M_v(t)\leq(p+3\gd)(1+\gd)\Xi+\gd\Xi\leq(p+\eps)\Xi.
$$
In particular $\xi_v\leq(p+\eps)\Xi<\Xi$ for $v\geq T$, so $\{v:\xi_v>(p+\eps)\Xi\}\subseteq[1,T)$ is finite and $\Xi=\max_{v<T}\xi_v$. For $\eps\geq1-p$ the stated set is empty, which proves its finiteness for every $\eps>0$.  Since $\gL(t)\to\Xi>(p+\eps)\Xi$, for all large $t$ the vertex attaining $\gL(t)$, i.e.\ the maximum degree, lies in $[1,T)$.  Finally, if the $k$-th largest value $\xi_{(k)}$ of $(\xi_v)$ exceeds $p\Xi$, choose $\eps$ with $(p+\eps)\Xi<\xi_{(k)}$: for large $t$ all vertices $v\geq T$ have $M_v(t)\leq(p+\eps)\Xi$, while the $k$ vertices in $[1,T)$ with the largest limits have $M_v(t)\to\xi_{(1)},\ldots,\xi_{(k)}>(p+\eps)\Xi$. Hence the $k$-th largest of $(M_v(t))_{v\leq t}$ converges to $\xi_{(k)}$, and the claim follows from the corresponding identity $D_{(k)}(t)+r/p=R_tM_{(k)}(t)$, where $D_{(k)}(t)$ and $M_{(k)}(t)$ denote the $k$-th largest degree and normalized degree, respectively.\qed

\begin{rem}
The threshold $p\Xi$ comes from the factor $p$ in the birth degree of a copy. The proof establishes the precise uniform bound
\[
 \lim_{T\to\infty}\sup_{v\geq T}\sup_{t\geq v}M_v(t)\leq p\Xi\text{ almost surely}.
\]
Consequently, every limiting degree strictly above $p\Xi$ is separated from all sufficiently late vertices. This separation is what allows us to pass from fixed-vertex convergence to convergence of the corresponding order statistics and, in part (c), to persistence of their vertex labels. The argument does not establish the same conclusion for order statistics at or below $p\Xi$.

When $p<p_c$, this result is compatible with the power-law tail of the stationary degree distribution $\pi$: that tail concerns the limiting law of a uniformly sampled vertex, whereas the maximum degree grows on the scale $t^p$. In particular, $D_{\max}(t)=(\Xi+o(1))t^p$ almost surely.
\end{rem}

\begin{lem}[Limits along the children of a fixed vertex]\label{lem:child-limits}
For every fixed $w$, there are almost surely infinitely many vertices $v\in\dN$ such that $u_v=w$, and  
$
 \xi_v\longrightarrow p\xi_w,
$
as $v\to\infty$.
\end{lem}
\begin{proof}
Uniform parent choices and the conditional Borel--Cantelli lemma give infinitely many children. The conditional variance estimate used in part~(a) gives
\[
 \E\bigl((\xi_v-M_v(v))^2\mid\cF_v\bigr)
 \le C v^{-p}M_v(v).
\]
Given $\cF_{v-1}$ and $u_v=w$, the binomial birth rule gives
\[
 \var(M_v(v)\mid\cF_{v-1},u_v=w)
 \le C v^{-p}\bigl(M_w(v-1)+v^{-p}\bigr),
\]
and
\[
 \abs{\E(M_v(v)\mid\cF_{v-1},u_v=w)-pM_w(v-1)}
 \le C\bigl(v^{-p}+v^{-1}M_w(v-1)\bigr).
\]
Here we used $R_v/R_{v-1}=1+p/(v-1)$. Since $M_w$ is bounded in $L^2$, these estimates and $\pr(u_v=w\mid\cF_{v-1})=(v-1)^{-1}$ imply
\[
 \E\left(\ind_{\{u_v=w\}}(\xi_v-pM_w(v-1))^2\right)
 \le C_w v^{-1-p}.
\]
The sum over $v$ is finite. Consequently, the errors tend to zero almost surely along the children, and $M_w(v-1)\to\xi_w$ completes the proof. A countable intersection makes the assertion simultaneous for every fixed label $w$.
\end{proof}

\subsection{Uniqueness and persistence of the vertex with largest degree} We now prove part (c). The key is to vary one copying mark while coupling all other randomness. When the parent is a neighbor of $v$, changing that mark from zero to one adds an edge at $v$. Its descendants contribute a strictly positive amount to $\xi_v$, while leaving the limits of all other vertices already present unchanged. Infinitely many opportunities to make this comparison, together with L\'evy's zero--one law, rule out both atoms in an individual limit and ties between two limits.

\begin{lem}\label{lem:extra-edge}
Fix a time $s\geq2$, a vertex $v<s$, and two graphs $G^{(0)}$ and $G^{(1)}$ on $[s]$ that differ only in the edge $\{v,s\}$, present in $G^{(1)}$ and absent in $G^{(0)}$.  The two duplication processes started from $G^{(0)}$ and $G^{(1)}$ at time $s$ can be coupled so that, for all $t\geq s$, $G^{(0)}_t\subseteq G^{(1)}_t,$ and 
 \[
 \deg^{(1)}_t(v)=\deg^{(0)}_t(v)+F_t,
 \qquad
 \deg^{(1)}_t(w)=\deg^{(0)}_t(w)\quad\text{for all }w<s,\ w\neq v,
\]
where $(F_t)_{t\geq s}$ is a Markov chain with $F_s=1$ which, given $F_t$, increases by one with probability $pF_t/t$ and is otherwise unchanged, and $F_t/R_t\to\vartheta$ almost surely with $\vartheta>0$ almost surely.
\end{lem}

\begin{proof}
Run both processes with the same parent choices $u_t$, the same parent-link indicators, and the same copying marks for the edges present in both graphs, and give the edges of $G^{(1)}_t\setminus G^{(0)}_t$ their own independent copying marks.  Since a new vertex copies each edge of its parent independently, each process is a copy of the model, and by induction $G^{(0)}_t\subseteq G^{(1)}_t$. The edges of $G^{(1)}_t$ not in $G^{(0)}_t$ are the copies of copies of $\{v,s\}$, which we call \emph{extra}.  An extra edge is created only when the parent is an endpoint of an extra edge and the corresponding extra mark succeeds, so every extra edge has the form $\{v,y\}$ or $\{y,y'\}$ with $y,y'\geq s$, and no vertex $w<s$ other than $v$ is ever incident to an extra edge. This gives the statement about $\deg_t(w)$.  Let $F_t$ be the number of extra edges at $v$.  It increases by one at step $t+1$ exactly when the parent is one of the $F_t$ extra neighbors $y$ of $v$ and the extra mark of $\{y,v\}$ succeeds, which has conditional probability $pF_t/t$ given $\cF_t$, since the parent is uniform and the extra marks are fresh.  Thus $N_t:=F_t/R_t$ is a nonnegative martingale with increments bounded in absolute value by $1/R_{t+1}$ and
\[
 \E\bigl((N_{t+1}-N_t)^2\mid\cF_t\bigr)  \leq\frac{pF_t}{tR_{t+1}^2}  \leq\frac{pN_t}{tR_{t+1}}.
\]
Since $\E N_t=N_s=1/R_s$, these increment variances are summable. Hence $N_t$ converges almost surely and in $L^2$ to a limit $\vartheta\geq0$. As in the proof of part (a),
\[
 \E(\vartheta\mid\cF_t)=N_t,
 \qquad
 \E(\vartheta^2\mid\cF_t)\leq N_t^2+CN_tt^{-p}.
\]
Here the lower bound needed to apply the conditional Paley--Zygmund inequality comes from $F_t\geq1$: we have $t^pN_t\geq1/C_p$.
Therefore $\pr(\vartheta>0\mid\cF_t)\geq c_0>0$ for every $t\geq s$, and L\'evy's zero--one law implies $\vartheta>0$ almost surely.
\end{proof}

\subsection{Proof of Theorem~\ref{thm:max-degree}\textup{(c)}}\label{ssec:pf-max-degree-c}
Fix $v\neq w$ and $a\in\dR$, and consider the events $A_1:=\{\xi_v=a\}$ and $A_2:=\{\xi_v=\xi_w\}$.  We show that $\pr(A_i)=0$. The first gives the absence of atoms, the second the distinctness.  Let $s_0>\max\{v,w\}$.  Call a step $s>s_0$ \emph{sensitive} if the parent $u_s$ is a neighbor of $v$.  By part (a), $\deg_t(v)\to\infty$ almost surely. Thus we have $\pr(s+1\text{ sensitive}\mid\cF_s)=\deg_s(v)/s\geq1/s$ for all large $s$, and the conditional Borel--Cantelli lemma gives almost surely infinitely many sensitive steps.

Let $\eta_{s,v}$ denote the copying mark for the edge from the new vertex $s$ to $v$ at a sensitive step. Define $\cH_s:=\cF_{s-1}\vee\gs(u_s, \text{ all marks of step }s\text{ except }\eta_{s,v})$, where we omit nothing when $s$ is not sensitive. Thus $\cH_s$ reveals the entire step except the single mark that decides whether the new edge $\{s,v\}$ is added.  Then $\cF_{s-1}\subseteq\cH_s\subseteq\cF_s$, so $(\cH_s)$ is a filtration with $\bigvee_s\cH_s=\cF_\infty$. L\'evy's zero--one law gives $\pr(A_i\mid\cH_s)\to\ind_{A_i}$ almost surely.  Fix a sensitive step $s$.  Given $\cH_s$, the graph $G_s$ is either $G^{(0)}$ or $G^{(1)}=G^{(0)}\cup\{\{v,s\}\}$ according to the value of the independent mark $\eta_{s,v}\sim\ber(p)$, and by the Markov property
\[
 \pr(A_i\mid\cH_s)=p \pr_{G^{(1)}}(A_i)+(1-p) \pr_{G^{(0)}}(A_i),
\]
where $\pr_{G}$ denotes the law of the process restarted from $G$ at time $s$.  Under the coupling of Lemma~\ref{lem:extra-edge}, $\xi^{(1)}_v=\xi^{(0)}_v+\vartheta$ and $\xi^{(1)}_w=\xi^{(0)}_w$ with $\vartheta>0$, so the events $\{\xi^{(0)}_v=a\}$ and $\{\xi^{(1)}_v=a\}$ are disjoint, and so are $\{\xi^{(0)}_v=\xi^{(0)}_w\}$ and $\{\xi^{(1)}_v=\xi^{(1)}_w\}$.  Hence $\pr_{G^{(0)}}(A_i)+\pr_{G^{(1)}}(A_i)\leq1$, and therefore $\pr(A_i\mid\cH_s)\leq\max\{p,1-p\}<1$ at every sensitive step.  Since there are infinitely many sensitive steps, $\ind_{A_i}\leq\max\{p,1-p\}$ almost surely, i.e.\ $\pr(A_i)=0$. Since $a$ was arbitrary, the law of $\xi_v$ has no atoms. Taking a countable intersection over all pairs $v\neq w$ gives pairwise distinctness simultaneously for all vertex limits.

Finally, choose $0<\eps<1-p$ and the finite cutoff $T$ from part (b). The maximum $\Xi$ is attained on the finite set $[1,T)$, and the values $\xi_v$, $v<T$, are almost surely distinct, so the maximizer $V^*$ is unique. Since $M_{V^*}(t)\to\Xi$ while $\limsup_t\max_{v\neq V^*}M_v(t)\leq\max\{\max_{v<T,v\ne V^*}\xi_v,(p+\eps)\Xi\}<\Xi$, the vertex $V^*$ is the unique vertex of maximum degree for all large $t$. For any fixed $k$ with $\xi_{(k)}>p\Xi$, choose $\eps>0$ with $(p+\eps)\Xi<\xi_{(k)}$. Part (b) confines the top $k$ degrees to a finite set for all sufficiently large times. Within that set, the limiting values are distinct, so convergence fixes their relative order and hence the labels of the top $k$ vertices.\qed

\bibliographystyle{abbrv}
\bibliography{copymodel}

\appendix
\section{Alternative proof of positivity of the tail constant, for $0<p\le 1/2$}\label{appendix: c>0 for p<1/2}
In the case when $0<p\le 1/2$ and $0<r\le 1$, it is much simpler to show that the tail constant is positive than the argument in Section~\ref{sec: proof of C>0 p<pc}. We present it here. Firstly, note that for $0<p\le 1/2$, $\gk$ which is the solution to $p^\gk=1-\gk p$, is at least 1. Thus, $x^\gk$ a convex function. 
Thus, by convexity for every $k\ge 1$, we have 
$
w^k-z^k \ge (w-z)((z+1)^k-z^k)\text{ for } w\ge z\ge 0.
$
Let $d_{z,\gk}:=(z+1)^\gk - z^\gk$. By the uniform-threshold coupling from Section~\ref{ssec:coupling} and the above inequality, we have
\begin{align*}
\E\left(W^\gk - Z^\gk \mid X=x, Z=z\right)&\ge d_{z,\gk} \E\left(W-z\mid X=x,Z=z\right).
\end{align*}

\begin{claim} 
Let $W,Z,X$ be as defined above, and $U\sim \unif(0,1)$. We have
\[
\E (W-z\mid X=x,Z=z) = \E \left(\left(X+r/p\right)U^{-p}-r/p-z\mid X=x,Z=z\right).
\]
\end{claim}
\begin{proof}
    We already proved that given $Z=z$ (and $X=x$), we have $(W \mid Z=z, W\ne z)\equald H_{z+1}$. We also established that $\E H_z = (r+z)/(1-p)$. So, we get 
    \begin{align*}
        \E (W-z \mid X=x, Z=z) &= \frac{b_z-a_z}{1-a_z} \E (H_{z+1}-z)
        = \frac{(b_z-a_z)(1+pz+r)}{(1-a_z)(1-p)}.
    \end{align*}
    Now, let us look at $\E\left(\left(X+r/p\right)U^{-p}-r/p-z\mid X=x,Z=z\right)$. We first note that given $Z=z,$ the expression in the expected value is in $[0,1]$ by definition of $Z.$ Moreover, given $X=x, Z=z$, we have $U\sim\unif(\ga,\gb)$ where 
    \begin{align*}
        \ga:=\left(\frac{x+r/p}{z+1+r/p} \right) ^{1/p} \text{ and } \gb:=\left(\frac{x+r/p}{z+r/p} \right)^{1/p}.
    \end{align*}
    In particular, the conditional expected value becomes
    \begin{align*}
        \left(x+r/p\right) \frac{\gb^{1-p}-\ga^{1-p}}{(\gb-\ga)(1-p)} &-r/p-z\\
        &= \frac{1}{1-p}\cdot \frac{(z+r/p)^{-(1-p)/p} - (z+1+r/p)^{-(1-p)/p} }{\left (z+r/p\right)^{-1/p}- (1+z+r/p)^{-1/p}} - (z+r/p)\\
        &=  \frac{z+r/p}{1-p}\cdot \frac{1-a_z^{1-p}}{1-a_z} - (z+r/p)\\
        &= \frac{z+r/p}{(1-a_z)(1-p)}\cdot \left(p-(a_z^{-p}-1+p)\cdot a_z\right)\\
        &= \frac{(pz+r)-(1+pz+r)\cdot a_z}{(1-a_z)(1-p)}
        = \frac{(b_z-a_z)(1+pz+r)}{(1-a_z)(1-p)}.
    \end{align*}
    Here, we used the fact that
    \[
    a_z= \left(\frac{z+r/p}{1+z+r/p}\right)^{1/p} \text{ and }
    b_z= \frac{pz+r}{pz+1+r}.
    \]
    This proves the claim.
\end{proof}
In particular, we have
\begin{align*}
\E\left(W^\gk \mid X=x,Z=z\right)
&\ge z^\gk+d_{z,\gk} \E((x+r/p)U^{-p}-r/p-z\mid X=x,Z=z)\\
&\ge\E\left(((x+r/p)U^{-p}-r/p)^\gk\mid X=x,Z=z\right).
\end{align*}
The final inequality follows from convexity. For $\gk>1$, if $q\in [z,z+1]$, we have
$
z^\gk + d_{z,\gk} (q-z) \ge q^\gk.
$ 
Taking expectation w.r.t.~$Z$, we obtain
\begin{align*}
    \E\left(W^\gk\mid X=x\right)
    &\ge \E\left(((x+{r}/{p})U^{-p}-{r}/{p})^\gk\right)
    \ge  \E\left(x^\gk U^{-p\gk}\right)
    ={x^\gk}/(1-p\gk)={x^\gk}p^{-\gk}.
\end{align*}
Here we used the fact that  $U^{-p}\ge 1$ a.s.~and $p^\gk=1-p\gk$.
Finally, $\E(X\mid Y=y)=r+py$, and Jensen's inequality gives the strict bound
\[
 \E(W^\gk\mid Y=y)
 \ge (r+py)^\gk p^{-\gk}
 =(y+r/p)^\gk>y^\gk.
\]
On the other hand, $\E((pU^{-p}y)^\gk)=y^\gk$. Proposition~\ref{prop:stationary-goldie-condition} permits integrating the difference under $\pi$, even though the two stationary $\gk$-moments separately diverge. Each conditional difference is strictly positive. Thus,
$
 \E\bigl(W^\gk-(pU^{-p}Y)^\gk\bigr)>0
$
and Goldie's formula \cite[equation~(2.18)]{goldie91} gives $C>0$.

\end{document}